\documentclass[10pt]{amsart}

\usepackage{amsmath, amssymb, amsfonts, graphicx, epsfig}

\swapnumbers

\newcommand{\degrad}[1]{\widetilde{\delta}(#1)}
\newcommand{\HT}[2][\ast]{\mathrm{HT}_{#1}(\mathcal{#2})}

\newcommand{\HTsub}[3][\ast]{\mathrm{HT}_{#1}(\mathcal{#2}_{#3})}
\newcommand{\unHT}[2][\ast]{\mathrm{\widetilde{HT}}_{#1}(#2)}
\newcommand{\CT}[2][\ast]{\mathrm{CT}_{#1}(#2)}

\newcommand{\pairing}[2]{\langle\,#1\,,\,#2\,\rangle}

\newcommand{\mir}[1]{\overline{#1}}
\newcommand{\lk}[1]{\mathcal{#1}}

\newcommand{\F}{\mathbb{F}}
\newcommand{\Z}{\mathbb{Z}}

\newcommand{\R}{\mathbb{R}}

\newcommand{\disc}[2][p]{A_{#1}(#2)}
\newcommand{\area}[1]{[#1]}
\newcommand{\ring}[1]{\mathbb{P}_{#1}}
\newcommand{\field}[1]{\mathbb{F}_{#1}}
\newcommand{\cross}[1]{\mathrm{\textsc{cr}}(#1)}
\newcommand{\regions}[1]{\mathfrak{F}_{L}}
\newcommand{\graph}[1]{\Gamma_{L}}
\newcommand{\circles}[2]{\mathrm{\textsc{Cir}}(#1_{#2})}

\newcommand{\lra}{\longrightarrow}

\newcommand{\pic}[1]{\includegraphics[scale=0.33]{#1}}
\newcommand{\treefig}[2]{\includegraphics[scale=#1]{#2}}

\newtheorem{theorem}{Theorem}[section]
\newtheorem{defn}[theorem]{Definition}
\newtheorem{lemma}[theorem]{Lemma}
\newtheorem{cor}[theorem]{Corollary}
\newtheorem{prop}[theorem]{Proposition}

\newcommand{\POz}{P.\ Ozsv{\'a}th\,}
\newcommand{\ZSz}{Z.\ Szab{\'o}\,}

\title{Totally Twisted Khovanov Homology}
\author{Lawrence P. Roberts}

\begin{document}
\maketitle

\section{Introduction}

\noindent In \cite{Khov} M. Khovanov introduced his well known construction of a homology theory for a link, $\lk{L}$, in $S^{3}$ whose Euler characteristic encodes a version of the Jones polynomial for $L$ (see also \cite{Bar1}, \cite{Viro}). This construction used the exponentially many ways a link diagram can be resolved, in a manner analogous to Kauffman's state summation approach to the Jones polynomial. However, Thisthlethwaite and others had established a more efficient means for computing the Jones polynomial by using the spanning trees of the Tait graphs for the link diagram. Trying to repeat this process for the Khovanov homology lead to papers by A. Champanerkar \& I. Kofman, \cite{Cha2}, and S. Wehrli, \cite{Wehr}. Both papers show that, in principle, the Khovanov homology can be computed from a complex whose generators are the spanning trees for one of the Tait graphs of $L$ by demonstrating that the Khovanov homology deformation retracts to a sub-complex whose generators are in one-to-one correspondence with the spanning trees. However, while these constructions identified the generators of the complex with spanning trees, they only demonstrated the existence of the differential for a spanning tree complex, remaining mute about how to fully and explicitly compute it. In \cite{Cha2} A. Champanerkar \& I. Kofman describe parts of the differential, but not the entire structure. \\
\ \\
\noindent In this paper, we describe a new structure in (reduced, characteristic 2) Khovanov homology which leads to a homology theory the author calls totally twisted Khovanov homology (due to formal analogies with the ``totally twisted'' Heegaard-Floer homology). It arises by deforming the Khovanov differential, and is 
a singly graded theory. It shares many of the properties of Khovanov homologies. In particular, it provides a new invariant homology theory for links. However, there is one notable difference: taken with the correct coefficients, the totally twisted homology deformation retracts to a spanning tree complex with a completely explicit differential. This story, with the results of the authors computerized comparison to Khovanov homology, are presented in the first section of the paper. The details of the proof of invariance follow in the remaining sections. The paper concludes by proving the basic results common to most knot homology theories using the spanning tree perspective. \\
\ \\
\noindent{\bf Acknowledgments:} The author did not discover the idea of ``twisting'' a link homology chain complex to deform it to a spanning tree complex. As he understands it, the idea first emerged in unpublished work of \POz and \ZSz in the context of Heegaard-Floer homology. The author learned the Heegaard-Floer idea from John Baldwin while at the Mathematical Sciences Research Institute for the program on Homology theories of knots and links in the spring of 2010. While at MSRI, he stumbled on to the constructions in this paper while trying to understand what he was being told, completing the proof of invariance in fall of 2010. John Baldwin and Adam Levine have used this idea, in conjunction with a construction of C. Manolescu, to describe Ozsv\'ath and Szab\'o's knot Floer homology using spanning trees of a link diagram, \cite{Bald}. The author would like to thank John Baldwin for those conversations, as well as \POz and \ZSz for the great idea. He would also like to thank Liam Watson, Matt Hedden, and Tom Mark for listening as he worked out some of the details while at MSRI. Finally, the author would like to thank MSRI for the great semester.    
 
\section{The construction}\label{sec:const}

\subsection{Preliminaries}
\noindent Throughout $\lk{L}$ will be an {\em oriented} link in $S^{3}$, equipped with a marked point $p \in \lk{L}$. We will study $\lk{L}$ through a link diagram: a generic projection, $L$, of $\lk{L}$ into $S^{2}$, taking $p$ to a non-crossing point. We will always use Roman letters to denote diagrams. 

\begin{defn} For an oriented link diagram, 
\begin{enumerate}
 \item $\cross{L}$ denotes the set of crossings in $L$,
 \item $n_{\pm}(L)$ is the number of right-handed/left-handed crossings.
 \end{enumerate}
\end{defn} 
\noindent  It will be convenient to let $\Gamma_L$ be the image of the projection of $L$ in $S^{2}$. $\Gamma_{L}$ is the four valent graph found by stripping $L$ of its crossing data.
\begin{defn}
The components of $S^{2} \backslash \Gamma_{L}$  are the {\em faces} of $L$. The set of faces will be denoted $\regions{L}$. Let  $\ring{L} = \Z/2\Z[x_{f}|f\in \regions{L}]$ be the polynomial ring which assigns a formal variable to each $f \in \regions{L}$. The field of fractions of $\ring{L}$ will be denoted $\field{L}$.
\end{defn}
\noindent We will make essential use of the requirement that every coefficient ring have characteristic 2. We will usually label the faces of $L$ with the corresponding formal variable.\\
\ \\
\noindent {\bf Example:} As an illustration, to which we return repeatedly, consider the following diagram for the two component unlink:\\
\begin{center}
\treefig{0.5}{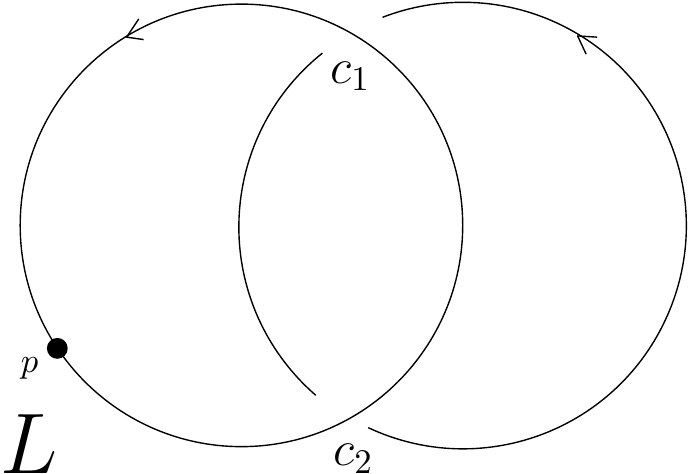} \hspace{1in} \treefig{0.5}{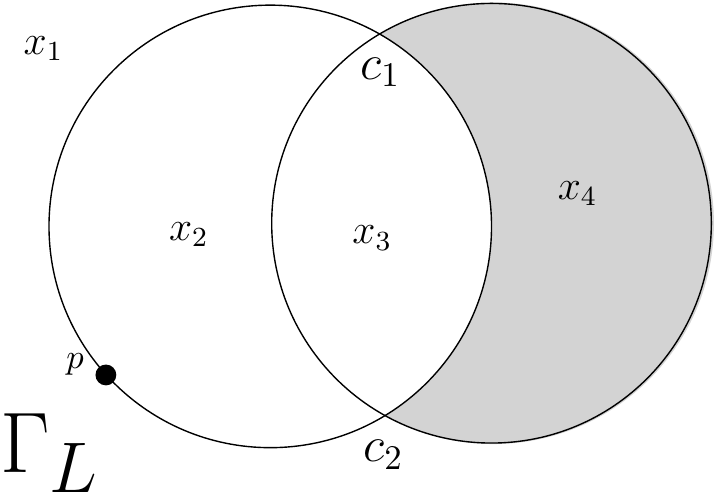}
\end{center}

\noindent On the left is the oriented diagram $L$ with a marked point $p$, and two labeled crossings. We will identify $\cross{L}$ with $\{c_{1}, c_{2}\}$. On the right is the four valent graph $\Gamma_{L}$. The vertices are identified with the crossings, and the marked point is not at a crossing. The shaded region is one face of $L$ labeled with the formal variable $x_{4}$.  There are three other faces in $\regions{L}$, labeled by their formal variables $x_{1}, x_{2}, x_{3}$. Thus, $\ring{L} = \Z/2\Z[x_{1}, x_{2}, x_{3}, x_{4}]$ and $\field{L}$ is the field of binary rational functions in four variables. \\
\ \\   
\subsection{Resolutions}

\begin{center}
\begin{figure}
\treefig{0.5}{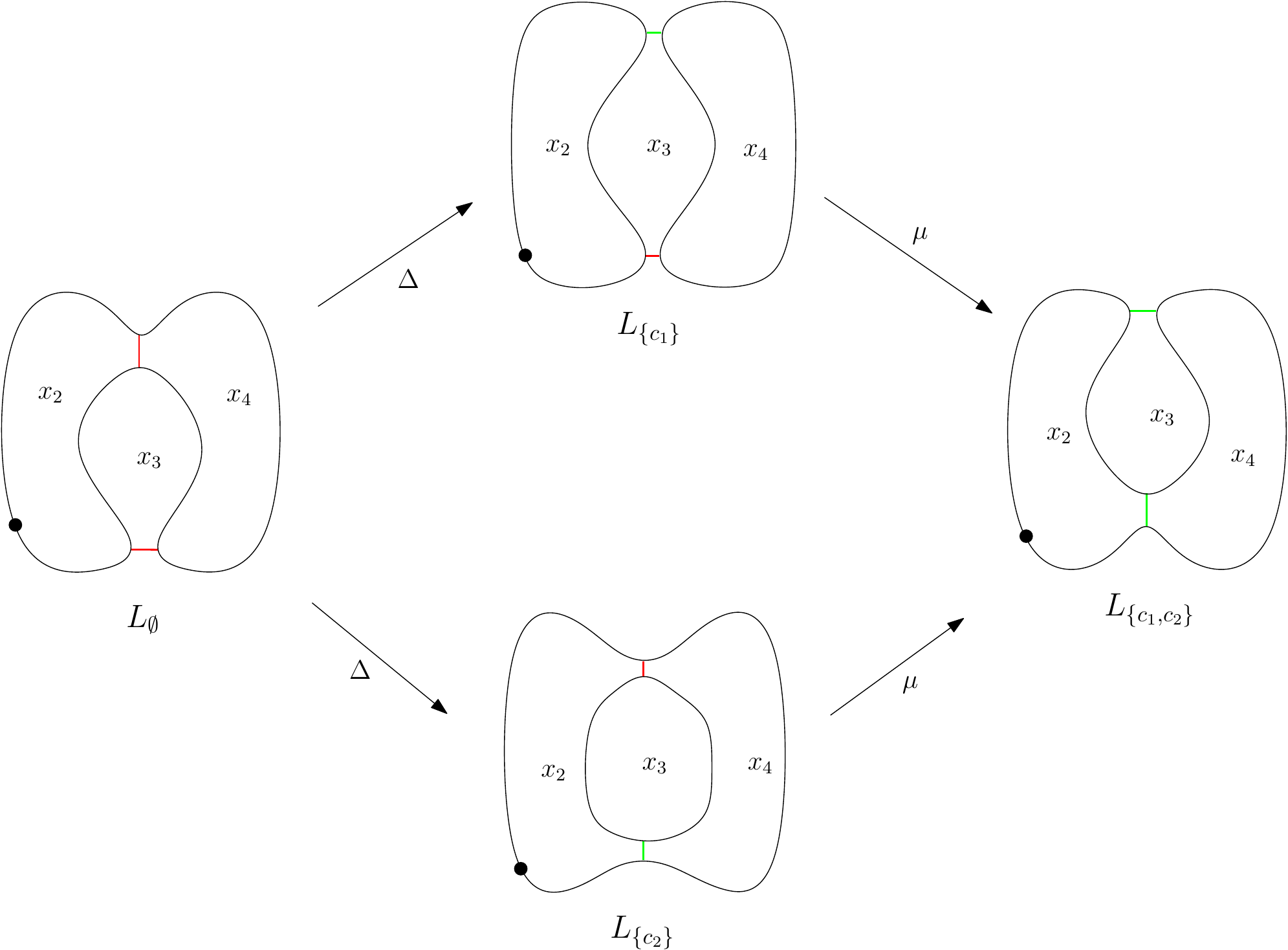}
\caption{The four resolutions for the unlink considered in the text. From left to right we have the resolutions with $\widetilde{\delta}$ equal to $0, 1$ and $2$ respectively. The marked point is depicted with a darkened disc. The labels for each of the faces in each resolution are given as a formal variable in $\ring{L}$. In $L_{\{c_{2}\}}$ the unmarked circle has area $x_{3}$ while the bounded region inside the marked circle has area $x_{2} + x_{3} + x_{4}$. The $\mu$ and $\Delta$ maps indicate whether we use a multiplication or comultiplication when constructing the Khovanov complex. }
\label{fig:resolutions}
\end{figure}
\end{center}

\noindent To each subset $S \subset \cross{L}$, we define the resolution, $L_{S}$, of $L$ to be the {\em diagram} obtained by locally resolving each crossing in $\cross{L}$ according to the rule\footnote{In the more common language of 0 and 1-resolutions used for Khovanov homology, the $0$, $1$-values correspond to the indicator function for $S$ as a subset of $\cross{L}$}:
$$
\pic{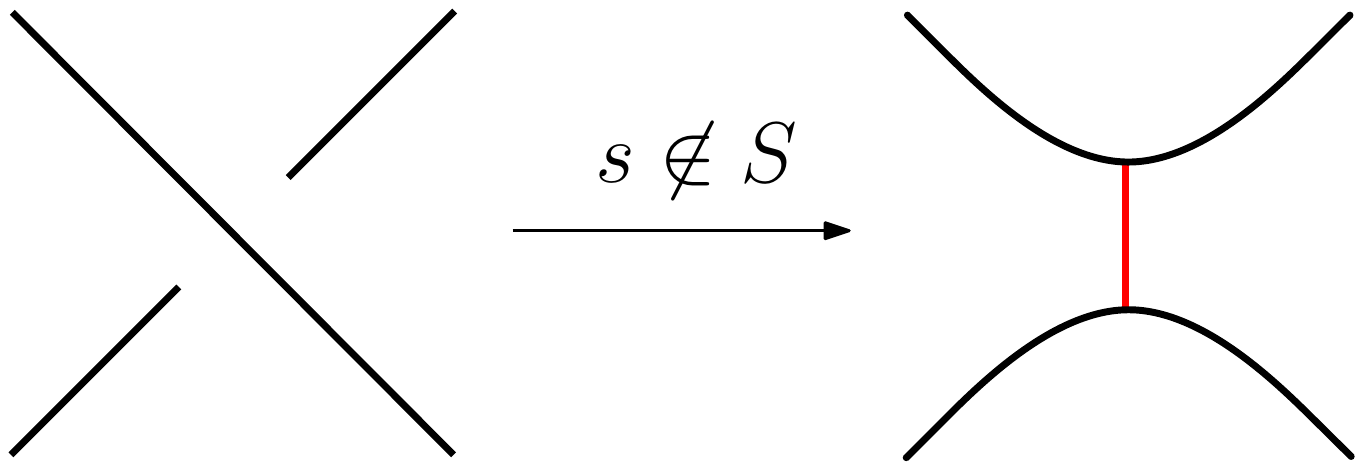} 
\hspace{.5in} 
\pic{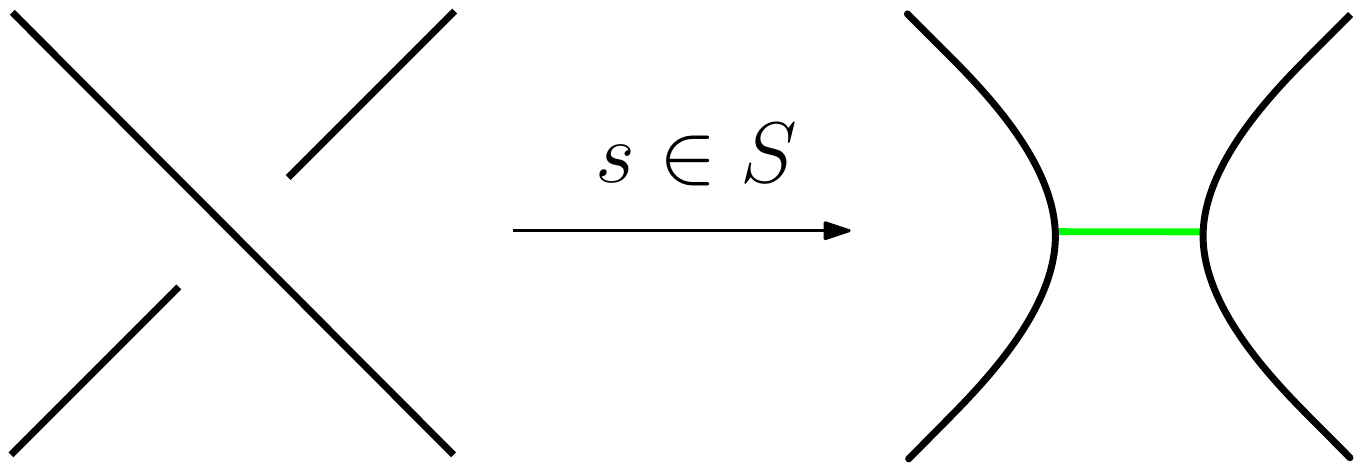} 
$$
The diagram $L_{S}$ consists of disjoint circles embedded in $S^{2}$. One of the circles contains the image of the marked point $p$ and will be called the {\em marked circle} for $L_{S}$. We decorate the resolution $L_{S}$ with colored arcs as in the diagram -- for those crossings $s \not\in S$ we add a red arc while for $s \in S$ we use a green arc. While the colored arcs will be useful below, we do not consider them part of the resolution.\\
\ \\
\noindent We can partially order the resolutions by asserting $L_{S} \leq L_{S'}$ when $S \subset S'$.  
\begin{defn}\ \\[-0.2in]
\begin{enumerate}
\item $R(L)$ is the set of resolutions $L_{S}$, $S \subset \cross{L}$. It is in one to one correspondence with the power set of $\cross{L}$
\item Given a resolutions $L_{S}$, $\degrad{L_{S}} = |S|$
\item For $i \in \Z$, $R_{i}(L)$ is the subset of $R(L)$ consisting of those $L_{S}$ with $\degrad{L_{S}} = i$. 
\end{enumerate}
\end{defn}
\noindent $R_{i}(L)$ consists of those resolution diagrams adorned with exactly $i$ green arcs.\\
\ \\
\noindent {\bf Example(cont.):} The four resolutions for our example unlink are depicted with their colored arcs in Figure \ref{fig:resolutions}.

\subsection{Discs and formal areas}

\noindent We use the additional arcs to divide the regions of $S^{2} \backslash L_{S}$ into a collection of components in one to one correspondence with $\regions{L}$. We call these components the faces of $L_{S}$ and label them with the same formal variable as the corresponding of $L$ it intersects. If $R$ is a region in $S^{2}$ which is a union of faces of $L_{S}$ we can assign it a formal area in $\ring{L}$
$$
[R] = \sum_{\big\{f \in \regions{L} \big| f \cap R \neq \emptyset\big\}} x_{f}
$$
When referring to regions by their numerical indices we will use the following shorthand: $[i_{1}\ldots i_{k}]$ will denote $[\big\{i_{1},\ldots, i_{k}\big\}]$. In particular both equal $x_{i_{1}} + \ldots + x_{i_{k}}$. Thus $[2]$ is $x_{2}$ and will be preferred to $[\{2\}]$.\\
\ \\
\noindent Using this data and the marked point $p$, we can assign a ``formal area'' to each of the circles in a resolution $L_{S}$. Every unmarked circle $C$ in $L_{S}$ bounds two discs in $S^{2}$ which are unions of faces in $L_{S}$. These two discs can be distinguished by which one contains the marked point $p$. 
\begin{defn}
Given $S \subset \cross{L}$, let $\circles{L}{S}$ be the set of circles in the resolution $L_{S}$. Given an unmarked circle $C \in \circles{L}{S}$ let  
\begin{enumerate}
\item $\disc{C}$ be the component of $S^{2}\backslash C$ which does not contain $p$. $\disc{C}$ will be called the {\em interior} of $C$.
\item the {\em exterior} of $C$ in $L_{S}$ is the component of $S^{2}\backslash{C}$ which contains $p$.   
\item $\area{C}$ will denote $[A_{p}(C)]$, the formal area of the interior of $C$, i.e. the sum of the formal variables assigned to the faces contained in $C$.
\end{enumerate}
\end{defn}

\subsection{Complexes associated to resolutions}

\noindent As with reduced Khovanov homology we will associate an algebraic object to each resolution $L_{S}$ of a diagram $L$. In this paper, we assign 
a Koszul complex to $L_{S}$. In our setting we need Koszul complexes over $\ring{L}$ which take the following form.

\begin{center}
\begin{figure}
\treefig{0.75}{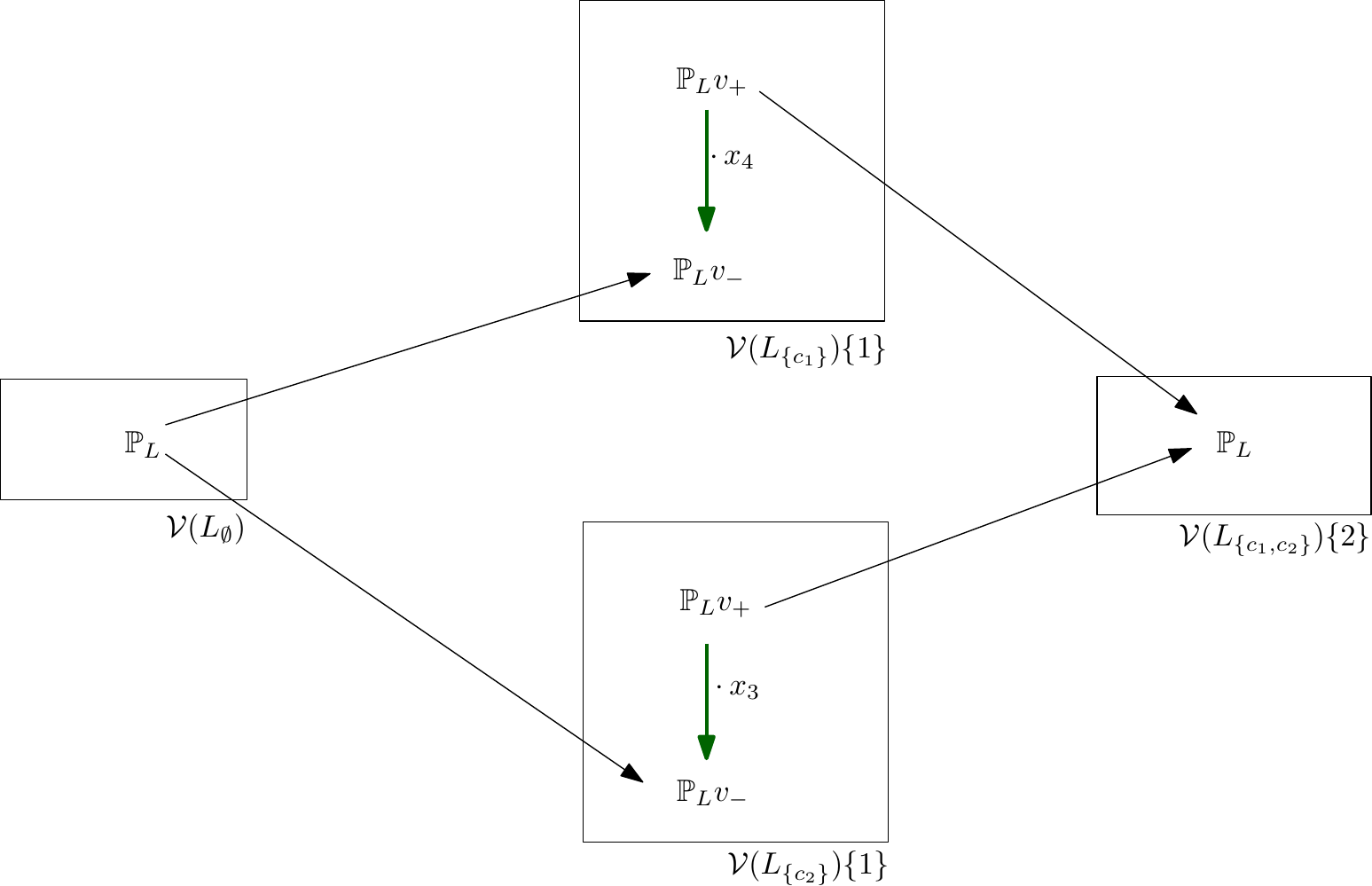}
\caption{In each box is an example complex of the form $\mathcal{V}(L_{S})$ for the resolutions in Figure \ref{fig:resolutions}. The relative positions of the complexes are identical to the relative positions of the resolutions in that figure. The thick, green vertical arrows correspond to the differentials in the Koszul complex for that resolution, using the formal areas of the unmarked circles in that resolutions. Furthermore, we have amassed the complexes for each resolution into an example of a totally twisted Khovanov complex. }
\label{fig:twisted}
\end{figure}
\end{center}

\begin{defn}
Let $R$ be a ring of characteristic 2. For $r \in R$ $\mathcal{K}_{R}(r)$ is the complex
$$
0 \longrightarrow R\,v_{+} \stackrel{\cdot r}{\longrightarrow} R\,v_{-} \longrightarrow 0
$$
where $v_{\pm}$ occur in gradings $\pm 1$. Given elements $r_{1}, \ldots, r_{k} \in \ring{L}$ the Koszul complex $\mathcal{K}_{R}(r_{1}, \ldots, r_{k})$ is the tensor product complex
$$
\mathcal{K}(r_{1}) \otimes_{R} \mathcal{K}(r_{2}) \otimes_{R} \cdots \otimes_{R} \mathcal{K}(r_{k})
$$ 
\end{defn}
\noindent Our definition differs from the normal one in three respects: (1) we only define the complex in characteristic 2 and thus do not specify a sign convention for the differential in the tensor product complex, (2) the gradings are different than usual but will be more convenient for our purposes, and (3) we choose to explicitly label a basis for the rank 1 free modules in $\mathcal{K}(r)$. \\
\ \\
\noindent We will need to shift the gradings arising in the Koszul complexes above. We use the following notation throughout the paper
\begin{defn}
Let $M = \oplus_{\vec{v} \in \Z^{k}} M_{\vec{v}}$ be a $\Z^{k}$-graded $R$-module, then $M\{\vec{w}\}$ is the $\Z$-graded module with $(M\{\vec{w}\})_{\vec{v}} \cong M_{\vec{v}-\vec{w}}$. 
\end{defn}
\noindent Consequently, $\mathrm{grad}(m\{\vec{w}\}) = \mathrm{grad}(m) + \vec{w}$ on homogeneous elements $m \in M$.\\
\ \\
\noindent We will use the ring associated to the faces of $L_{S}$ as well as the marked point to associate a Koszul type complex to the resolution $L_{S}$:
\begin{defn}
Given a link diagram $L$, and a subset $S \subset \cross{L}$, let $\circles{L}{S} = \big\{ C_{0}, C_{1}, \ldots, C_{k} \big\}$ with $C_{0}$ being the marked circle. Define
$$\mathcal{V}(L_{S}) \cong \ring{L}v_{0}\,\otimes_{\ring{L}}\,\mathcal{K}_{\ring{L}}(\area{C_{1}}, \ldots, \area{C_{k}})$$
where $v_{0}$ is in grading $0$. The grading on any $\mathcal{V}(L_{S})$ will be called a $q$-grading. 
\end{defn}
\noindent To summarize, we use the formal areas of the unmarked circles in $L_{S}$ as the sequence of elements of $\ring{L}$ in forming a Koszul complex\footnote{One needs a convention for which circle to choose as $C_{1}$, $C_{2}$, etc. Such a convention can be found in \cite{Bar1} and is the one used for computations described later in this paper}. However, we do not use the rank 2 complex for $C_{0}$. Instead, let $\mathcal{K}_{C_{0}}$ be the trivial complex $0 \rightarrow \ring{L}v_{0} \rightarrow 0$, supported in degree $0$. This is the first factor in the tensor product defining $\mathcal{V}(L_{S})$.   
\begin{defn}
Let $C$ be an unmarked circle. The differential in $\mathcal{K}([C])$ will be denoted $\partial_{C}$. Thus the differential in $\mathcal{V}(L_{S})$ will be  $\partial_{\mathcal{V}(L_{S})} = \sum_{i > 0} \partial_{C_{i}}$. These change the $q$-grading by the $-2$.
\end{defn}
\noindent As usual, we represent the basis for $\mathcal{V}(L_{S})$ through decorations on the diagram $L_{S}$: each basis element can be represented as $L_{S}$ with circle $C_{i}$, $i > 0$, adorned with either a $+$ or a $-$, depending on whether the element $v_{+}$ or $v_{-}$ is in $i^{th}$-factor of the basis element, \cite{Bar1}. We will call these basis elements, and their corresponding diagrams, {\em pure states}, while elements of $\mathcal{V}(L_{S})$ will be called states.\\ 
\ \\
\noindent {\bf Example:} In Figure \ref{fig:resolutions}, $L_{c_{2}}$ consists of two circles: $C_{0}$ and $C_{1}$. Then $\mathcal{K}(C_{1})$ is the complex $\mathcal{K}_{\ring{L}}([3])$ and $\mathcal{V}(L_{S})$ is the complex $\ring{L}\,v_{0}\,\otimes_{\ring{L}}\,\mathcal{K}_{\ring{L}}([3])$. Each of complexes $\mathcal{V}(L_{S})$ associated to the four resolutions in Figure \ref{fig:resolutions} are depicted depicted in Figure \ref{fig:twisted} (with the factor from the marked circle suppressed).

\noindent Before proceeding we note that if we mod out by the ideal $(x_{f} | f \in \regions{L})$, then the complex $\mathcal{V}(L_{S})$ is precisely the module associated to $L_{S}$ by Khovanov in defining the characteristic 2, reduced Khovanov homology, \cite{Khov}, \cite{Bar1}. This motivates the use of $v_{\pm}$ and the definition of the $q$-grading, as these will then be identical to Khovanov's. 

\subsection{Totally twisted Khovanov homology:} We can now describe the {\em totally twisted Khovanov complex}. The reader familiar with reduced Khovanov homology in characteristic 2 will immediately recognize the procedure as that of \cite{Khov}, but with significant additions. \\
\ \\
\noindent If we follow Khovanov's construction, we should combine the complexes $\mathcal{V}(L_{S})$ into a bigraded module $\widetilde{KT}(L)$ by taking $\widetilde{KT}(L) = \oplus K^{i}(L)$ where for each $i \in \Z$  
$$
K^{i}(L) = \bigoplus_{S \subset \cross{L},\ |S| = i}  \mathcal{V}(L_{S})\{|S|\}
$$
is a $q$-graded chain complex. Since $\mathcal{V}(L_{S})$ has a differential, $\partial_{\mathcal{V}(L_{S})}$, $\widetilde{KT}(L)$ inherits a differential, $\partial_{\mathcal{V}}$, which reduces the $q$-grading by 2. The other grading on $\widetilde{KT}(L)$ is $K^{i}(L) \rightarrow i$, and will be called the $h$-grading. The homogeneous elements of $\widetilde{KT}(L)$ in $h$-grading $i$ and $q$-grading $j$ will be denoted $\widetilde{K}^{(i,j)}$.  Shifts in the bigrading will be denoted by $\{(\Delta h, \Delta q)\}$. \\
\ \\ 
\noindent  As a bigraded module $\widetilde{KT}(L)$ is identical to that used for the (unshifted, reduced) Khovanov chain complex (tensored with $\ring{L}$). Khovanov  discovered a $(+1,0)$ differential, $\partial_{KH}$, on this bigraded module, which defines the reduced, unshifted Khovanov complex. Note that $(\widetilde{KT}(L), \partial_{KH})$ is a direct sum of chain complexes, one in each $q$-grading. \\
\ \\
\noindent $\partial_{KH}$ is constructed from the multiplication, $\mu: V \otimes V \rightarrow V$ and the comultiplication $\Delta: V \rightarrow V \otimes V$ for a certain commutative Frobenius algebra, $V = \Z/2\Z\,v_{+} \oplus \Z/2\Z\,v_{-}$, graded as above. On the graded basis, these maps are given by
$$
\mu : \left\{\begin{array}{l} 
	v_{+} \otimes v_{+} \rightarrow v_{+} \\
	v_{+} \otimes v_{-} \rightarrow v_{-} \\
	v_{-} \otimes v_{-} \rightarrow 0 \\
	\end{array}\right.
\hspace{.5in}
\Delta : \left\{\begin{array}{l} 
	v_{+} \rightarrow v_{+} \otimes v_{-} + v_{-} \otimes v_{+} \\
	v_{-} \rightarrow v_{-} \otimes v_{-} \\
	\end{array}\right.
$$
Each map shifts the $q$-grading by $-1$, so to obtain a differential preserving the $q$-grading the image needs to be shifted by $1$. \\
\ \\
\noindent If we tensor $V$ with $\ring{L}$ we obtain the module underlying the complex $\mathcal{K}(\area{C})$ for $C$ an unmarked circles in some resolution $L_{S}$. $\mu$ and $\Delta$ can be extended to maps on the modules $\mathcal{K}(\area{C})$, but when we try to extend to the chain complexes we need to account for the different areas incorporated in each circle. Nevertheless, in section \ref{sec:chainmaps} we prove 

\begin{prop} If $C$ is a circle in $L_{S \cup \{i\}}$ formed by merging the circles $C_{1}$ and $C_{2}$ in $L_{S}$, then
$$
\mu : \mathcal{K}(\area{C_{1}}) \otimes \mathcal{K}(\area{C_{2}})\longrightarrow \mathcal{K}(\area{C})\{1\}
$$
is a chain map on the Koszul complexes over $\ring{L}$. Likewise, if $C_{1}$ and $C_{2}$ arise in $L_{S \cup \{i\}}$ from dividing a circle $C$ in $L_{S}$, then
$$
\Delta : \mathcal{K}(\area{C}) \longrightarrow \big(\mathcal{K}(\area{C_{1}}) \otimes \mathcal{K}(\area{C_{2}})\big)\{1\}
$$
is a chain map over $\ring{L}$. 
\end{prop} 

\noindent On a summand $\mathcal{V}(L_{S})$ with $|S| = i$, $\partial_{KH}$ is a sum of chain maps $\mathcal{V}(L_{S}) \rightarrow \mathcal{V}(L_{S \cup \{i\}})\{1\}$ for each $i \not\in S$. The addition of $i$ corresponds to changing $L_{S}$ at one crossing, and thus either merging two circles in $L_{S}$ -- in which case the map uses $\mu$ on the corresponding factors and the identity on the others -- or we divide a circle, in which case we use $\Delta$ on the factor corresponding to the splitting circle. The remaining factors are mapped by the identity. Thus $\partial_{KH}$ is a chain map when considering unmarked circles.\\
\ \\
\noindent  For the marked circle, we treat $v_{0}$ as a shifted $v_{-}$; so
$$
\mu : \left\{\begin{array}{l} 
	v_{+} \otimes v_{0} \rightarrow v_{0} \\
	v_{-} \otimes v_{0} \rightarrow 0 \\
	\end{array}\right.
\hspace{.5in}
\Delta : v_{0} \rightarrow v_{0} \otimes v_{-} \\
$$
since dividing the marked circle results in a marked circle and a new unmarked circle. These also extend to chain maps on the Koszul complexes above.\\
\ \\
\noindent Together these results imply

\begin{theorem}
Let $\partial_{KH} : \widetilde{KT}^{\ast, \ast}(L) \longrightarrow \widetilde{KT}^{\ast+1,\ast}(L)$ be the Khovanov differential, and let $\partial_{\mathcal{V}} : \widetilde{KT}^{\ast, \ast}(L) \longrightarrow \widetilde{K}^{\ast, \ast - 2}(L)$ be the Koszul differential. Then $\widetilde{\partial} = \partial_{KH} + \partial_{\mathcal{V}}$ is a differential on $\widetilde{KT}(L)$. 
\end{theorem}

\noindent Thus, $(\widetilde{KT}(L), \widetilde{\partial})$ is a chain complex, the {\em unshifted totally twisted Khovanov complex} for $L$. \\
\ \\
\noindent{\bf Example (cont.):} Figure \ref{fig:twisted} depicts the unshifted totally twisted Khovanov complex for the two component unlink we introduced earlier. Each of the horizontal arrows is an isomorphism in this complex. \\
\ \\
\noindent{\it Gradings:} The (unshifted) reduced Khovanov complex over $\ring{L}$ is $(\widetilde{KT}(L), \partial_{KH})$. Since $\partial_{KH}$ is a $(+1,0)$-differential on the bigraded module, the homology of this complex is also bigraded, with the $h$-grading being the homology grading. With the addition of $\partial_{\mathcal{V}}$ the homology is no longer naturally bigraded. Instead, we can equip $\widetilde{KT}(L)$ with a single grading:
\begin{defn}
The $\delta$-grading on the bigraded module $\widetilde{KT}^{\ast,\ast}(L)$ is  $\delta: \widetilde{K}^{i,j} \rightarrow 2i - j$. 
\end{defn}
\noindent $\partial_{KH}$ is a $(+1,0)$ map, and thus changes $\delta$ by $+2$. $\partial_{\mathcal{V}}$ is a $(0,-2)$ map, so it also shifts the $\delta$-grading by $+2$. Thus, $\delta$ provides a grading to the complex $(\widetilde{KT}(L), \widetilde{\partial})$. The $\delta$-grading will be written as a subscript, to distinguish it from the $q$ and $h$-gradings\footnote{This definition is $-2$ times the definition of the $\delta$-grading found in other papers on Khovanov homology}.\\
\ \\
\noindent{\it Shifting:} In addition, the homology of the unshifted complex is not quite an invariant of $\mathcal{L}$. This is also true for the reduced Khovanov homology $(\widetilde{KT}(L), \partial_{KH})$. To define a complex whose homology is an invariant of $\lk{L}$, Khovanov shifts the whole (bigraded) complex $(\widetilde{KT}(L), \partial_{KH})$ by $\big\{\big(-n_{-}(L), n_{+}(L) - 2n_{-}(L)\big)\big\}$. We will make the same shift, which changes $\delta$ by 
$2(-n_{-}(L)) - (n_{+} - 2n_{-})(L)$ $= -n_{+}(L)$. Once we make this shift, the resulting complex will be (almost) an invariant of the link $\lk{L}$.  

\begin{defn} The {\em unshifted totally twisted Khovanov complex} for a link diagram $L$ is the complex $(\widetilde{KT}_{\ast}(L), \widetilde{\partial})$ where
$$
\widetilde{KT}_{\ast}(L) = \bigoplus_{S \subset \cross{L}} \mathcal{V}(L_{S})\{|S|\}[|S|]
$$
and $\widetilde{\partial} = \partial_{KH} + \partial_{\mathcal{V}}$, equipped with the $\delta$-grading. \\
\ \\
\noindent The {\em totally twisted Khovanov complex}, $\underline{KT}_{\ast}(L)$ is $\widetilde{KT}_{\ast}(L)[-n_{+}(L)]$, the complex resulting from shifting the $\delta$-grading by $-n_{+}(L)$. We will denote the homology of this complex by $\underline{HT}_{\ast}(L)$.
\end{defn}

\noindent {\bf Comment:} The name ``totally twisted'' comes from assigning a formal variable to every face in $L$. The construction works equally well if we only assign formal variables to some of the faces of $L$, which is equivalent to modding out by the ideal generated by the remaining faces. Therefore, between the reduced Khovanov complex over $\ring{L}$ and $KT_{\ast}(L)$ there are many twisted complexes, one for each subset of faces, with $KT_{\ast}(L)$ being the most twisted.

\subsection{Invariance}

\noindent In section \ref{sec:invariance} we prove the fundamental result for the invariance of the totally twisted Khovanov homology:

\begin{theorem}
Let $L$ be a digram for $\lk{L}$ with marked point $p$, and let $L'$ be another diagram with the {\em same marked point}. Then $L$ and $L'$ differ by Reidemeister I, II, and III moves conducted away from the marked point. Furthermore, the chain complexes $KT_{\ast}(L)$ and $KT_{\ast}(L')$ are stably chain homotopy equivalent. 
\end{theorem}

\noindent In that section we also explain the notion of stable isomorphism needed in this paper. This notion relates the rings $\ring{L}$ and $\ring{L'}$ even though $\regions{L}$ and $\regions{L'}$ may not be in one-to-one correspondence. \\
\ \\
\noindent It should be noted that the author does not know if a similar result, or just isomorphism of the corresponding homology modules, occurs when changing the marked point. However, we will change coefficients in a moment, and then the invariance under change of marked point can be established. 

\subsection{The spanning tree deformation}

\begin{center}
\begin{figure}
\treefig{0.75}{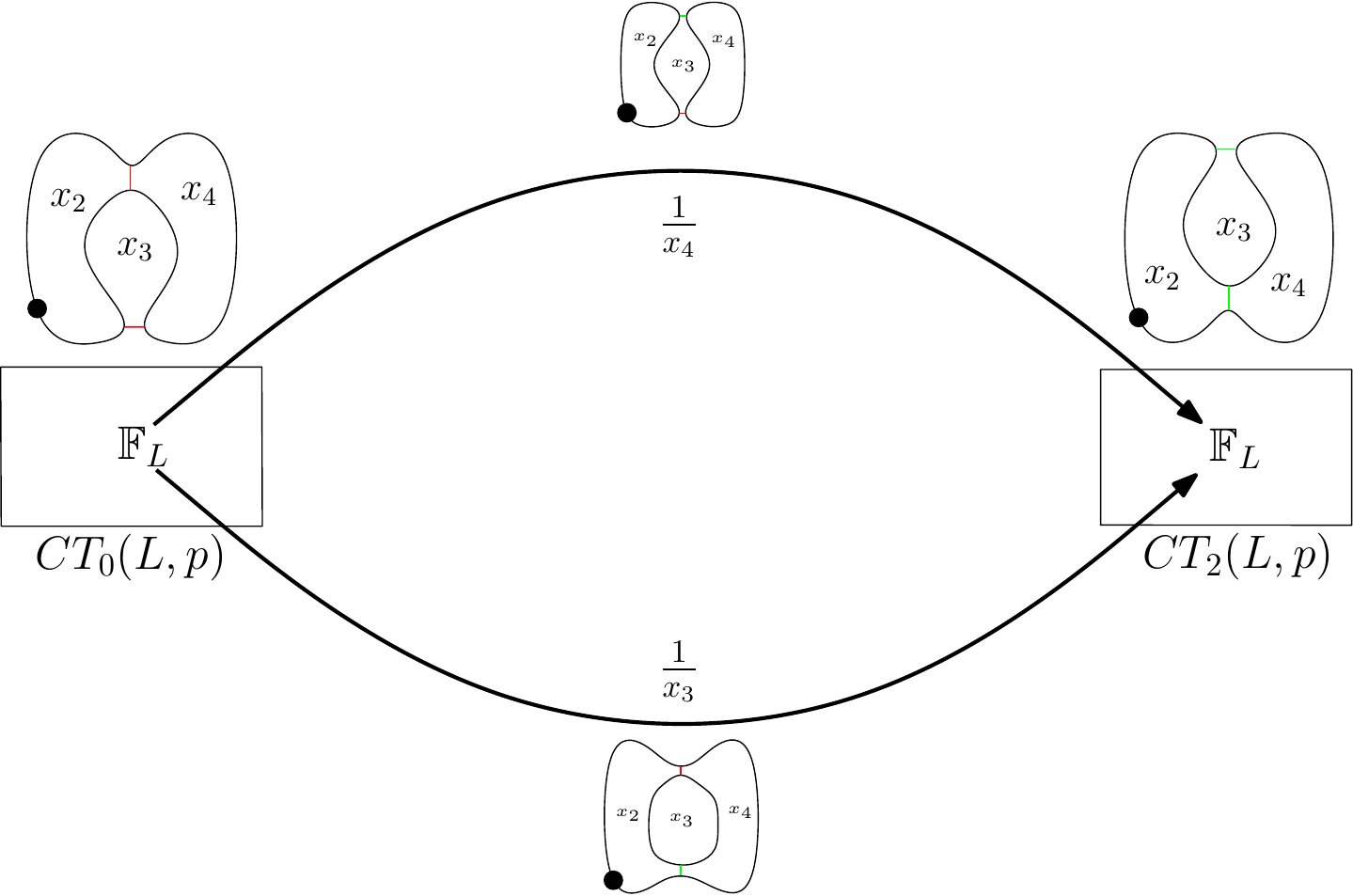}
\caption{The unshifted spanning tree complex for the unlink example in Figure \ref{fig:twisted}. Over $\field{L}$ the vertical arrows in Figure \ref{fig:twisted} are isomorphisms. Simplifying both of these leaves the left and right diagrams, as these consist only of a single marked circle. The expense is a change in coefficient, written next to the top and bottom arrows. The differential in the new complex is the sum of the maps corresponding to the two arrows. The top arrow comes from changing the resolution at $c_{1}$ first, followed by $c_{2}$. This alteration cleaves from the marked circle the region marked with $x_{4}$ and then rejoins it. Call this $B_{T,T'}$, so $[B_{T,T'}] = x_{4}$. This means $W_{T,T'}$ is the region cleaved off by changing the resolution at $c_{2}$ first and then at $c_{1}$. This is the region contained in the inner circle, so $[W_{T,T'}] = x_{3}$. } 
\label{fig:tree}
\end{figure}
\end{center}

\noindent As stated in the introduction, our interest in the totally twisted Khovanov homology comes from the existence of an homotopy equivalent spanning tree model for the chain complex $KT_{\ast}(L)$. However, to realize this model it is necessary to use elements of the field $\field{L}$ as coefficients. Doing so makes each term in $\partial_{\mathcal{V}}$ into a vector space isomorphism. For example, in Figure \ref{fig:twisted} the thickened vertical arrows are isomorphisms over $\field{L}$. We can use these isomorphisms to simplify the chain complex using an analog of Gaussian elimination:\\
\begin{quotation}
For a chain complex, $C$, over a field, if $\partial\,v = \lambda\,w + z$ with $\lambda \neq 0$ (and $w$ and $z$ linearly independent), there is a chain homotopy equivalent complex defined on $C/\mathrm{Span}\{v,w\}$ with differential $\partial'$ defined by:
\begin{quotation}
If $\partial\,u = \nu\,w + \eta\,v + r$ with $r$ linearly independent of $v$ and $w$, then $\partial'\,u = r - \nu\lambda^{-1}\,z$ 
\end{quotation}
where, for grading reasons one or other must be zero and we allow both $\nu, \eta = 0$.
\end{quotation} 
The result of applying this simplification to both vertical arrows in Figure \ref{fig:twisted} appears in Figure \ref{fig:tree}. For example, the bottom arrow
in Figure \ref{fig:tree} arises from the formula with  $v = v_{+}$ in $\mathcal{V}(L_{\{c_{2}\}})$, $w = v_{-}$ and $\lambda = x_{3}$\\
\ \\
\noindent For more general diagrams we use a standard result about Koszul complexes: that $\mathcal{K}(\area{C_{1}}, \ldots, \area{C_{k}})$ is acyclic over $\field{L}$ (\cite{Mats}). Consequently, when we simplify along non-zero components of $\partial_{\mathcal{V}}$ only those resolutions with $\mathcal{V}(L_{S}) = \field{L}v_{0}$ will remain to contribute. These are precisely the resolutions which consist of a single circle. \\
\ \\
\noindent Resolutions consisting of a single circle are in one-to-one correspondence with pairs of complementary spanning trees for the Tait graphs of $L$. These graphs are found by first bicoloring the faces of $L$ in a checkerboard fashion. Taking all the black faces as vertices, each crossing in $\cross{L}$ provides an edge since it abuts one or two black faces. The Tait graphs for $L$ are the two planar graphs obtained by repeating this constructions for both the black and white faces. Furthermore, the marked point on $L$ abuts one black and one white region. These regions identify a root vertex in each graph. Resolutions $L_{S}$ which consist of a single circle divide $S^{2}$ into two discs, one of each color, which are composed of the black faces and the white faces. If we take only those crossings which are resolved in $L_{S}$ as a merging of two black faces we obtain a spanning tree for the black graph. Likewise, if we take complementary set of edges, we obtain a spanning tree for the white graph. \\
\ \\
\noindent  Thus, in the homotopy equivalent complex, only those resolutions contribute generators which correspond to the two rooted spanning trees. We now  explicitly describe the chain complex in terms of the rooted spanning trees.\\
\ \\
\noindent Within $R(L)$ we distinguish those resolutions which result in a single circle in $S^{2}$:

$$
O(L) = \big\{S \subset \cross{L} \big| L_{S} \mathrm{\ is\ connected}\big\}
$$
\ \\
\noindent Furthermore, let $O_{i}(L) = O(L) \cap R_{i}(L)$. Elements of $O(L)$ will be typically be denoted by $T$, or a decorated variant. Given $T \in O(L)$
we let 
$$
O(T,L) = \big\{T' \in O(L) \big| T \subset T'\big\}
$$
\ \\
\noindent and $O_{i}(T,L) = O(T,L) \cap O_{i}(L)$. If $\degrad{T} = i$ then $O_{i+k}(T,L)$ are those resolutions such that $L_{T'}$ is a single (marked) circle and $T'\backslash T$ is a $k$ crossing subset of $\cross{L}\backslash T$. 

\noindent We now define the chain complex: take $$\CT[i]{L} = \mathrm{Span}_{\field{L}}\big\{ T \in O_{i}(L)\big\}$$ for each $i \in \Z$. Note that if $L$ is a split diagram, then $\CT{L} \cong 0$ as there are no such resolutions. There is boundary map 

$$\partial_{i,L}: \CT[i]{L} \longrightarrow \CT[i+2]{L}$$
\ \\
\noindent For each $T' \in O_{i + 2}(S,L)$, $T' \backslash T = \{c_{1}, c_{2}\}$ for two crossings $c_{1}, c_{2} \in \cross{L}\backslash{S}$. In $L_{T}$ these are depicted with resolution arcs, $a_{1}$ and $a_{2}$.  For $L_{T'}$ to be a single circle, $a_{1}$ and $a_{2}$ must have interlocking feet along the circle $L_{T}$. Since all the arcs are disjoint, one of the arcs must lie in each region of $S^{2} \backslash L_{T}$. Between $T$ and $T'$ there are two elements of $b, w \in R(L)$, with $b= T \cup \{c_{1}\}$ and $w = T \cup \{c_{2}\}$.  $L_{b}$ consists of two circles, found by surgering the arc $a_{1}$. One of these circles contains the marked point while the other determines a subset $B_{T,T'} \subset S^{2}$ {\em disjoint} from the marked point $p$, which is cleaved from $L_{T}$ by the change in resolution at $c_{1}$. Likewise, $L_{w}$ consists of two circles, one marked and the other containing a subset $W_{T,T'}$ disjoint from $p$. For the case of our extended example, see Figure \ref{fig:tree}. \\
\ \\
\noindent To $B_{T,T'}$ we assign the formal area $[B_{T,T'}]$  in $\field{L}$ which is the sum of the formal variables for the faces in $B_{T,T'}$. Similarly we can define an ares $[W_{T,T'}]$ for $W_{T,T'}$. These define the boundary map $\partial_{i,L}$

\begin{equation}
\partial_{i,L}\,T = \sum_{T' \in O_{i+2}(T,L)} \pairing{T}{T'} T'
\label{eqn:diff}
\end{equation}
\ \\
\noindent where
\begin{equation}
\pairing{T}{T'} = \frac{1}{[B_{T,T'}]} + \frac{1}{[W_{T,T'}]}
\end{equation}

\noindent This differential emerges from the reduction process previously described, but given its form we can also verify that it is a boundary map directly. This argument is combinatorial and is presented in section \ref{sec:boundary}. Moreover, we can verify that this complex {\em is invariant under changes of the marked point}. Together, these statements form

\begin{theorem}\label{thm:complex}
Let $L$ be the diagram for an oriented link with a marked point $p$. Let $\CT{L, p} = \bigoplus_{i \in \Z} \CT[i]{L, p}$, and $\partial_L$ be the map $\oplus \partial_{i,L}$. Then $\big(\CT{L}, \partial_L\big)$ is a chain complex, whose homology will be denoted $\unHT{L}$. The isomorphism type of the chain complex is invariant under changes of the marked point.  
\end{theorem}

\noindent The arguments for invariance of the twisted Khovanov homology are readily extended to $\field{L}$, and thus apply to the homotopy equivalent spanning tree complex. Consequently, we obtain the main result of this paper:

\begin{theorem}\label{thm:invariance}
Let $L$ be the diagram for an oriented link, $\lk{L}$. Then the (stable) isomorphism class of  $\unHT{L}[-n_{+}]$ is an oriented link invariant, denoted $\HT{L}$.
\end{theorem}

\subsection{Properties}

\noindent The reduced Khovanov homology has Euler characteristic equal to the Jones polynomial (for a suitable convention on the coefficients of the polynomial). The spanning tree complex also has a classical knot invariant as its Euler characteristic:

\begin{theorem}
For $\lk{L}$, a link in $S^{3}$, let
$$
P(t) = \sum_{j \in \Z} rk_{\F_{L}}\left(\HT[j]{L}\right)\delta^{j}
$$
then $det(\mathcal{L}) = |P(i)|$ where $i =\sqrt{-1}$.
\end{theorem}

\noindent The strange form of the Euler characteristic comes from $\partial_{L}$ being a $+2$-differential. However, using a $+2$-differential ensures that our gradings will occur in $\Z$ and not $\frac{1}{2}\Z$. This theorem can be derived either from relationship with Khovanov homology, where one interprets the addition of $\partial_{\mathcal{V}}$ as corresponding to evaluating the Jones polynomial at $-1$, or from the known relationship between the spanning trees of the Tait graph and the determinant of a link. We will opt for the latter in section \ref{sec:euler}.\\ 
\ \\
\noindent Second, we highlight the long exact sequence arising from the resolutions of a crossing in $L$. 

\begin{prop}
Let $L$ be oriented link diagram with crossing $c$. Let $L_{0}$ be the resolution of $L$ at $c$ according to the rule $c \not\in S$, and let $L_{1}$ be the resolution of $L$ at $c$ according to the rule $c \in S$. If $c$ is a positive crossing, and if $e=n_{+}(L) - n_{+}(L_{1})$ (for any orientation on $L_{1}$), then\\
\begin{equation}
\cdots \rightarrow \HTsub[i+e-1]{L}{1}\otimes \F_{L} \rightarrow \HT[i]{L} \rightarrow \HTsub[i+1]{L}{0}\otimes \F_{L} \stackrel{\tau_{c,\ast}}{\longrightarrow} \HTsub[i+e+1]{L}{1}\otimes \F_{L} \rightarrow \cdots
\end{equation}
\ \\
\noindent However, if $c$ is negative, and $f = n_{+}(L) - n_{+}(L_{0})$, then \\
\begin{equation}
\cdots \rightarrow \HTsub[i-1]{L}{1}\otimes \F_{L} \rightarrow \HT[i]{L} \rightarrow \HTsub[i+f]{L}{0}\otimes \F_{L} \stackrel{\tau_{c,\ast}}{\longrightarrow} \HTsub[i+1]{L}{1}\otimes \F_{L} \rightarrow \cdots
\end{equation}
\end{prop}

\noindent The groups in these long exact sequences are tensored with $\field{L}$ to have all groups over the same field. The proof of the proposition as well as the details of the actions of $\field{L_{0}}$ and $\field{L_{1}}$ are in section \ref{sec:LES}.  \\
\ \\
\noindent The long exact sequence can be used to replicate an argument of Manolescu and Ozsv\'ath, \cite{MaOz}. In particular,  

\begin{theorem}\label{thm:quasi}
If $\lk{L}$ represents a (quasi-)alternating link with a connected diagram, then $\HT[i]{L} \cong 0$ when $i \neq \sigma(L)$ and has rank $\mathrm{det(L)}$ when $i = \sigma(L)$, the signature of $\lk{L}$. 
\end{theorem}

\noindent We will use the convention that the signature of the right handed trefoil is $-2$. This result can also be proved from the more detailed connection with Khovanov homology given below. We also note the the spanning tree homology has two properties similar to those for other knot homologies. 

\begin{theorem}
Let $\lk{L}$ be an oriented link, then $\HT[i]{L} \cong \HT[-i]{\mir{L}}$.
\end{theorem}

\begin{theorem}
Let $\lk{L}_{1}$, $\lk{L}_{2}$ be two non-split oriented links, and let $\lk{L} = \lk{L}_{1} \# \lk{L}_{2}$, in some manner. Then
$$
\HT[k]{L} \cong \oplus_{i+j=k} \HTsub[i]{L}{1} \otimes \HTsub[j]{L}{2}
$$ 
where $\cong$ denotes stable equivalence.
\end{theorem}

\noindent Again, similar results should be provable directly from the totally twisted Khovanov homology. However, each of these three properties is proved in section \ref{sec:prop} using the spanning tree formalism. 

\subsection{Relationship with Khovanov homology}

\noindent There is a more precise relationship between the (characteristic 2, reduced) Khovanov homology and the totally twisted Khovanov homology, which extends also to the spanning tree homology. The $q$-grading, after we add $\partial_{\mathcal{V}}$, defines a filtration on the totally twisted Khovanov homology. We can examine the induced Leray spectral sequence to derive a relationship with Khovanov homology.

\begin{theorem}\label{thm:spectral}
The spectral sequence induced by the filtration from the $q$-index has $E^{0}$ page isomorphic to the $\delta$-graded reduced Khovanov homology $KH^{\ast,\ast}(L) \otimes \field{L}$ and converges, in finitely many steps, to $\HT{L}$. 
\end{theorem}

\noindent {\bf Proof:} On the unshifted, twisted Khovanov complex, the map $(i,j) \rightarrow j$ is a filtration and $\delta(i,j) = 2i - j$ is a grading. To compute the $E^{0}$-page we ignore the portion of the differential which changes the $j$-value. For us, this is the $(1,0)$ portion of the differential, $\partial_{KH}$. Consequently, the $E^{0}$ page is just the reduced Khovanov homology over the field $\field{L}$. The complex is bounded, so the corresponding spectral sequence converges to the total homology of the complex. The total homology when using the $\delta$ - grading is isomorphic to $\HT{L}$, since we are working over a field.  Finally, when we shift the Khovanov complex by $\{(-n_{-}, n_{+} - 2n_{-})\}$ to obtain the invariant homology, the $\delta$-grading shifts by $-n_{+} +2n_{-} + 2(-n_{-}) = -n_{+}$, which is how we calculated the shift to apply to $\widetilde{KT}_{\ast}(L)$ to obtain $KT_{\ast}(L)$. Thus in the spectral sequence above, we may use the Khovanov shifts on each bi-graded page, and this appropriately shifts the total grading so that the direct sum of the pieces with $\delta= k$ converges to $KT^{k}(L)$. $\Diamond$

\begin{cor}
If $KH_{s}(L)$ is the portion of the (reduced) Khovanov homology over $\field{L}$ in $\delta$-grading $s$, then $$\mathrm{rk}\,KH_{s}(L) \geq \mathrm{rk}\,\HT[s]{L}$$
\end{cor}

\noindent This corollary implies bounds om the Khovanov width of links in $S^{3}$, a result we examine in more depth in a later paper. For now, it is a natural question whether more can be said. To this end we relate some computations and the results of computer calculations of $\HT{L}$ and compare them with the results for characteristic 2, reduced Khovanov homology.

\subsection{\bf Computations:}\ \\ 

\subsubsection{\bf Unlinks:}\label{example:unlink} To finish our extended example, we compute the totally twisted homology $\underline{HT}(\lk{L})$ and the spanning tree homology $\HT{L}$. The spanning tree homology is trivial. Both of the arrows in Figure \ref{fig:tree} are multiplication by a non-zero element of $\field{L}$; the boundary operator is multiplication by their sum: $\frac{x_{3}+x_{4}}{x_{3}x_{4}}$. Since this element is also non-zero, the corresponding boundary map is an isomorphism, and the homology is trivial. This confirms the intuition that the homology should be trivial since the simplest diagram for the two component unlink is disconnected and thus unable to support any spanning trees for both its Tait graphs. On the other hand, for the twisted homology $\underline{HT}(\lk{L})$, the situation is more complicated. For the complex in Figure \ref{fig:twisted}, we can reduce two of the horizontal isomorphisms -- one out of $L_{\emptyset}$ and the other into $L_{\{c_{1}, c_{2}\}}$ -- to be left with either
$$
0 \longrightarrow \ring{L} \stackrel{\cdot x_{3}}{\longrightarrow} \ring{L} \longrightarrow 0
$$
with homology $\Z/2\Z[x_{1},x_{2},x_{4}]$, as a $\ring{L}$-module where $x_{3}$ acts by multiplication by $0$, or, if we choose the other pair,
$$
0 \longrightarrow \ring{L} \stackrel{\cdot x_{4}}{\longrightarrow} \ring{L} \longrightarrow 0
$$
with homology $\ring{L}$-module $\Z/2\Z[x_{1},x_{2},x_{3}]$, where $x_{4}$ acts by multiplication by $0$. Furthermore, in the usual geometrically split 
diagram for the two component unlink, we would have two circles, enclosing regions $y_{1}$ and $y_{2}$, say, with that enclosing $y_{1}$ being the marked circle. The homology module would then be the homology of the complex
$$
0 \longrightarrow \Z/2\Z[y_{1},y_{2}] \stackrel{\cdot y_{2}}{\longrightarrow} \Z/2\Z[y_{1},y_{2}] \longrightarrow 0
$$
The homology of this complex is $\Z/2\Z[y_{1}]$ as a $\Z/2\Z[y_{1},y_{2}]$-module with $y_{2}$ acting by multiplication by $0$. These complexes and their homologies will be related through the notion of stable isomorphism over polynomial rings, under which they are all equivalent to the last module.\\
\ \\
\noindent Using the invariance results above, the completely split diagram for the $n$-component unlink will give a complex isomorphic to $\mathcal{K}(y_{2},\ldots, y_{n})$ over $P = \Z/2\Z[y_{1}, \ldots, y_{n}]$, where $y_{1}$ corresponds to the bounded region enclosed by the marked circle. Since $y_{2}, \ldots, y_{n}$ is a regular sequence over $P$, standard results about Koszul complexes (see \cite{Mats}, Theorem 16.5 for example) imply 
$$
\underline{HT}_{k}(\lk{L}) \cong \left\{\begin{array}{cl} \Z/2\Z[y_{1}] & k = 0 \\ 0 & k \neq 0 \end{array} \right.
$$
which determines the stable isomorphism class of the totally twisted Khovanov homology. In fact, we can take this homology to be $\Z/2\Z$ in grading $0$, as a $\Z/2\Z$-module, considered up to stable equivalence.\\
\ \\
\subsubsection{\bf Results for knots:} A computer can quickly calculate the ranks of the spanning tree homology for knots up to 15 crossings. There is a mild difficulty in that computations over $\field{L}$ are not very efficient due to the large number of variables that may need to be tracked. The details of the workaround and results of these computer surveys will be highlighted in the sequel to this paper, as will various generalizations. Here we relate some of the data with an eye to understanding $\HT{K}$ more fully, where $K$ is a {\em knot}. In particular, theorem \ref{thm:spectral} suggests that we compare $\HT{K}$ to the $\delta$-graded, characteristic 2, reduced Khovanov homology, which we will denote $KH^{r}_{\delta}(L)$.\ \\
\\ 
\noindent{\bf Caution:} In the Khovanov homology literature, there is another $\delta$-grading which is related but not identical to ours. In particular, ours is two times the negation of that $\delta$-grading. The choice here ensures that the gradings are all integers, although with a $+2$-differential, and corresponds to the number of elements in $S$ for the resolution $L_{S}$. Below we implicitly convert the usual $\delta$-grading to the convention employed in this paper.\\
\ \\

\noindent We will describe the homology $\HT{K}$ by it Poincare\'e polynomial: 
$$
\sum_{j \in \Z} rk_{\F_{L}}\left(\HT[j]{L}\right)\delta^{j}
$$
which indicates the ranks and the gradings. Since we are currently interested in stable equivalence over fields, the ``graded ranks'' are all that remain.  
\begin{enumerate}
\item Proposition \ref{thm:quasi} shows that for (quasi-)alternating knots and links, $\HT{L}$ has the same rank in each $\delta$-grading as $KH^{r}_{\delta}(K)$. The analogous theorem for $KH^{r}_{\delta}(L)$ was proven for alternating links by E. S. Lee, and for quasi-alternating links
by Manolescu and Ozsv\'ath, \cite{MaOz}. 
\item Computer computations verify that the rank of $\HT{K}$ in each $\delta$-grading is the same as that for $KH^{r}_{\delta}(K)$ for every non-alternating knot with 13 or fewer crossings.
\item The torus knots $T_{5,3}$, $T_{7,3}$, and $T_{5,4}$ also have the rank of $\HT{K}$ in each $\delta$-grading the same as that for $KH^{r}_{\delta}(K)$. For $T_{5,3}$ the common Poincar\'e polynomial is
$$
HT_{\ast}(T_{5,3}): \hspace{0.5in} 4\delta^{-8} + 3\delta^{-6} 
$$
found from a diagram yielding 27 spanning trees as generators. Likewise, for $T_{7,3}$ both theories have Poicar\'e polynomial 
$$
HT_{\ast}(T_{7,3}): \hspace{0.5in}4\delta^{-12} + 4\delta^{-10} + \delta^{-8}
$$
found from a diagram yielding 841 spanning trees as generators. For $T_{5,4}$ the homology is
$$
HT_{\ast}(T_{5,4}): \hspace{0.5in} 4\delta^{-12} + 4\delta^{-10} + 5 \delta^{-8}
$$
found from a diagram yielding 1805 spanning trees as generators. It is worth describing the chain complex for $T_{5,4}$. The table below lists the non-zero number of generators in each $\delta$-grading for the diagram used in the computation:
\begin{center}
\begin{tabular}{|c|c|c|c|c|c|}
\hline
$\delta$ & -12 & -10 & -8 & -6 & -4 \\
\hline
\# gens & 125 & 500 & 700 & 400 & 80 \\
\hline
\end{tabular}
\end{center}

\item In \cite{Rasm}, there are several examples of knots which have identical $\delta$-graded, reduced Khovanov homology and knot Floer homology (for an analogous choice of $\delta$-grading), as well as several examples where they differ. Many of these examples are in the list of knots with 15 or fewer crossings, but several are not. Those where we have computed both homologies do not provide any examples of knots where the ranks of $\HT{K}$ differ from those of $KH^{r}_{\delta}(K)$. Where they agree we have the following results:
\begin{center}
\begin{tabular}{|l|c|l|}
\hline
Knot & \# Trees & $\HT{K} \cong KH^{r}_{\delta}(K)$ \\
\hline
\hline
$11_{n}19$ & $65$ & $3\delta^{2} + 8 \delta^{4}$ \\
$11_{n}38$ & $75$ & $8\delta^{-2} + 5\delta^{0}$ \\
$11_{n}57$ & $95$ & $12\delta^{-6} + 5 \delta^{-4}$\\
$11_{n}79$ & $87$ & $16\delta^{-2} + \delta^{0}$ \\
$12_{n}121$ & $89$ & $8\delta^{0} + 7 \delta^{2}$\\
$12_{n}502$ & $153$ & $12\delta^{-8} + 3\delta^{-6}$\\
$12_{n}591$ & $175$ & $4\delta^{-8} + 11\delta^{-6}$\\
$12_{n}725$ & $75$ & $4\delta^{-10} + 9\delta^{-8}$\\
$12_{n}749$ & $167$ & $4\delta^{-4}+11 \delta^{-2}$\\
$13_{n}192$ & $87$ & $5\delta^{4} + 8\delta^{6}$\\
$14_{n}21882$ & $249$ & $4\delta^{-2} + 13\delta^{0}$\\
$15_{n}4863$ & $109$ & $7\delta^{6} + 8\delta^{8}$\\
$15_{n}41127$ & $1203$ & $8\delta^{-4} + 8 \delta^{-2} + \delta^{0}$\\
$15_{n}80764$ & $401$ & $\delta^{0} + 8\delta^{2} + 8\delta^{4}$\\
\hline
\end{tabular}
\end{center}
\ \\
\noindent For the knots listed in \cite{Rasm} where the knot Floer homology differs from $KH^{r}_{\delta}(K)$, we consider the $(2,5)$-cable of the positive trefoil, $C_{2,5}T$, and the $(2,7)$-cable of the positive trefoil, $C_{2,7}T$. These are the knots $13_{n}4639$ and $13_{n}4587$, respectively. For these we have
$$
HT_{\ast}(C_{2,5}T): \hspace{.5in} 4\delta^{-8} + 7\delta^{-6} + 8\delta^{-4}
$$
and
$$
HT_{\ast}(C_{2,7}T): \hspace{.5in} 4\delta^{-10} + 5\delta^{-8} + 8\delta^{-6}
$$
In addition, the torus knot $T_{4,5}$ also has different $\delta$-graded knot Floer and Khovanov homology. Nevertheless, all of these still have identical Poincare\'e polynomials for the spanning tree and reduced, characteristic 2, $\delta$-graded Khovanov homologies.
\end{enumerate}

\noindent For {\em knots}, these results raise the following question:
\begin{quote}
 For a knot $K$, is $\HT{K}$ (stably) isomorphic to $KH^{r}_{\delta}(K)$?
\end{quote} 
\noindent An affirmative answer means that the spanning tree homology provides an explicit spanning tree model for this version of Khovanov homology. A negative answer would mean that $\HT{K}$ is a new, and computable knot invariant. Even though the author has not found a counter-example, the negative answer seems more likely given the results for links. 

\subsubsection{\bf Results for Links:} One can make the same comparison for links as for knots. It is straightforward to see that the homologies for the two component unlinks are different (see Example \ref{example:unlink}). However, even among non-split links there are discrepancies between the homologies as soon as non-alternating links appear in the link tables. Of the 1424 links with 11 or fewer crossings found on the KnotAtlas website, 200 have different Poincar\'e polynomials for the spanning tree and reduced, characteristic 2 Khovanov homology. \\
\ \\
\noindent For  instance,  $L6_{n}1$, $L7_{n}1$, and $L8_{n}8$ are depicted from left to right below:\\
\ \\
\begin{center}
\treefig{0.3}{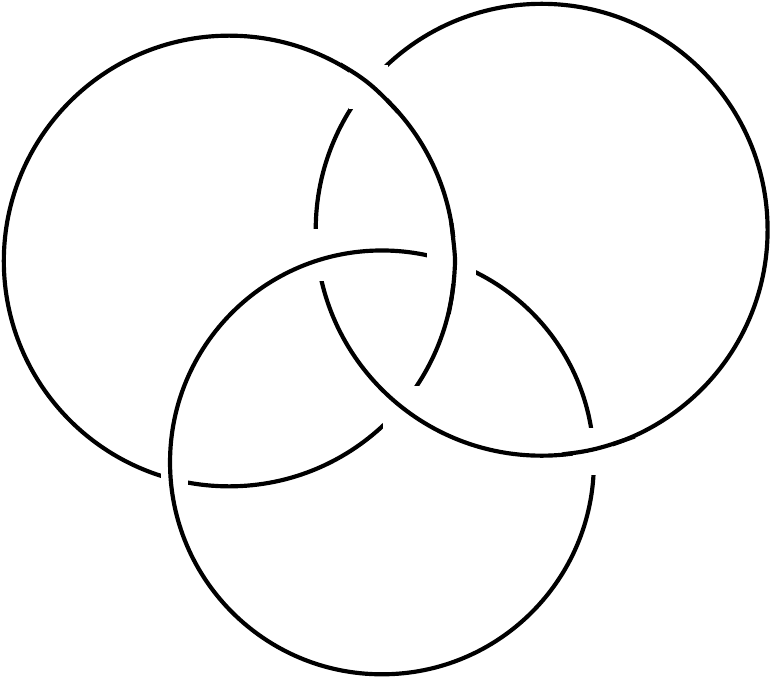} \hspace{.5in}
\treefig{0.15}{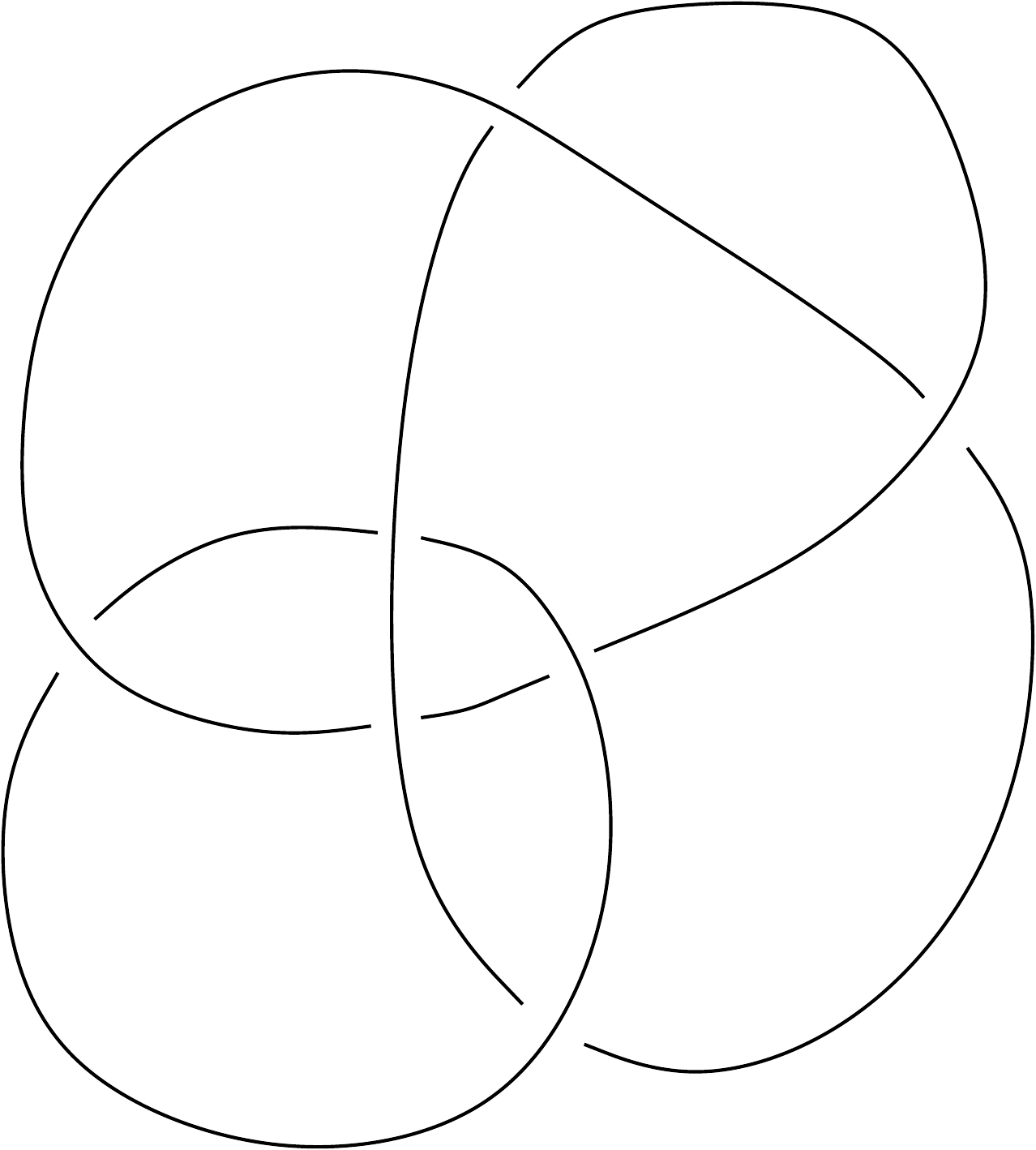} \hspace{.5in}
\treefig{0.2}{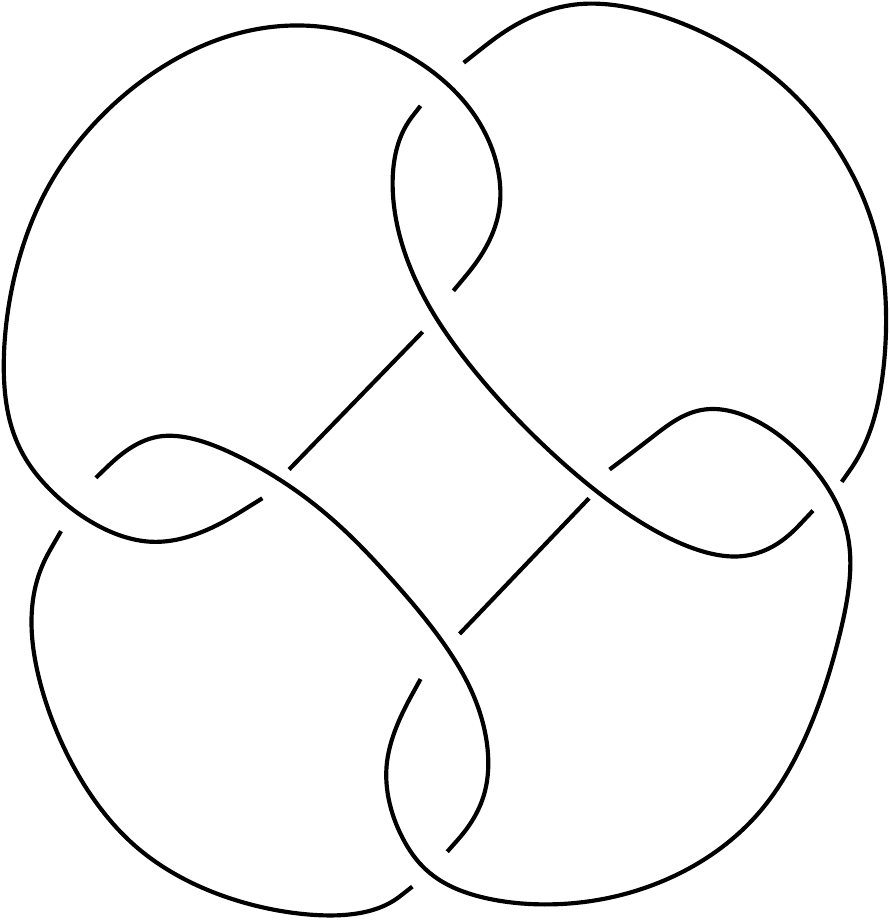}
\end{center}
\ \\
\noindent $\mathrm{HT}_{\ast}(L6_{n}1)$ has Poincar\'e polynomial $4\delta^{0}$, whereas for the $\delta$-graded, reduced, characteristic 2 Khovanov homology it is $\delta^{-2} + 5\delta^{0}$. Thus there is a higher differential in the spectral sequence in theorem \ref{thm:spectral}. For $L7_{n}1$ a similar reduction occurs: the spanning tree homology has polynomial $4\delta^{5}$ whereas the Khovanov homology has $\delta^{3} + 5\delta^{5}$. For $L8_{n}8$ we have $4\delta^{-1} + 4 \delta^{1}$ for the spanning tree homology, whereas the Khovanov homology is $6\delta^{-1} + 6\delta$. For $L8_{n}6$ we have spanning tree homology $4\delta^{2} + 4\delta^{4}$ whereas the Khovanov homology is $5\delta^{2} + 5\delta^{4}$. No obvious property of a link explains these reductions. $L6_{n}1$ has three components whereas $L7_{n}1$ has only two; nevertheless, they have the same rank difference when compared to their Khovanov homologies. $L8_{n}8$ and $L8_{n}6$ both have nullity 1, but their rank comparisons are different. On the other hand, $L6_{n}1$ has nullity $0$ but the same rank reduction as $L8_{n}8$. The author has no ready explanation for the occurrence of the higher differentials for links, nor an explanation for why they occur frequently, whereas no higher differentials have been found in the spectral sequence for knots.

\section{Totally twisted Khovanov homology}\label{sec:chainmaps}

\begin{theorem}
Let $\partial_{KH} : \widetilde{KT}^{\ast, \ast}(L) \longrightarrow \widetilde{KT}^{\ast+1,\ast}(L)$ be the Khovanov differential, and let $\partial_{\mathcal{V}} : \widetilde{KT}^{\ast, \ast}(L) \longrightarrow \widetilde{KT}^{\ast, \ast - 2}(L)$ be the map $\oplus \partial_{\mathcal{V}(L_{S})}$. Then $\partial_{KH} + \partial_{\mathcal{V}}$ is a boundary map on $\widetilde{KT}(L)$. 
\end{theorem}
 
\noindent {\bf Proof:} To show that $\big(\partial_{KH} + \partial_{\mathcal{V}}\big)^{2} \equiv 0$ we need $\partial_{KH}\circ\partial_{\mathcal{V}} = \partial_{\mathcal{V}} \circ \partial_{KH}$ since $\ring{L}$ has characteristic 2. However,
$$\partial_{KH} : \mathcal{V}(L_{S}) \longrightarrow \bigoplus_{i \in \cross{L}\backslash S} \mathcal{V}(L_{\{i\} \cup S})\{+1\}$$
so it suffices to verify that $\partial_{KH, S}: \mathcal{V}(L_{S}) \longrightarrow \mathcal{V}(L_{\{i\} \cup S})\{+1\}$ is a chain map for each $i\in \cross{L}\backslash S$. Furthermore, since $\mathcal{V}(L_{S})$ is itself a tensor product, we can verify that $\partial_{KH}$ is a chain map through three lemmas addressing its affect on each of the factors. These verify that $\partial_{KH}$ is a chain map for factors corresponding to circles which are not merging or dividing, and for factors corresponding to circles that merge or that divide. Furthermore, implicitly each argument will allow one of the circles to be the marked circle. 

\begin{lemma}
If $C$ is a circle which is unaffected by changing the resolution from $L_{S}$ to $L_{\{i\}\cup S}$, then $\partial_{KH}$ induces a chain map on $\mathcal{V}_{C}$. 
\end{lemma} 

\noindent {\bf Proof:} The map $\partial_{KH}$ induces on $V_{C}$ is the identity. Furthermore, $\area{C}$ is the same for both $L_{S}$ and $L_{S \cup\{i\}}$ since the circle is unchanged. Thus, the map $\partial_{KH}$ induces on $\mathcal{V}_{C}$ is a chain map. $\Diamond$.\\
\ \\

\begin{center}
\begin{figure}

\treefig{0.7}{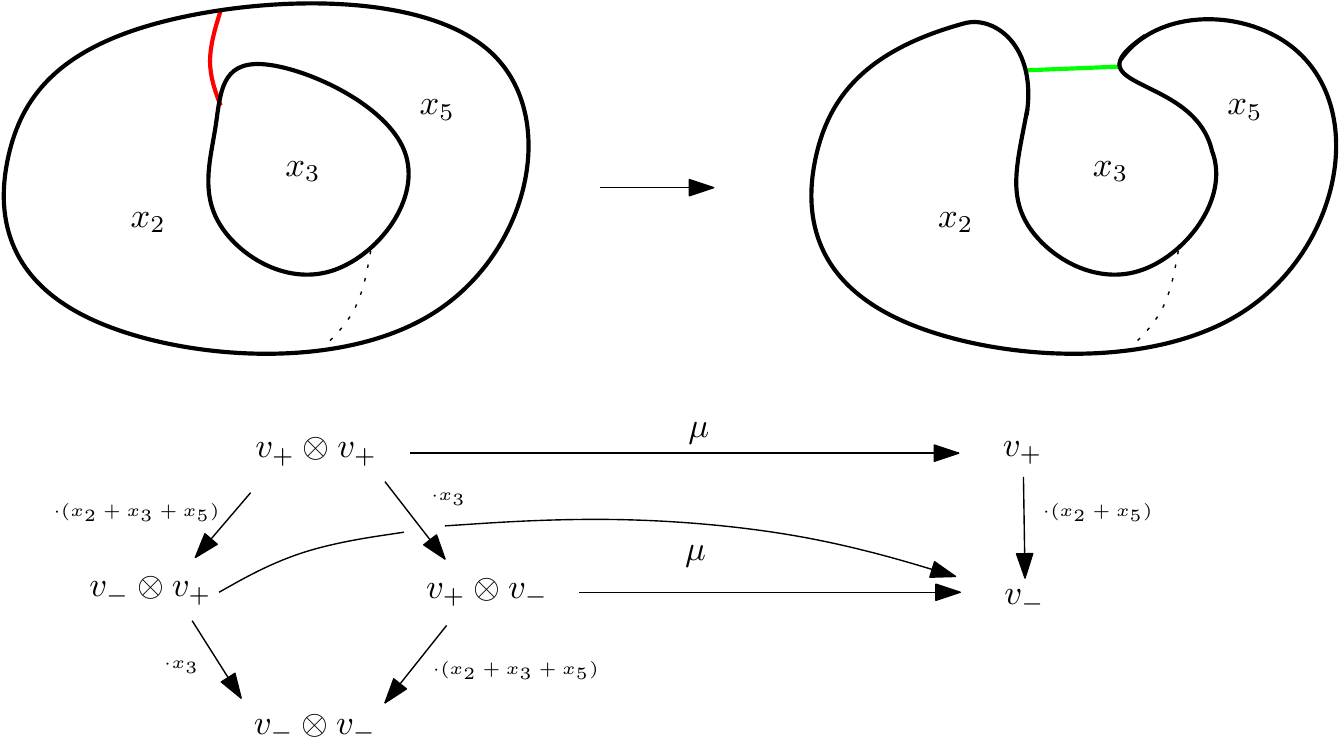} 

\vspace{0.5in} 

\treefig{0.7}{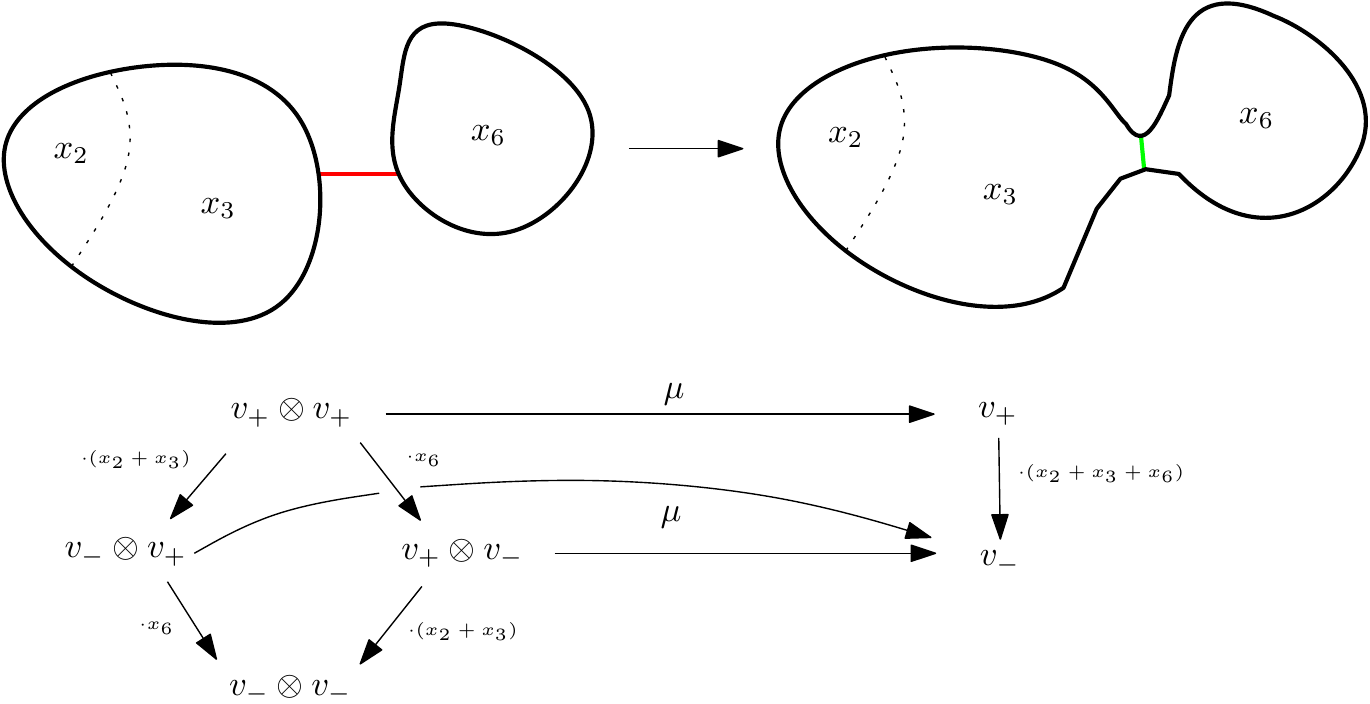}

\caption{Examples of the two cases considered in the proof of Lemma \ref{lem:mu}. The first figure depicts case 1) $\disc{C_{1}} \subset \disc{C_{2}}$ along with a simplified version of the chain complexes and maps involved, while the bottom figure depicts case 2), $\disc{C_{1}} \cap \disc{C_{2}} = \emptyset$. }
\label{fig:merge}
\end{figure}
\end{center}

\begin{lemma} \label{lem:mu}
Let $C_{1}$ and $C_{2}$ be circles in $L_{S}$ which merge into a single circle $C$ in $L_{\{i\}\cup S}$. The map 
$$
\mu : \mathcal{V}_{C_{1}} \otimes \mathcal{V}_{C_{2}} \longrightarrow \mathcal{V}_{C}\{1\}
$$
is a chain map.
\end{lemma}

\noindent{\bf Proof:} When $C_{1}$ and $C_{2}$ are both unmarked circles, there are two cases to consider: 1) $\disc{C_{1}} \subset \disc{C_{2}}$ (or $\disc{C_{2}} \subset \disc{C_{1}}$), and 2) $\disc{C_{1}} \cap \disc{C_{2}} = \emptyset$. Examples of the argument in the two cases are given in Figure \ref{fig:merge}. In each case we start by computing $\area{C}$. When $\disc{C_{1}} \cap \disc{C_{2}} = \emptyset$ then merging $C_{1}$ and $C_{2}$ produces a circle $C$ with $\disc{C} = \disc{C_{1}}\#\disc{C_{2}}$, the boundary connect sum of $\disc{C_{1}}$ and $\disc{C_{2}}$. But then any face $f \in \regions{L}$ with $f \cap \disc{C} \neq \emptyset$ has either $f \cap \disc{C_{1}} \neq \emptyset$ or $f \cap \disc{C_{2}} \neq \emptyset$, but {\em not both}. Thus $\area{C} = \area{C_{1}} + \area{C_{2}}$. In the other case, if $\disc{C_{1}} \subset \disc{C_{2}}$, merging $C_{1}$ and $C_{2}$ results in $\disc{C} = \disc{C_{2}}\backslash \disc{C_{1}}$, but $[\disc{C_{2}}\backslash\disc{C_{1}}] + \area{C_{1}} = \area{C_{2}}$. Since we are working in characteristic 2 
$$
\area{C} =  [\disc{C_{2}}\backslash\disc{C_{1}}] = \area{C_{1}} + \area{C_{2}}
$$
Thus in both cases $\area{C} = \area{C_{1}} + \area{C_{2}}$. With this result, we can easily verify that $\partial_{\mathcal{V}_{C}} \circ\mu = \mu \circ \big(\partial_{\mathcal{V}_{C_{1}}}\otimes \mathbb{I} + \mathbb{I} \otimes \partial_{\mathcal{V}_{C_{2}}}\big)$.\\
\ \\
\noindent The map $\partial_{\mathcal{V}_{C}} \circ \mu$ computed on generators of $\mathcal{V}_{C_{1}} \otimes \mathcal{V}_{C_{2}}$ equals
$$
\begin{array}{l}
v_{+} \otimes v_{+} \stackrel{\mu}{\longrightarrow} v_{+} \stackrel{\partial_{\mathcal{V}_{C}}}{\longrightarrow} \area{C}v_{-}\\
\ \\
v_{+} \otimes v_{-}, v_{-} \otimes v_{+} \stackrel{\mu}{\longrightarrow} v_{-} \stackrel{\partial_{\mathcal{V}_{C}}}{\longrightarrow} 0 \\
\ \\
v_{-} \otimes v_{-} \stackrel{\mu}{\longrightarrow} 0 \stackrel{\partial_{\mathcal{V}_{C}}}{\longrightarrow} 0
\end{array}
$$
On the other hand, if we apply $\partial = \partial_{\mathcal{V}_{C_{1}}}\otimes \mathbb{I} + \mathbb{I} \otimes \partial_{\mathcal{V}_{C_{2}}}$ first, and then $\mu$ we obtain
$$
\begin{array}{l}
v_{+}\otimes v_{+} \stackrel{\partial}{\rightarrow} \left(\area{C_{1}}\,v_{-}\otimes v_{+} + \area{C_{2}}\,v_{+} \otimes v_{-}\right) \stackrel{\mu}{\rightarrow} \big(\area{C_{1}} + \area{C_{2}}\big) v_{-}\\
\ \\
v_{+} \otimes v_{-} \stackrel{\partial}{\rightarrow} \area{C_{1}}\,v_{-} \otimes v_{-} \stackrel{\mu}{\rightarrow} 0\\
\ \\
v_{-} \otimes v_{+} \stackrel{\partial}{\rightarrow} \area{C_{2}}\,v_{-} \otimes v_{-} \stackrel{\mu}{\rightarrow} 0\\
\ \\
v_{-} \otimes v_{-} \stackrel{\partial}{\rightarrow} 0 \stackrel{\mu}{\rightarrow} 0
\end{array}
$$
Since $\area{C} = \area{C_{1}} + \area{C_{2}}$ these two maps equal. Now let $C_{1}$ be the marked circle. Then $\mathcal{V}_{C_{1}}$ is spanned by $v_{0}$. $v_{0}$ behaves identically to $v_{-}$ in $\mu$, but both of the above maps have image equal to 0 if one of the generators equals $v_{-}$. Thus $\mu$ is also a chain map when merging the marked circle. $\Diamond$\\

\begin{lemma} 
Let $C$ be an unmarked circle in $L_{S}$ which divides into two circles $C_{1}$ and $C_{2}$ in $L_{\{i\}\cup S}$. The map 
$$
\Delta : \mathcal{V}_{C} \longrightarrow \left(\mathcal{V}_{C_{1}}\otimes \mathcal{V}_{C_{2}}\right)\{1\}
$$
is a chain map.
\end{lemma}

\noindent{Proof:} $V_{C}$ is a Frobenius algebra with multiplication $\mu$ and co-unit $\epsilon: V_{C} \longrightarrow \F$ given by $\epsilon(v_{-}) = 1$ and $\epsilon(v_{+}) = 0$. That $V_{C}$ is a Frobenius algebra implies that $\epsilon \circ \mu$ is a non-degenerate bilinear form which induces an isomorphism $\lambda : V \rightarrow V^{\ast}$. The co-multiplication $\Delta$ for a Frobenius algebra is the map obtained from $\mu: V \otimes V \rightarrow V$ by dualizing , $V^{\ast}\otimes V^{\ast} \leftarrow V^{\ast}$, and then identifying $V^{\ast}$ with $V$ using $\lambda^{-1}$ to obtain $V \rightarrow V \otimes V$. Let $f_{\pm}$ be the basis dual to $v_{\pm}$ in $V^{\ast}$. The map $\partial_{\mathcal{V}}$ induces a map $\partial^{\ast}$ which can be computed as $\partial^{\ast}(f_{-}) = \area{C} f_{+}$, and otherwise $0$. $\lambda(v_{\pm}) = f_{\mp}$, so $\lambda$ is a chain map $\mathcal{V} \rightarrow \mathcal{V}^{\ast}_{C}$. Since $\mu$ is a chain map, it follows easily that $\Delta$ is likewise a chain map. $\Diamond$ \\
\ \\
\noindent Of course, this last lemma can also be verified directly, using the same method as in lemma \ref{lem:mu} and the result on $\area{C}$. \\
\ \\
\noindent These three lemmas imply that the building blocks of $\partial_{KH}$ are chain maps on the factors in $\mathcal{V}(L_{S})$. Consequently, $\partial_{KH} + \partial_{\mathcal{V}}$ is a differential on $\widetilde{K}$. $\Diamond$\\
\ \\
\section{Preliminaries for proving invariance}

\noindent In this section we try to show that $KT_{\ast}(L)$ is a link invariant if $\lk{L}$. However, due to the presence of the additional marked point, the result is not quite so strong. Instead we will show that

\begin{theorem}
Let $L$ be the diagram for a link $\lk{L}$ in $S^{2}$, equipped with a marked point $p$. The (stable) chain homotopy type of $KT_{\ast}(L)$ is an invariant of $L$ under Reidemeister moves and planar isotopies in $S^{2}\backslash \{p\}$  
\end{theorem}

\noindent In other words, as long as the isotopies do not cross the marked point, the homology is an invariant. The author does not know if the twisted homology is invariant under changes of marked point. However, for a different set of coefficients, we will be able to prove this. \\
\ \\
\noindent Different projections on $\lk{L}$ can have different numbers of faces, and thus the corresponding twisted complexes occur with non-isomorphic coefficient rings, $\ring{L}$. We will need an appropriate algebraic equivalence to relate the complexes for different projection. We start by describing this equivalence. We then examine some useful technical results. After the technical results, we prove the theorem by adapting the usual proofs of invariance for Khovanov homology, \cite{Khov}, \cite{Bar2}.  \\
 
\subsection{Stable equivalence over polynomial rings}

\noindent Let $W$ be a vector space over $\F$, and let $P_{W} = \mathrm{Sym}(W)$ be its symmetric algebra. $P_{W}$ is an integral domain, so we may find it field of fractions, $F_{W}$. A basis for $W$ identifies $P_{W}$ with a commutative polynomial ring generated by the basis elements, and $F_{W}$ with the corresponding field of rational functions. Thus, we will sometimes refer to $F_{W}$ as $\mathrm{Rat}(W)$ when we wish to emphasize this connection. Any linear map $A : W \longrightarrow W'$ induces a map $\mathrm{Sym}(A):P_{W} \longrightarrow P_{W'}$. If $A$ is also an injection, then $A$ also induces a map $\mathrm{Rat}(A):F_{W} \longrightarrow F_{W'}$ since $\mathrm{Sym}(A)$ has trivial kernel. When $A$ is an isomorphism then $\mathrm{Sym}(A)$ and $\mathrm{Rat}(A)$ are also isomorphisms in the appropriate category.\\
\ \\
\noindent If $M$ is a module over $P_{W}$, and $A : W \longrightarrow W'$, then $M \otimes_{P_{W}} P_{W'}$ is a module over $P_{W'}$ where $(p \cdot m) \otimes p' = m \otimes (\mathrm{Sym}(A)(p)\cdot p')$ and the action of $P_{W'}$ occurs on the second factor. Likewise if $V$ is a vector space over $F_{W}$, then $V \otimes_{F_{W}} F_{W'}$ is a vector space over $F'$.  

\begin{defn}
Let $W$ and $W'$ be two $\F$-vector spaces. A module $M$ over $P_{W}$ is stably isomorphic to a module $M'$ over $P_{W'}$ if there is an $\F$-vector space $W''$, and injections $i,i':W,W' \hookrightarrow W''$ which induce an isomorphism $$M \otimes_{P_{W}} P_{W''} \cong M' \otimes_{P_{W'}} P_{W''}$$ as $P_{W''}$ modules. When we wish to identify $W''$ and the injections we will say that $W$ is stably isomorphic to $W'$ through $(W'', i, i')$. 
\end{defn}
We will consider this as a relation on pairs $(M,W)$, although we will often omit reference to $W$ when it is clear in the context.

\begin{lemma}
$(M, W)$ stably isomorphic to $(M',W')$ is an equivalence relation on modules $M$ over the rings $P_{W}$, when $W$ is an $\F$-vector space. 
\end{lemma}

\noindent{\bf Proof:} The identity and symmetry of the relation are clear in the definition. We need only verify transitivity. Suppose that $(M,W)$ is stably isomorphic to $(M',W')$ through $(\widetilde{W}_{1}, i, i'_{1})$ and $(M',W')$ is stably isomorphic to $(M'', W'')$ through $(\widetilde{W}_{2}, i'_{2}, i'')$. Let $\widetilde{W}$ be the quotient of $\widetilde{W}_{1} \oplus \widetilde{W}_{2}$ by the subspace $\big\{\,i'_{1}(w) \oplus \vec{0} - \vec{0} \oplus i'_{2}(w) \,\big|\, w \in W' \,\big\}$. Then projection  onto the quotient composed with $i \oplus \vec{0}$ is an injection $ W \hookrightarrow \widetilde{W}$, since it maps entirely into the first factor, while no non-trivial element in the subspace is entirely in the first factor. Likewise, $\vec{0} \oplus i''$ induces an injection  $W'' \hookrightarrow \widetilde{W}$ is an injection. Furthermore, $i'_{1} \oplus 0$ and $0 \oplus i'_{2}$ induce injections of $W'$ into $\widetilde{W}$ with the same image. We now consider $T = \big(M \otimes_{P_{W}} P_{\widetilde{W}_{1}}\big) \otimes_{P_{\widetilde{W}_{1}}} P_{\widetilde{W}}$. On the one hand, $T$ is isomorphic to $M \otimes_{P_{W}} \big(P_{\widetilde{W}_{1}} \otimes_{P_{\widetilde{W}_{1}}} P_{\widetilde{W}} \big)$, but since $W_{1}$ injects into $\widetilde{W}$, $P_{\widetilde{W}_{1}} \otimes_{P_{\widetilde{W}_{1}}} P_{\widetilde{W}} \cong P_{\widetilde{W}}$. Here the action of $P_{W}$ on $P_{\widetilde{W}}$ is given by symmetric power of the composition of the inclusion maps. On the other hand stable equivalence implies $T$ is isomorphic to $T'= \big(M' \otimes_{P_{W'}} P_{\widetilde{W}_{1}}\big) \otimes_{P_{\widetilde{W}_{1}}}P_{\widetilde{W}}$ where the action of $P_{\widetilde{W}_{1}}$ on $P_{\widetilde{W}}$ is by the inclusion $Id \oplus \vec{0}$ followed by projection, and the action of $P_{W'}$ on $P_{\widetilde{W}_{1}}$ is by $i'_{1}$.  Re- organizing the tensor product using associativity, as before we obtain an isomorphism with $M' \otimes_{P_{W'}} P_{\widetilde{W}}$ where the action of $P_{W'}$ on $P_{\widetilde{W}}$ is given by the symmetric power of $pr \circ (i'_{1} \oplus \vec{0})$. \\
\ \\
\noindent We can perform the same argument starting with $\big(M'' \otimes_{P_{W''}} P_{\widetilde{W}_{2}}\big) \otimes_{P_{\widetilde{W}_{2}}} P_{\widetilde{W}}$. This is isomorphic to $M'' \otimes_{P_{W''}} P_{\widetilde{W}}$ with action given by the symmetric power of the composition $W'' \hookrightarrow \widetilde{W}_{2} \hookrightarrow \widetilde{W}$. It is, as above, also isomorphic to $M' \otimes_{P_{W'}} P_{\widetilde{W}}$ with action of $P_{W'}$ on $P_{\widetilde{W}}$ given by the symmetric power of $pr \circ (\vec{0} \oplus i_{2}')$.  However, $pr \circ (\vec{0} \oplus i_{2}') = pr \circ (i_{1}' \oplus \vec{0})$, so this is also isomorphic to $M' \otimes_{P_{W'}} P_{\widetilde{W}}$ where the action of $P_{W'}$ on $P_{\widetilde{W}}$ is given by the symmetric power of $pr \circ (i'_{1} \oplus \vec{0})$. From the preceding paragraph, we can conclude that $M \otimes_{P_{W}} P_{\widetilde{W}}$ using the composition $I: W \hookrightarrow \widetilde{W}_{1} \hookrightarrow \widetilde{W}$ is isomorphic, as a $P_{\widetilde{W}}$-module, to $M'' \otimes_{P_{W''}} P_{\widetilde{W}}$ using the inclusion $I'': W'' \hookrightarrow \widetilde{W}_{2} \hookrightarrow \widetilde{W}$. In particular, $(M,W)$ and $(M'', W'')$ 
are stably isomorphic through $(\widetilde{W}, I, I'')$. $\Diamond$\\
\ \\
\begin{lemma}
If $M$ is a free module over $P_{W}$, $M'$ is a free module over $P_{W'}$, and $(M,W)$ is stable isomorphic to $(M',W')$ then $\mathrm{dim}_{P_{W}} M = \mathrm{dim}_{P_{W'}}M'$
\end{lemma}

\noindent{\bf Proof:} Let $\big\{\,e_{i} \in M\,\big|\,i \in \Lambda\,\big\}$ be a basis for $M$ over $P_{W}$. Then since $\left(\sum a_{i}e_{i} \right) \otimes w$
equals $\sum e_{i} \otimes \mathrm{Sym}(I_{1})(a_{i})\cdot w = \sum (e_{i} \otimes 1) \cdot (\mathrm{Sym}(I_{1})(a_{i})\cdot w) $, we see that $\big\{\, e_{i} \otimes 1\,\big|\, i \in \Lambda\,\big\}$ is a basis for $M \otimes_{P_{W}} P_{\widetilde{W}}$ over $P_{\widetilde{W}}$, consequently it has the same rank. Performing the same calculation for $M'$ we see that the ranks must be equal. $\Diamond$.\\
\ \\

\noindent A similar equivalence relation holds for vector spaces over $F_{W}$. When two vector spaces are stably isomorphic their dimensions over their respective fields satisfy .  If $V$ is a graded vector space over $F_{W}$, stable isomorphism induces stable isomorphism in each grading, and the rank equality holds in each grading.  \\
\\ 
\noindent We will use stable isomorphism to relate chain complexes of modules defined over $P_{W}$ for differing vectors spaces $W$. To do this we note the following commutative algebra result:

\begin{lemma}
Let $I: W \hookrightarrow W'$ be an injection of finite dimensional vector spaces. Then $P_{W'}$ is free, hence flat, as a $P_{W}$-module
\end{lemma} 

\noindent{\bf Proof:} We can select a basis for $W$, $\{w_{1}, \ldots, w_{k}\}$ and consider $S= \{I(w_{1}), \ldots, I(w_{k})\}$ in $W'$. $S$ is linearly independent over $\F$, and thus can be extended to a basis for $W'$ by appending some vectors $\{y_{k+1}, \ldots, y_{l}\}$. With these choices, there are ring isomorphisms $P_{W} \cong \F[w_{1}, \ldots, w_{k}]$ and $P_{W'} \cong \F[I(w_{1}), \ldots, I(w_{k}), y_{k+1}, \ldots, y_{l}]$. As a module over $P_{W}$, $P_{W'}$ is thus isomorphic to $P_{W}[y_{k+1}, \ldots, y_{l}]$ with basis given by the monomials in $y_{k+1}, \ldots, y_{l}$. Thus $P_{W'}$ is free over $P_{W}$. $\Diamond$.\\
\ \\
\noindent Consequently, we can define stable isomorphism for chain complexes where every chain group is free. 
\begin{defn}
Let $\mathcal{C}$ be a chain complex with chain groups free over $P_{W}$ and let $\mathcal{C}'$ be similarly defined for $P_{W'}$. We will say that $\mathcal{C}$ is stably isomorphic to $\mathcal{C}'$ if there are injections $I,I': W, W' \hookrightarrow W''$ such that $\mathcal{C} \otimes_{I} P_{W''}$ is chain isomorphic to $\mathcal{C}' \otimes_{I'} P_{W''}$ as chain complexes over $P_{W''}$. Likewise, we will say that  $\mathcal{C}$ is stably chain homotopic to $\mathcal{C}'$ if there is a $W''$ where $\mathcal{C} \otimes_{I} P_{W''}$ is chain homotopic to $\mathcal{C}' \otimes_{I'} P_{W''}$
\end{defn}
\noindent Due to the flatness, if $\mathcal{C}$ is stably isomorphic to $\mathcal{C}'$ then $H_{i}(\mathcal{C})$ is stably isomorphic to $H_{i}(\mathcal{C}')$ through $(W'', I, I')$. Consequently, there is a well-defined notion of stable rank, and stably isomorphic complexes will have identical Euler characteristics.\\
\ \\
\begin{defn}
Let $v \in W$ with $W$ an $\F$-vector space. $\mathcal{K}_{W}(v)$ is the complex
$$
0 \longrightarrow P_{W} \stackrel{\cdot v}{\longrightarrow} P_{W} \longrightarrow 0
$$
supported in gradings $+1$ and $-1$. Let $v_{1}, \ldots, v_{k}$ be vectors in $W$. Then 
$$
\mathcal{K}_{W}(v_{1}, \ldots, v_{k}) = \mathcal{K}_{W}(v_{1}) \otimes_{P_{W}} \mathcal{K}_{W}(v_{2}) \otimes_{P_{W}} \cdots \otimes_{P_{W}} \mathcal{K}_{W}(v_{k})
$$ 
denotes the Koszul complex for $v_{1}, \ldots, v_{K}$. 
\end{defn}
\ \\
\noindent{\bf Example:} Let $L$ be a link diagram, and let $W$ be the vector space over $\Z/2\Z$ generated by the faces of $L$. Then $P_{W} = \ring{L}$. Then the Koszul complexes in the definition correspond to those in the definition of the totally twisted Khovanov homology in section \ref{sec:const}.\\
\ \\
\noindent The following is a straightforward exercise in definitions
\begin{prop}\label{prop:chain}
Let $I: W \hookrightarrow W'$ be an injection of vector spaces, and $v_{1}, \ldots, v_{k} \in W$. Then $\mathcal{K}(v_{1}, \ldots, v_{k})$ is stably chain isomorphic to $\mathcal{K}(I(v_{1}), \ldots, I(v_{k}))$. 
\end{prop}
\ \\
\noindent In addition, some of our chain complexes will be defined over fields with large isomorphism groups. These allow us to construct ``new'' chain complexes. let $C_{\ast}$ be a chain complex over a field $\F$, and let $\sigma: \F \rightarrow \F'$ be a field homomorphism, then we have new chain complex over $F'$: $C'_{\ast}  = C_{\ast} \otimes_{\F} \F'$ where the action of $\F$ on $\F'$ is given by $(\lambda, f) \rightarrow \sigma(\lambda)\cdot f$. In effect this construction just applied $\sigma$ to all the coefficients. That is, if we have a basis $\{x_{i}\}$ given in $C_{\ast}$ and $\partial x' = \sum f_{i} x_{i}$, then $C_{\ast}'$ will be spanned by the same basis elements but with differential map $\partial' x = \sum \sigma(f_{i}) x_{i}$. $\partial'$ is easily verified to be a differential from $\sigma$ being a field map.  In particular, when $\F = \mathrm{Rat}(W)$ and $F'=\mathrm{Rat}(W')$ with $I : W \hookrightarrow W'$, $C_{\ast}$ and $C'_{\ast}$ will be stably isomorphic. We will also denote the tensor product as $C_{\ast} \otimes_{\sigma} \F'$ when we wish to emphasize the homomorphism. 

\section{Invariance for the totally twisted Khovanov homology}\label{sec:invariance}

\noindent Given a set, $S$, and a field, $\F$, let $W_{\F,S}$ be the vector space over $\F$ generated by the elements of $S$. As $\F$ will usually be $\Z/2\Z$, or at least be clear from the context, we will often omit the field subscript when using this notation. We will write $x_{s}$ for the basis vector corresponding to $s \in S$. Let $\F[S] = \mathrm{Sym}(W_{\F,S})$ be the polynomial ring with coefficients in $\F$, and let $\F(S)$ denote the field of fractions of $\F[S]$, which we will identify with the field of rational functions in the formal variables $x_{s}$. \\
\ \\
\noindent We will use a special notation for certain elements in $\F[S]$ (and by extension $\F(S)$). For each $T \subset S$ let
$$
[T] = \sum_{s \in T} x_{s}
$$
We will also denote this element by $x_{T}$ when it is convenient. Note that for $i_{1},\ldots, i_{k} \in S$ we will shorten $[\{i_{1}i_{2}\ldots i_{k}\}]$ to $[i_{1}i_{2}\ldots i_{k}]$, so both will denote  $x_{i_{1}} + x_{i_{2}} + \ldots + x_{i_{k}}$. 

\subsection{Two auxiliary constructions}

\subsubsection{Dissection:} 
\begin{center}
\begin{figure}
\treefig{0.7}{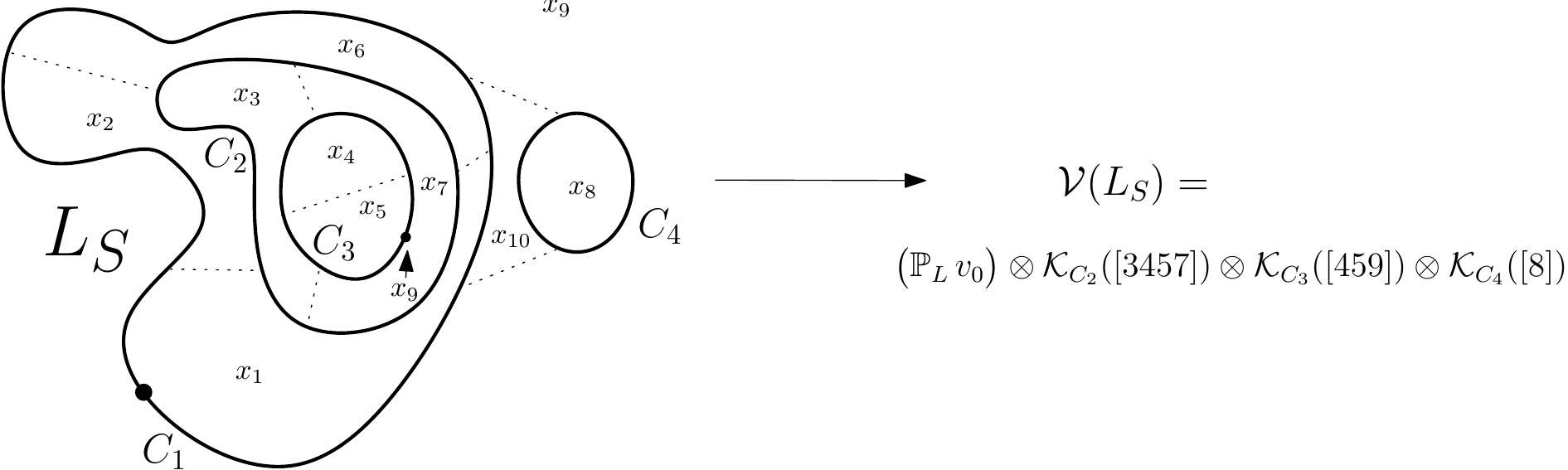} 
\caption{The complex assigned to a resolution in which some regions have been dissected and some edges have received weights.}
\label{fig:auxiliary}
\end{figure}
\end{center}

In constructing $\ring{L}$ we associated a formal variable to each face in $\regions{L}$. If we consider modules over these rings up to stable equivalence, we can make this more flexible. In particular, we can partition the faces in $L$ and assign formal variables to each of the new components. Usually we will do this by embedding arcs in $S^{2}$ with endpoints on $L$ and interiors mapped to $S^{2}\backslash \Gamma_{L}$. To make this precise, pick a face $R_{i} \in \regions{L}$ and partition it into a set of sub-components $S_{1,i}, \ldots, S_{k,i}$. Associate a formal variable $y_{j,i}$ to each $S_{j,i}$. Then we can form a complex $KT_{\ast}(L;S)$ as above, but with all the formal areas, $\area{C}$ measured in $\ring{L,S} = \Z/2\Z[x_{1}, x_{2}, \ldots, x_{i-1}, y_{1,i}, \ldots y_{k,i}, x_{i+1},\ldots, x_{m}]$. This gives a complex over $\ring{L,S} = \F[\big(\regions{L}\backslash \{R_{i}\}\big) \cup S]$, but we can define a map $I : \ring{L} \hookrightarrow \ring{L,S}$ by the decomposition relations:\\
$$
x_{i} \stackrel{I}{\longrightarrow} y_{1,i} + y_{2,i} + \ldots + y_{k,i}
$$
\ \\
\noindent corresponding to $R_{i} = S_{1,i} \cup \cdots S_{k,i}$. For an example, see Figure \ref{fig:auxiliary}.

\begin{lemma}\label{lem:decomp}
The complex $KT_{\ast}(L) \otimes_{I} \ring{L,S}$ is chain isomorphic to $KT_{\ast}(L;S)$. In particular, $HT_{\ast}(L)$ will be stably isomorphic to $HT_{\ast}(L;S)$.
\end{lemma}

\noindent{\bf Proof:} First, note that $\ring{L} \otimes_{I} \ring{L,S} \cong \ring{L,S}$ where $r \otimes y = 1 \otimes I(r)y$. Let $v$ and $w$ be generators for the chain groups of $KT_{\ast}(L)$. Then the coefficient $\langle \partial_{KT(L)} v, w \rangle \otimes 1$ will equal $\langle \partial_{KT(L;S)} v, w\rangle$. In particular, if $\langle \partial_{KT(L)} v, w \rangle = x_{R_{i}} + \sum$, where $\sum$ does not contain any $x_{R_{i}}$ summand, then $\langle \partial_{KT(L)} v, w \rangle \otimes 1 = \big(x_{R_{i}} + \sum\big) \otimes 1 = 1 \otimes (\sum_{j} y_{j,i} + \sum)$ which is $\langle \partial_{KT(L;S)} v, w\rangle$. Thus the differential on $KT_{\ast}(L) \otimes \ring{L,S}$ is identical to the differential in $KT_{\ast}(L;S)$ under the isomorphism induced by the tensor products. The result then follows. $\Diamond$

\subsubsection{Edge Weighting}

\noindent In addition to dissecting the regions, we will also find it useful to attach a weight to an edge, which can then be incorporated into the formal area of a circle in resolutions $L_{S}$. We depict this by adding a dot to the edge and labeling the dot with an element of $\ring{L}$ or with a new formal variable (and then enlarging $\ring{L}$). Diagrammatically we will draw an arrow from the weight to the point, to distinguish it from the weights assigned to the faces. For an example, see Figure \ref{fig:auxiliary}. \\
\ \\
A weight $w$ assigned to a point on an edge should be interpreted as adding $w$ to the weights of each of the regions adjacent to the edge. Thus, if $C$ is a circle in a diagram for $L$ and we add a weight $w$ to a point in an edge then the area of $C$ in the new diagram is 
\begin{enumerate}
\item $\area{C}$ if the closure of $\disc{C}$ contains neither or both of the faces adjacent to the edge,
\item $\area{C} + w$ if $\disc{C}$ contains one of the adjacent faces but not the other
\end{enumerate}
Each edge occurs in only one circle $C$ in a resolution $L_{S}$. Given $L_{S}$, write states as $v \otimes v'$ where $v'=v_{\pm}$ is the decoration on the circle containing the weighted edge. If $\partial_{\mathcal{V}(L_{S})}$ is the Koszul differential without the edge weight, then the complex with the edge weight has differential $\partial_{\mathcal{V}(L_{S})} + w\,D$ where $D(v \otimes v_{+}) = v \otimes v_{-}$ and $D(v \otimes v_{-}) = 0$. In particular. \\
\ \\ 
\noindent If there are multiple points on the same edge, we can coalesce the points into one point and add their weights to get the weight of the new point. Once we have the adjusted area we can form the twisted Khovanov complex as before. For an example, see Figure \ref{fig:auxiliary}. In particular, the area for $C_{3}$ is $x_{4} + x_{5} + x_{9}$ due to the edge weight, while that for $C_{2}$ is $x_{3}+ x_{4} + x_{5} + x_{7}$ with no edge contribution, since $C_{3}$ contains both regions which receive the extra $x_{9}$ term.

\begin{lemma}\label{lem:edgeWeight}
Let $L$ be an edge weighted diagram for a link $\mathcal{L}$, and let $e$ be an edge assigned weight $w$. Let $L'$ be the edge weighted diagram identical to $L$ away from $e$, but with weight $0$ on $e$. If $w$ is linearly independent from the other areas and weights in $L$, then the twisted complex for $L$ is stably chain isomorphic to that of $L'$.  
\end{lemma}

\noindent{\bf Proof:} Let $A$ and $B$ be the variables for the adjacent regions, then $x_{A} \lra x'_{A} + w$ and $x_{B} \lra x'_{B} + w$ will give the desired
change of variables to identify the complex for $L$ and $L'$. The linear independence guarantees that $x'_{A} + w$ is not equal to the adjusted area assigned to any other face. Indeed the adjusted areas remain linearly independent vectors over $\F$. Consequently, the map is an injection $V \rightarrow V'$ and the conclusion holds.

\subsection{Deforming the chain complex}

\subsubsection In proving that stable chain homotopy class of $CT_{\ast}(L)$ is invariant under the Redemeister moves we will make use of the following formulation of a well-known lemma in homological algebra, which follows from a graded version of Gaussian elimination.

\begin{lemma}\label{lem:deform}
Let $(M,d)$ be a differential module over a ring $R$. Suppose $M \cong M_{1} \oplus M_{2} \oplus M_{3}$ as an $R$-module and that  $d=\big[L_{ij}\big]_{i,j = 1,2,3}$ with respect to this decomposition. If $L_{32}$ is an $R$-isomorphism, then there is a sub-module $D \subset M$ with 
$d|_{D} : D \longrightarrow D \subset M$ such that
\begin{enumerate}
\item $(D,d|_{D})$ is a deformation retract of $(M,d)$, and
\item $(D,d|_{D})$ is isomorphic to $(M_{1}, d')$ where $d' = L_{11} - L_{12} \circ L_{32}^{-1} \circ L_{31}$
\end{enumerate}
\end{lemma}

\noindent We will refer to the process of cutting down from $(M,d)$ to $(D,d|_{1})$ as {\em reduction}, and we will say that we are {\em canceling} $L_{32}$. It will often be the case that $L_{11}^{2} = 0$, and so $L_{11}$ is a differential on $M_{1}$. For this reason, we will call $- L_{12} \circ L_{32}^{-1} \circ L_{31}$ a {\em perturbation} term and $d'$ the perturbed differential. As the lemma is well known, we will omit some of the computations underlying the proof in favor of indicating which computations should be performed. \\
\ \\
\noindent{\bf Proof:} Let $Q = M_{2} + d(M_{2})$. Then $(Q, d|_{Q})$ is a sub-complex of $M$. The quotient complex $M/Q$ can be described through $M/Q \cong (M/M_{2})\big/(Q/M_{2}) = (M_{1} \oplus M_{3})\big/(d(M_{2})/(d(M_{2}) \cap M_{2})$. Since $L_{32}$ is an isomorphism, quotienting by $d(M_{2})$ results in a module isomorphic to $M_{1}$. The differential on this module can be computed by examination $d(x \oplus 0 \oplus 0) = L_{11}(x) \oplus L_{21}(x) \oplus L_{31}(x)$ $= L_{11}(x) \oplus L_{31}(x)$ modulo $M_{2}$. $L_{31}(x) = L_{32}(y)$, and, modulo $d(M_{2})$, $L_{32}(y) + L_{12}(y) \equiv 0$. Consequently, 
under $\pi: M \rightarrow M/Q$, $d(x)$ is mapped to $L_{11}(x) - L_{12}(y) = (L_{11} - L_{12}L_{32}^{-1}L_{31})(x)$. Thus the quotient $M/Q$ is the complex in item ii) of the lemma. In particular, this shows that $d'$ is a differential, which can also be verified directly using the nine relations between the $L_{ij}$ found from $d^{2} \equiv 0$. \\
\ \\
\noindent On the other hand, there is a chain map $(M_{1}, d') \rightarrow (M,d)$ defined as $\iota(x) = x \oplus -L_{32}^{-1}L_{31}(x) \oplus 0$. That this is a chain map is an exercise in using the entries of $d^{2} \equiv 0$. Furthermore, $\iota$ is injective, due to the form of the first summand. We let $D$ be $\mathrm{im}\,\iota$, so $(D, d|_{D})$ is chain isomorphic to $(M_{1}, d')$. It remains to verify that $(D, d|_{D})$ is a deformation retract. However,
$$
\pi \circ \iota = \mathrm{Id}_{M_{1}}
$$
and if we let 
$$
H = \left[\begin{array}{ccc} 0 & 0 & 0 \\ 0 & 0 & -L_{32}^{-1} \\ 0 & 0 & 0 \end{array}\right]
$$
it is then easy to verify (once again using $d^{2} = 0$) that 
$$
\iota \circ \pi - \mathrm{Id}_{M} = d\,H + H\,d
$$
so that $(D, d|_{D})$ is a deformation retract of $(M,d)$. In particular, $D$ is chain homotopy equivalent to $M$. $\Diamond$\\
 
\subsection{Invariance under the first Reidemeister move}\ \\

\noindent{\bf Convention:} All gradings in the following sections are $\widetilde{\delta}$ gradings.

\begin{prop}
Let $c$ be a crossing in an oriented link diagram, $L$, which can be removed by a local Reidemeister I move. If $L'$ is the diagram after the move, then $KT_{\ast}(L)$ is (stably) chain homotopy equivalent to $KT_{\ast}(L')$.
\end{prop}

\begin{center}
\begin{figure}
\treefig{0.7}{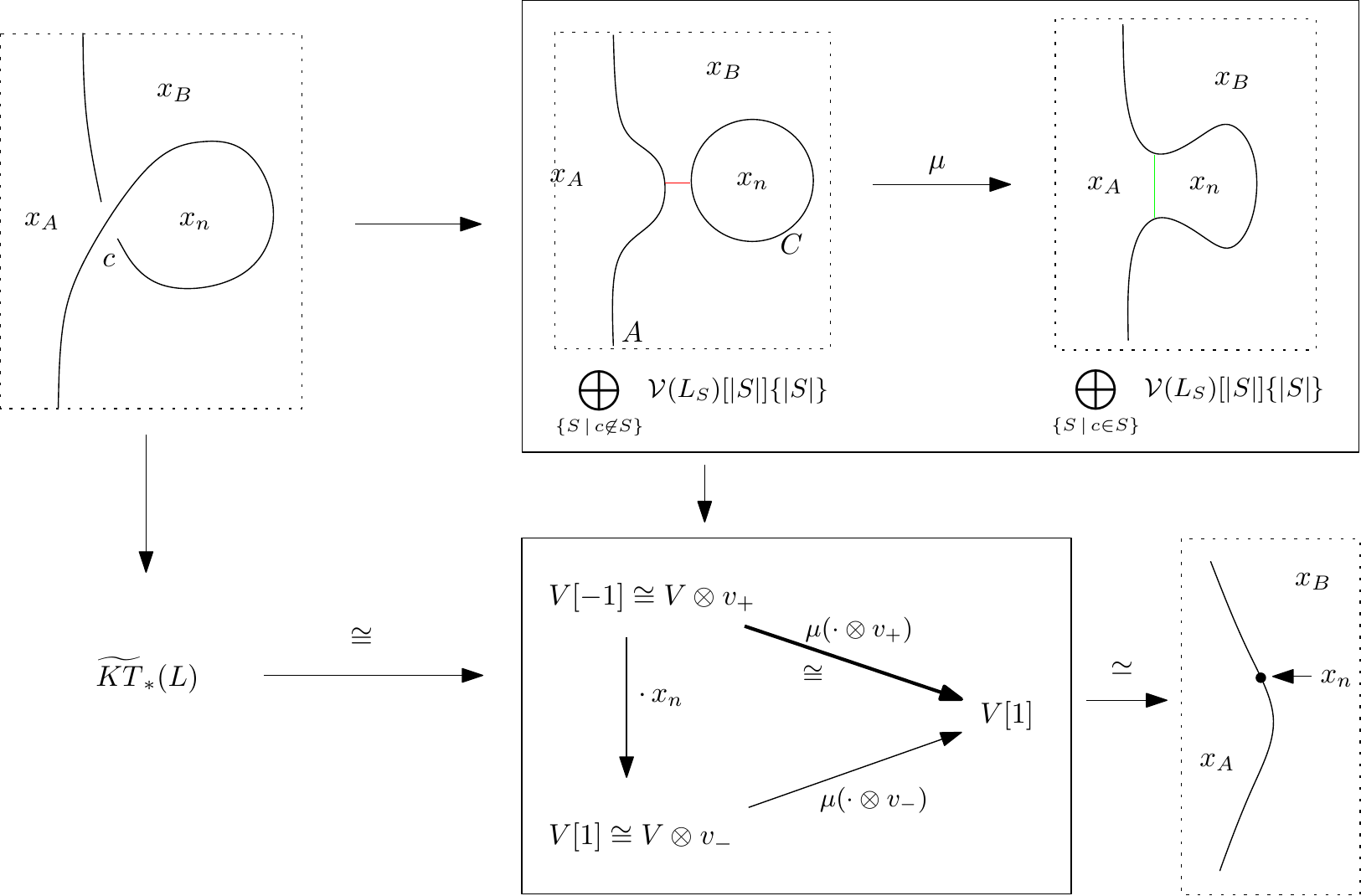} 
\caption{A schematic representation of the proof of invariance under the first Reidemeister move, for a positive crossing. The complex is reduced along the thickened arrow, which is an isomorphism of the respective sub-modules of the chain modules. The proof in the text additionally shows that the final diagram represents a complex stably isomorphic to the diagram without the edge weighting. }
\label{fig:RI}
\end{figure}
\end{center}

\noindent {\bf Proof:} There are two cases to consider, based on the handedness of the crossing. Let $L'$ be the diagram resulting after the Reidemeister move. The argument for Case I below is directly reflected in the diagrams in Figure \ref{fig:RI}.\\
\ \\
\noindent {\bf Case I: $c$ is right handed}.  We can divide $\widetilde{KT}_{\ast}(L)$ as a direct sum of $\bigoplus_{c \not\in S} \mathcal{V}(L_{S})$ and $\bigoplus_{c \in S} \mathcal{V}(L_{S})$. When $c \not\in S$, there is
a complete circle, $C$, in the local diagram used in the Reidemeister move. Thus $\mathcal{V}(L_{S}) \cong \mathcal{V}(L') \otimes \mathcal{K}(C)$. Let $x_{n}$ be the formal variable associated with the region, $\disc{C}$, so that $\partial_{C}$ is multiplication by $x_{n}$. The complex, $\widetilde{KT}_{\ast}(L)$ can be decomposed further. If $V = \mathcal{V}(L')$, then $\bigoplus_{c \in S} \mathcal{V}(L_{S})$ $\cong V[1] $, where the grading shift comes from the additional resolution at $c$ when compared with $L'$, and $\bigoplus_{c \not\in S} \mathcal{V}(L_{S}) \cong V_{+} \oplus V_{-}$ where  $V_{\pm} = V \otimes v_{\pm}$ and the last factor is that from $C$. Then $\widetilde{KT}_{\ast}(L) \cong V_{+} \oplus V_{-} \oplus V[1] \cong V[-1] \oplus V[1] \oplus V[1]$ since $V_{+}$ has Khovanov bigrading shifted by $(0,1)$ compared with $V$ and that corresponds to a $\widetilde{\delta}$-shift of $-1$. The twisted differential, when written to respect this decomposition, is the sum of two maps:
$$
\partial_{L} = \left[\begin{array}{ccc} \partial_{L'} & 0 & 0 \\ \cdot x_{n} & \partial_{L'} & 0 \\ 0 & 0 & \partial_{L'} \end{array}\right] \hspace{1in}
\partial_{KH,L} = \left[\begin{array}{ccc} \partial_{KH,L'} & 0 & 0 \\ 0 & \partial_{KH, L'} & 0 \\ \mu(\cdot\otimes_{+}) & \mu(\cdot \otimes v_{-}) & \partial_{KH, L'} \end{array}\right]
$$
In particular, $\mu(\cdot \otimes v_{+}) : V_{+}\cong V[-1] \stackrel{\mu}{\longrightarrow} V[1]$ is an isomorphism, since $v_{+}$ acts as the identity element in the Frobenius algebra. The inverse map $H$ is simply $\xi \longrightarrow \xi \otimes v_{+}$\\
\ \\
\noindent Consequently, we can cancel the map $\mu(\cdot \otimes v_{+})$ and obtain a deformation retraction of the complex, $KT_{\ast}(L)$ supported on $V_{-}$, with a new differential, $\partial'$. To compute this differential there are two cases to consider:
\begin{enumerate}
\item {\bf Case i:} The arc in the local diagram has been assigned a $v_{-}$: Such states are in the kernel of $\mu(\cdot \otimes v_{-})$. Consequently there will be no change in the differential. Thus $\partial' = \partial_{V}$ on such states. 
\item {\bf Case ii:} The arc in the local diagram has been assigned a $v_{+}$: We will write the generator as $v \otimes v_{+} \otimes v_{-}$ where $v$ incorporates the assignments for the circles not in the local diagram, $v_{+}$ is the marking on the arc, and $v_{-}$ is the marking on $C$. $\partial_{KH}(v \otimes v_{+} \otimes v_{-})$ equals $v \otimes \mu(v_{+} \otimes v_{-})$ $= v \otimes v_{-}'$ with the $v_{-}'$ assigned to the arc. However, any state, $\xi$ in $V[1]$ with a $v_{-}'$ on the local arc, is canceled by $H(\xi)$ in the complex for $L$. Thus in the deformed complex, the boundary of $v \otimes v_{+} \otimes v_{-}$ will be $\partial_{KT, V} \xi + \partial_{C} H(v \otimes v_{-}') = \partial_{L'} + \partial_{C} (v \otimes v_{-} \otimes v_{+})$ $= \partial  + x_{n} (v \otimes v_{-} \otimes v_{-})$. 
\end{enumerate}
Under the identification $V_{-} \cong V[1]$, we can drop the last $v_{-}$ factor. Given a full resolution, $L'_{S}$, let $C'$ be the circle containing the local arc and let $D_{C'}$ be the map on $V$ which takes a $+$ marker on $C'$ to $-$-marker, fixing the marking on the other circles, and a state with a $-$-marker on $V'$ to $0$. Then the perturbed differential on $V_{-}$ is $\partial_{V} + x_{n}D_{C'}$. This is the edge weighted differential occurring in the last diagram in \ref{fig:RI}. By lemma \ref{lem:edgeWeight} the last complex is stably chain isomorphic to $V _{-} = \widetilde{KT}_{\ast}(L')[1]$.\ \\
\ \\
\noindent  Finally, we address the grading shifts. We know that $n_{+}(L) = n_{+}(L') + 1$. Thus $KT_{\ast}(L) \cong \widetilde{KT}_{\ast}(L)[-n_{+}(L)]$ is stably chain homotopic to $\widetilde{KT}_{\ast}(L')[1][-n_{+}(L')-1] \cong \widetilde{KT}_{\ast}(L')[-n_{+}(L')]$. The latter is $KT_{\ast}(L')$ by definition,  and we have verified (stable) chain homotopy invariance in this case.\\
\ \\
{\bf Note:} There's a special case implicitly handled in the above argument. Namely, if the arc in the local diagram belongs to the marked circle. This is addressed by noting that if the circle $C'$ is assigned $v_{-}$, then by case i) deforming the differential has no effect on the image of that state. 
\\
\ \\  
{\bf Case II: $c$ is a left handed crossing}. We decompose $KT_{\ast}(L)$ as in the right handed case, although the gradings and differential are different. Indeed,  $\bigoplus_{c \in S} \mathcal{V}(L_{S}) \cong V_{+}[1] \oplus V_{-}[1] \cong V \oplus V[2]$ and
the map $\Delta: V \rightarrow V_{+}[1] \oplus V_{-}[1]$ followed by projection  onto $V_{-}[1]$ is surjective. There is still a map $V_{+}[1] \stackrel{\cdot x_{n}}{\longrightarrow} V_{-}[1]$ as before, and elements in $V$ with non-trivial image under the composition of $\Delta$ with projection to $V_{+}[1]$. These occur when the circle $C'$ containing the local arc in $V$ is adorned with a $v_{+}$. Thus, we may deform the complex to one supported on $V_{+}[1]$ with a perturbed differential. Namely, for states marked $v_{+} \otimes v_{+}$ (arc $\otimes$ $C$) in $V_{+}[1]$ the image under $\partial_{C}$ is
$x_{n}(v_{+}\otimes v_{-})$. This is canceled using $\Delta$ since $\Delta(x_{n} v_{+}) = x_{n}(v_{-} \otimes v_{+}) \oplus x_{n}(v_{+} \otimes v_{-})$. But this map also has the component $x_{n}(v_{-}\otimes v_{+})$ in $V_{+}[1]$. Consequently, the effect of the cancellation is to perturb the differential in $V_{+}[1]$ by adding terms $v_{+} \otimes v_{+} \longrightarrow x_{n}(v_{-} \otimes v_{+})$. Once again, this has the effect of adding $x_{n}$ to the  formal area of the region bounded by $C'$. The result of this deformation is a chain complex in the same
$\delta$-gradings as $\big(V \otimes v_{+}\big)[1]$. This complex is isomorphic, under the change of variables imposed by lemma \ref{lem:decomp}, to $V$, with the same grading since the tensor product with $v_{+}$ contributes a $-1$ to the $\delta$-grading. On the other hand, the number of positive crossings is not altered by the Reidemeister move, hence $KT_{\ast}(L)$ is (stably) chain homotopy equivalent to $KT_{\ast}(L')$.  $\Diamond$

\subsection{Invariance under the second Reidemeister move}

\begin{center}
\begin{figure}
\treefig{0.7}{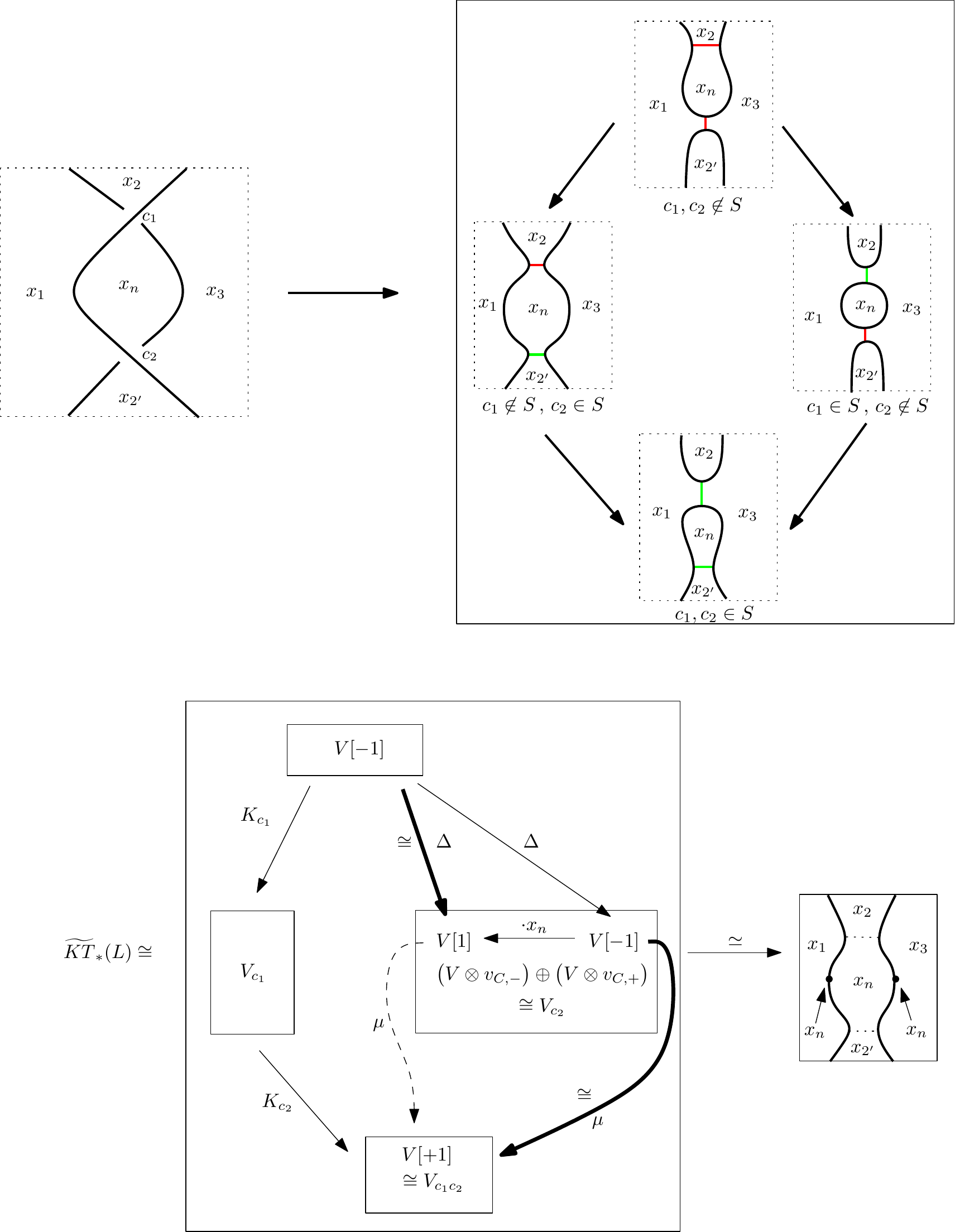} 
\caption{A schematic representation of the proof of invariance under the second Reidemeister move. The complex is reduced along the thickened arrows, which are both isomorphisms. This is essentially the proof in \cite{Khov}, but the twisting differential now alters the differential on the complex after the reduction. The alteration is depicted through the edge weighting in the last diagram. It is straightforward to see that the weights are linearly independent of the areas, and thus that the last diagram is stably isomorphic to the diagram for $L'$.}
\label{fig:RII}
\end{figure}
\end{center}

\begin{prop}
If $L$ is a diagram for $\lk{L}$ and $L'$ is another diagram differing from $L$ only by a local Redemeister II move, then $KT_{\ast}(L)$ is (stably) homotopy equivalent to $KT_{\ast}(L')$. 
\end{prop}

\noindent{\bf Proof:} The argument refers to Figure \ref{fig:RII}. We will use the notation $V_{c_{1}} = \bigoplus_{c_{1} \not\in S, c_{2} \in S} \mathcal{V}(L_{S})$. Note that $V_{c_{1}}$ is identical with $L'$, up to the area contributions. The proof consists of canceling portions of $\widetilde{KT}_{\ast}(L)$ using lemma \ref{lem:deform}, to obtain a perturbed differential on $V_{c_{1}}$ which can then be interpreted as a complex for $L'$ with different areas assigned to each region. To that end, note that as a module $\widetilde{KT}_{\ast}(L) \cong V_{\emptyset} \oplus V_{c_{1}} \oplus V_{c_{2}} \oplus V_{c_{1}c_{2}}$. \\
\ \\
\noindent  Furthermore, $V_{c_{2}}$ can be decomposed as $\big(V \otimes v_{+}\big) \oplus \big(V \otimes v_{-}\big)$, where the second factor is the decoration on the circle, $C$, in the local diagram for $V_{c_{2}}$ (in Figure \ref{fig:RII}, and  $V \cong V_{c_{1}c_{2}}$.  The map $V_{c_{2}} \longrightarrow V_{c_{1}c_{2}}[1]$ in $\widetilde{KT}_{\ast}(L)$ is identical to the map in the Khovanov complex: it is $\mu$ applied to the last factor and a factor corresponding to the arc $C$ merges into. In particular, on $V \otimes v_{+}$, $\mu$ restricts to an isomorphism $f: V \otimes v_{+} \longrightarrow V_{c_{1}c_{2}}[1]$. Consequently, given any state $s$ in $V_{c_{1}c_{2}}[1]$, there is a canceling state, $f^{-1}(s)$ in $V_{c_{2}}$. Likewise, the Khovanov division map $\Delta: V[-1] \cong V_{\emptyset} \longrightarrow V_{c_{2}}$ followed by projection to $V \otimes v_{-}$ is an isomorphism, since $s \longrightarrow L(s) \oplus (s \otimes v_{-}) $ where $L(s) \in V \otimes v_{+}$ may or may not be zero, depending on the state $s$. Call this isomorphism $g$. \\
\ \\
\noindent The twisted differential provides a map $\partial_{C}: V \otimes v_{+} \stackrel{\cdot x_{n}}{\longrightarrow} V \otimes v_{-}$ which takes $s \otimes v_{+}$ to $x_{n}(s \otimes v_{-})$ (there can, of course, be other non-zero summands in the twisting differential, but they will not be relevant here). \\
\ \\
\noindent We can cancel the isomorphisms $f$ and $g$ simultaneously, since their domains and images are disjoint. In particular, they can be considered two portions of a single $+2$ map formed from summands in the differential for $\widetilde{KT}_{\ast}(L)$. Doing so produces a deformation equivalent chain complex which can be described as $V_{c_{1}}$ with a perturbed differential. We now endeavor to compute this perturbation.  \\
\ \\
\noindent We start with a state, $s$, for $V_{c_{1}}$. To compute the perturbation applied to $s$, we first apply the Khovanov map $K_{c_{2}} : V_{c_{1}} \longrightarrow V_{c_{1}c_{2}}[1]$ to obtain new state data in $V_{c_{1}c_{2}}$. $K$ is either a copy of the Frobenius algebra maps $\mu$ or $\Delta$, but will depend on the resolution underlying the state $s$. The image $K(s)$ could be zero, but in any case, this new data is the image of $K_{c_{1}}(s) \otimes v_{+}$ under $f$. Thus canceling $f$ deforms the chain complex to one supported on the direct sum of $V_{c_{1}}$, $V \otimes v_{-}$ and $V$ with differential given by taking the original differential projected to these summands and adding $\partial_{V_{c_{2}}}(K_{c_{2}}(s) \otimes v_{+})$. However, the terms in $\partial_{V_{c_{2}}}(K_{c_{2}}(s) \otimes v_{+})$ all lie in $V \otimes v_{+}$, and thus will be canceled when we cancel $g$, except for $x_{n}\big(K_{c_{2}}(s) \otimes v_{-}\big)$. The new differential on $V_{c_{1}}$ is thus $\partial_{V_{c_{1}}}(s) \oplus x_{n}\big(K_{c_{2}}(s) \otimes v_{-}\big) \oplus 0$. Canceling the map $g : V[-1] \longrightarrow V \otimes v_{-}$ deforms this complex to one supported on $V_{c_{1}}$. The new differential arises from the fact that $x_{n}\big(K_{c_{2}}(s) \otimes v_{-}\big)$ lies in the image of $g$. Using the form for $g$ given above, we have $g^{-1}(x_{n}\big(K_{c_{2}}(s) \otimes v_{-}\big)$ $=x_{n} K_{c_{2}}(s)$  as a state in $V$.
\ \\
\noindent The portion of the differential for $\widehat{KT}_{\ast}(L)$ mapping $V[-1] \longrightarrow V_{c_{1}}$ is the Khovanov map $K_{c_{1}}$ coming from changing the resolution at $c_{1}$ in a diagram contributing to $V$. Consequently, the perturbation of the differential on $V_{c_{1}}$ after canceling $g$ is $x_{n} K_{c_{1}} \circ K_{c_{2}}$; the perturbed differential is $\partial_{V_{c_{1}}} + x_{n} K_{c_{1}}K_{c_{2}}$. We will now examine the maps $K_{c_{1}}$ and $K_{c_{2}}$ more closely.  
\ \\
\noindent For any resolution $L_{S}$ with $c_{1} \not\in S$, $c_{2} \in S$, we will let $C_{L}$ and $C_{R}$ be the circles containing the left and right arcs in the local picture for $V_{c_{1}}$ in \ref{fig:RII}. Note that $C_{L}$ and $C_{R}$ {\em may be the same circle in the larger diagram} $L_{S}$. The map $K_{c_{1}}K_{c_{2}}$ is just a Khovanov homology map, and can be computed by considering two cases.   \\
\ \\
\noindent{\bf Case i: $C_{L} = C_{R}$} Then $K_{c_{2}}$ is a copy of $\Delta$ and $K_{c_{1}}$ is a copy of $\mu$ applied to the circles resulting from the division at $c_{2}$ that gives $\Delta$. But in characteristic $2$, $\mu \circ \Delta : \mathcal{V} \rightarrow \mathcal{V}$ is the zero map. Consequently there is no perturbation to the differential applied to any state from this case. \\
\ \\
\noindent{\bf Case ii: $C_{L} \neq C_{R}$} Then $K_{c_{2}}$ is a merge map, and we use $\mu$. $K_{c_{1}}$ comes from dividing this same circle, and uses $\Delta$. Thus we need to compute $\Delta \circ \mu : \mathcal{V}_{C_{L}} \otimes \mathcal{V}_{C_{R}} \longrightarrow \mathcal{V}_{C_{L}} \otimes \mathcal{V}_{C_{R}}$. It is straightforward to verify that 
$$
\begin{array}{l} 
	v_{+} \otimes v_{+} \longrightarrow v_{-} \otimes v_{+} + v_{-} \otimes v_{+}  \\
	v_{+} \otimes v_{-}, v_{+} \otimes v_{-} \longrightarrow v_{-} \otimes v_{-} \\
	v_{-} \otimes v_{-} \longrightarrow 0 \\
	\end{array}
$$
The perturbation is $x_{n}$ times this map, and thus equals $x_{n}\big(D_{C_{L}} \otimes \mathbb{I} + \mathbb{I} \otimes D_{C_{R}}\big)$, where $D_{C}$ is the isomorphism $V \otimes v_{+, C} \rightarrow V \otimes v_{-,C}$ where $\mathcal{V}_{C} = \F\,v_{+,C} \oplus \F\,v_{-,C}$.  We can interpret this formula as endowing each of the two arcs in the local diagram with the additional weight $x_{n}$, and including that weight in the formal area used in the vertical differentials for $C_{L}$ and $C_{R}$. When $C_{L}$ and $C_{R}$ coincide, the weight is added twice, once for each arc, and thus cancels so the vertical differential doesn't change. When $C_{L}$ and $C_{R}$ are distinct,  $x_{n}$ is added to each of $\area{C_{L}}$ and $\area{C_{R}}$. Note also that if one or both arcs is contained in the marked circle, it will act as if adorned by $v_{-}$ and the additional contribution will not appear for that circle. \\
\ \\
\noindent To complete the proof, we will now show that this deformed complex is (stably) isomorphic to the complex for $L'$. This requires several steps. Due to the final form of the perturbation, the Khovanov portions of the differential are the same for the deformed complex and for $L'$. Likewise, all the (marked) resolutions are the same. We need only compare the formal areas of the discs appearing in each diagram for $L'$ with the edge weight adjusted areas from the deformed complex.
\begin{enumerate}
\item Let $Q$ be a circle in some resolution diagram $L'_{S}$, bounding a disc $D$ which completely covers the local diagram for the Redemeister move. In the deformed complex, its vertical differential is unchanged, as it does not include either arc. However, its formal area is that from the diagram for $L'$ plus $x_{n}$, since in the resolutions for $V_{c_{1}}$ there is an additional region corresponding to the crossings.  
\item Now suppose $Q = C_{L} = C_{R}$. $\disc{Q}$ can either contain the additional region $R_{n}$ or not. If $R_{n} \subset \disc{Q}$ then, in the deformed complex, we add $x_{n}$ to its area from $L'$ once due to the included region (this occurs even before the deformation). The perturbation also adds $x_{n}$ for each of the arcs in the boundary of $Q$, thus adding $x_{n}$ to the formal area {\em three} times. Modulo $2$ this is just adding it once. 
\item On the other hand, if $R_{n} \not\subset \disc{Q}$ then we do not need the first $x_{n}$ summand and the perturbation, as before, adds it twice. Hence the vertical differential in the deformed complex is identical with that from $L'$.
\end{enumerate}
There are several more cases to check, but we pause here to give an idea of where we are heading. Let the regions for $L'$ be numbered the same as those for $L$, but using $y_{1}, \ldots, y_{n-1}$ for the formal variables. Let $\overline{x} = x_{2} + x_{2}'$. The map $y_{1} \rightarrow x_{1} + x_{n}$, $y_{2} \rightarrow \overline{x} + x_{n}$ and $y_{3} \rightarrow x_{3} + x_{n}$ provides a change of variables which mimic the changes in areas described above. Circles which fully contain the local diagram will have area $y_{1} + y_{2} + y_{3} + A$ in $L'$, with $A$ denoting the area outside the local diagram, but $x_{1} + \overline{x} + x_{3} + x_{n} + A$ in the deformed complex.  Under the change of variables $y_{1} + y_{2} + y_{3} + A$ is sent to $(x_{1} + x_{n}) + (\overline{x} + x_{n}) + (x_{3} + x_{n}) + A$. Since there are an odd number of $x_{n}$'s this is the same as adding $x_{n}$ to the area. A circle which contains both arcs and the region $R_{n}$ has area $y_{2} + A$ in $L'$. This is mapped to $\overline{x} + x_{n} + A$ under the change of variables, which is exactly the area from the deformed complex. A circle which contains both arcs, but not $R_{n}$, will have area $y_{1} + y_{3} + A$, with no $y_{2}$ contribution. Under the map, this area becomes $x_{1} + x_{n} + x_{3} + x_{n} + A$. The even number of $x_{n}$'s means that the image of the area is unchanged. \\
\ \\
\noindent However, we also need to check that this continues to happen when $C_{L}$ and $C_{R}$ are distinct in $L_{S}'$. We must also use the same map from $y_{i}$'s to $x_{i}$'s as above. 
\begin{enumerate}
\item $C_{R} \not\subset \disc{C_{L}}$ and $C_{L} \not\subset \disc{C_{R}}$: Then $R_{1} \subset \disc{C_{L}}$ and $R_{3} \subset \disc{C_{R}}$. In the deformed complex $\disc{C_{R}}$ and $\disc{C_{L}}$ have area $x_{n}$ bigger than that in the complex for $L'$. The change of variables does just this: $y_{1} \rightarrow x_{1} + x_{n}$ and $y_{3} \rightarrow x_{3} + x_{n}$.
\item $C_{R} \subset \disc{C_{L}}$ (or vice-versa). Then $\disc{C_{R}} \subset \disc{C_{L}}$ as well, for otherwise their union is $S^{2}$ and contains the marked point. Furthermore, $R_{n} \subset \disc{C_{L}}$. In the complex for $L'$, $\disc{C_{L}}$ has area $y_{2} + y_{3} + A$, which is mapped to $\overline{x} + x_{3} + A$ under the change of variables. This is the effective area in the deformed complex. Indeed, in the the complex for $V_{c_{1}}$, $\disc{C_{L}}$ has area $(\overline{x} + x_{n}) + x_{3} + A$. Since $C_{L}$ contains one local arc, the deformation $\disc{C_{L}}$ adds $x_{n}$ to this area: $x_{n} + (\overline{x} + x_{n} + x_{3} + A)$ $= x_{2} + x_{3} + A$. For $\disc{C_{R}}$ the argument is the same as in the previous case since it does not contain $R_{n}$. 
\end{enumerate}
Thus the change of variables is a chain isomorphism from the complex for $L'$ to the deformation of the complex for $L$. Consequently, the two chain complexes are
stably homotopy equivalent. $\Diamond$\\
\ \\

\begin{center}
\begin{figure}
\treefig{0.7}{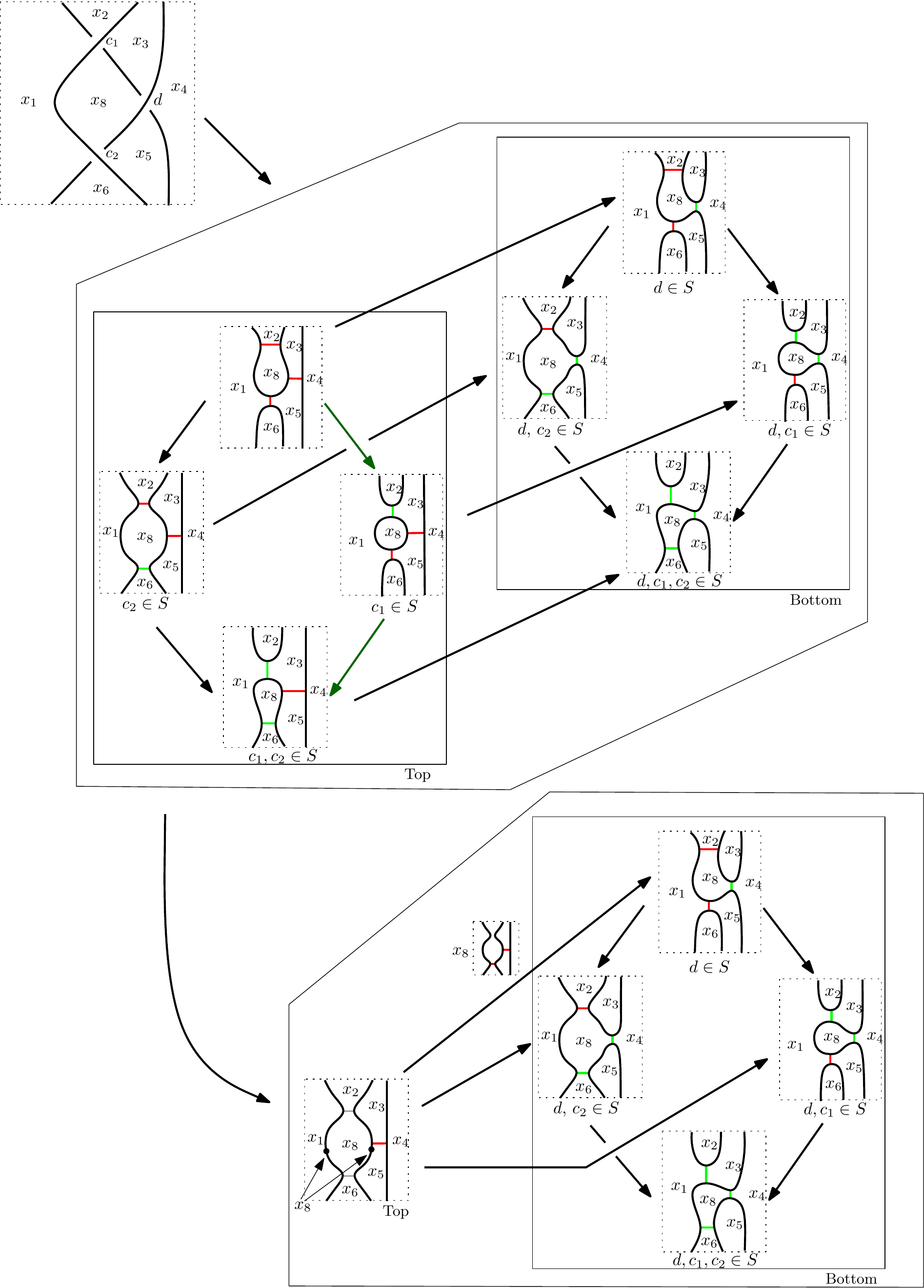} 
\caption{Diagrams for link $L$ in the proof of RIII invariance. Notice the small diagram over an arrow in the bottom picture. This depicts the surface used in 
the proof of RIII invariance.}
\label{fig:RIII_R}
\end{figure}
\end{center}

\begin{center}
\begin{figure}
\treefig{0.7}{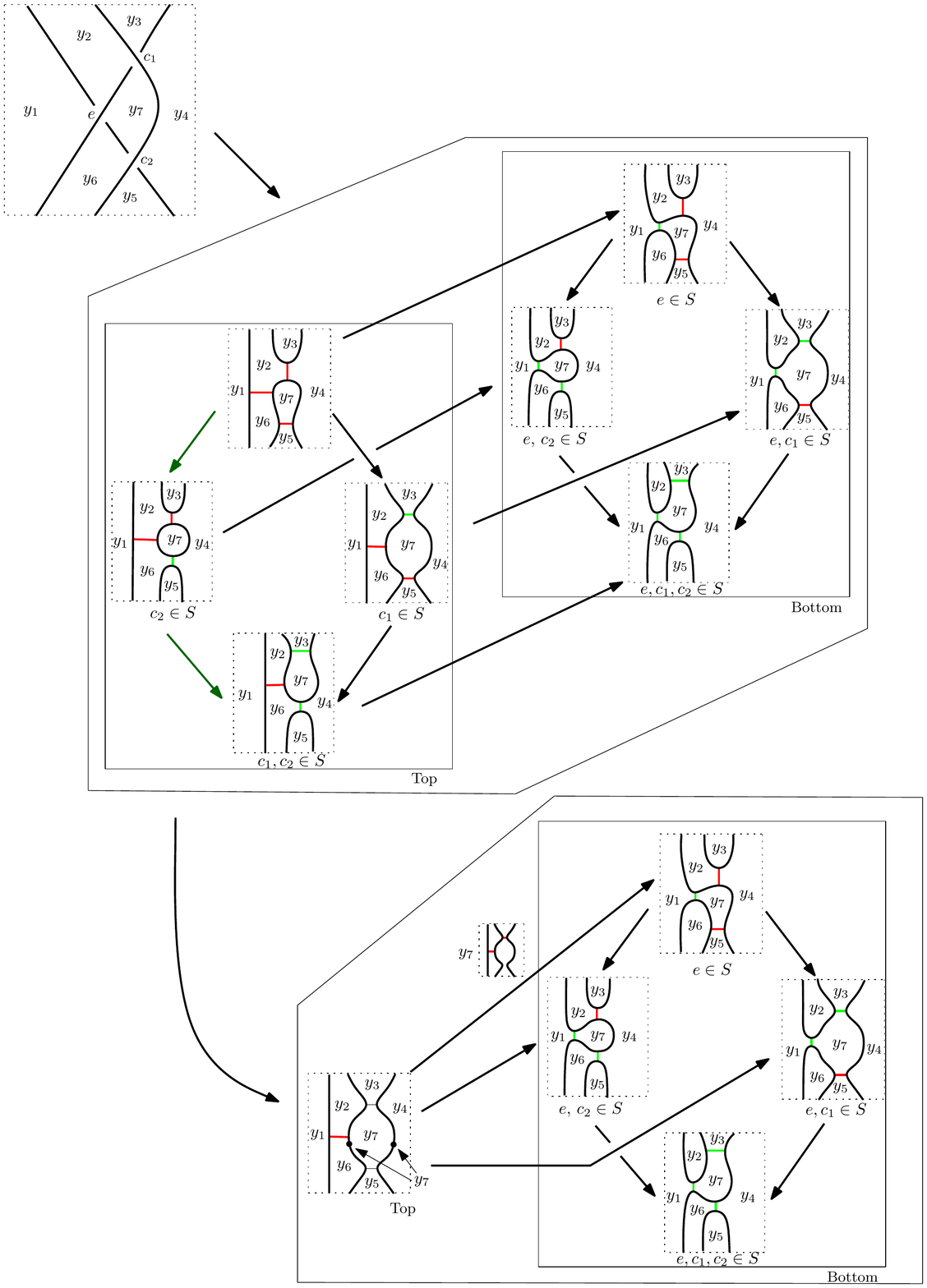} 
\caption{Diagrams for the link $L'$ in the proof of RIII invariance.}
\label{fig:RIII_L}
\end{figure}
\end{center}

\subsection{Invariance under the third Reidemeister move}

\begin{center}
\begin{figure}
\treefig{0.7}{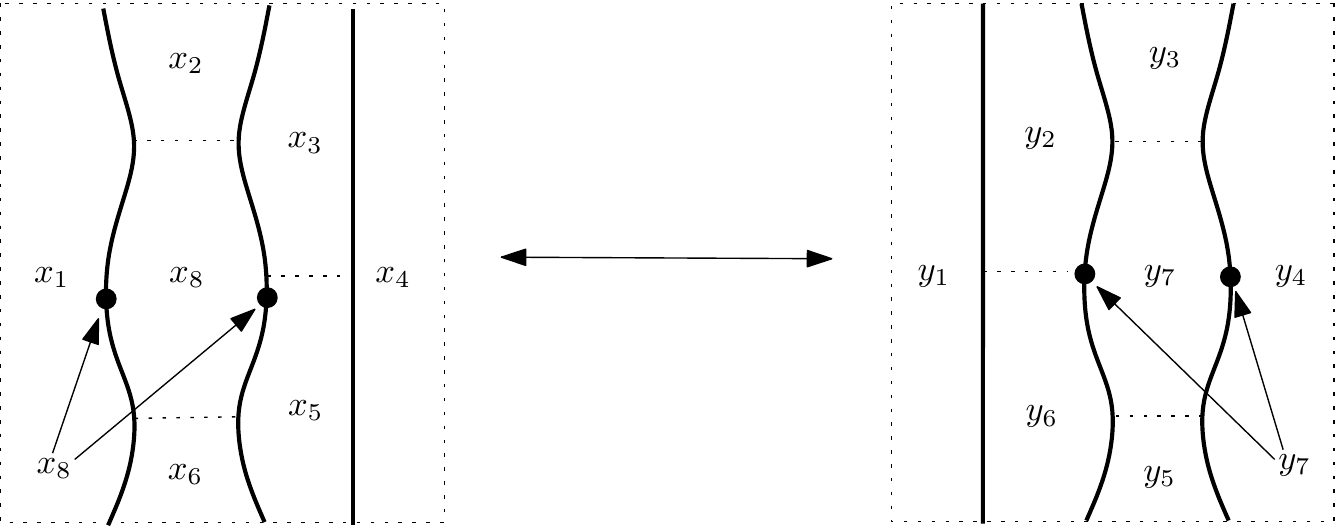} 
\caption{}
\label{fig:comp}
\end{figure}
\end{center}

\noindent We consider the case, as in Figures \ref{fig:RIII_R} and \ref{fig:RIII_L}. The resolutions, $S$, with $d,e \not\in S$ are in the upper layer of the cubes while those with $d,e \in S$ are in the bottom layer. As with the proof of $RII$-invariance, we group according to the local resolution pattern; thus $$\widetilde{KT}_{\ast}(L) \cong \mathrm{MC}\big(\bigoplus_{d \not\in S} \mathcal{V}(L_{S}) \rightarrow \bigoplus_{d \in S} \mathcal{V}(L_{S})\big)$$ with the homomorphism coming entirely from Khovanov maps. A similar result holds for $L'$ and $e$. \\
\ \\
\noindent For the two diagrams related by the RIII move, the bottom layers, when $e \in S$, $d \in S$, are identical except for the presence of $x_{8}$ in one and $y_{7}$ in the other. Otherwise, $x_{i}$ and $y_{i}$
represent the same region. However, the region corresponding to $x_{8}$ is always included in the same component as that for $x_{4}$, and the region
corresponding to $y_{7}$ is always in the same component as $y_{1}$ in their respective diagrams. Thus, defining $\Phi$ by $x_{i} \rightarrow y_{i}$, $i \neq 1,4,7,8$,
$x_{8} \rightarrow y_{7}$ and $x_{1} \rightarrow y_{1}+y_{7}$, $x_{4} \rightarrow y_{4} + y_{7}$ will be a chain isomorphism from the subcomplex $e \in S$ to that with $d \in S$.\\
\ \\
\noindent As it stands, this is not a chain isomorphism on the summands where $d \not\in S$, nor will it correctly map the connecting homorphism in the mapping cone for $L$ to that for $L'$. However, $\bigoplus_{d \not\in S} \mathcal{V}(L_{S})$ is isomorphic to the complex for a diagram with $d$ resolved so that a local $RII$-move can be performed. A similar observation  holds for $e$ and $L'$. Consequently, repeating the cancellations performed in the proof of $RII$-invariance will simplify the top layer (as is done in \cite{Bar2}) to a deformed complex $\mathcal{D}_{L}$. As the cancellations in the proof of $RII$-invariance are along the Khovanov maps, we obtain two maps from the deformed complex to $\mathcal{V}_{d,c_{1}}$ and $\mathcal{V}_{d,c_{2}}$. These maps are exactly the Khovanov maps obtained from the argument for $RIII$-invariance for Khovanov homology, and are known to be equal to the corresponding maps found from simplifying the diagrams for $L'$, compare \cite{Bar2}. Thus, if we can use the same change of variables in the new upper layer as the one we have used in the lower layer, these maps will remain equal.\\
\ \\
\noindent However, after the simplifications of the upper layer, we obtain a complex $\mathcal{D}_{L}$ in the upper layer which is the perturbed complex from the proof of RII invariance, and a new map to $\mathcal{V}_{d}$ (for Figure \ref{fig:RIII_R}). This map comes from states, $s$, in $\mathcal{V}_{c_{1}}$ mapped to $\mathcal{V}_{c_{1}c_{2}}$ and then canceled. The canceling state in $\mathcal{V}_{c_{2}}$ is $s \otimes v_{+}$, where the $v_{+}$ is for the local circle. This state has a non-zero differential to $s \otimes v_{-}$ which is canceled by the state $s$ in $\mathcal{V}$. However, $s$ in $\mathcal{V}$ may now have non-zero boundary in $\mathcal{V}_{d}$. Thus we obtain a map $s \rightarrow s'$ where $s' \in \mathcal{V}_{d}$. We can likewise simplify the complex for $L'$.\\
\ \\
\noindent For $L$, the map to $D \rightarrow \mathcal{V}_{d}$ is $x_{8}$ times the Khovanov map for a surface. This is depicted in \ref{fig:RIII_R} by the small diagram over an arrow in the bottom picture. The surface is the one obtained by attaching one handles to the black arc over the red arcs. For $L'$ the surface is likewise depicted in Figure \ref{fig:RIII_L}. In this case we obtain $y_{7}$ times the Khovanov map for the surface in the bottom picture Figure \ref{fig:RIII_L}, again depicted as a diagram over an arrow. \\
\ \\
\noindent These surfaces are evidently isotopic, preserving boundaries, since the attached one handles occur in isotopic positions; thus, the Khovanov maps are the same, \cite{Bar2}. Under the change of variables $\Phi$ above $x_{8} \rightarrow y_{7}$. Thus the three maps composing the connecting homomorphism of the mapping cone will be identified: two $\mathcal{D}_{L} \rightarrow \mathcal{V}_{d,c_{1}}, \mathcal{V}_{d, c_{2}}$ will be identified with $\mathcal{D}_{L'} \rightarrow \mathcal{V}_{e,c_{1}}, \mathcal{V}_{e,c_{2}}$ because they are Khovanov maps which do not depend on the formal variables, and the last the map $\mathcal{D}_{L} \rightarrow \mathcal{V}_{d}$ will be identified with $\mathcal{D}_{L'} \rightarrow \mathcal{V}_{e}$ because $\Phi(x_{8}) =  y_{7}$ and the planar isotopy identified the Khovanov maps. It remains to see that $\Phi$ also identifies the perturbed complexes $\mathcal{D}_{L}$ and $\mathcal{D}_{L'}$. These are depicted as edge weighted diagrams in Figure \ref{fig:comp}.\\
\ \\
\noindent In each of the diagrams in Figure \ref{fig:comp}, label the arcs $A_{1}$, $A_{2}$, $A_{3}$ from left to right. Let $D_{A_{i}}$ be the map on states which takes $v_{+} \rightarrow v_{-}$ for the circle $C$ containing the arc $A_{i}$. From the argument for Redemeister II invariance, the differential in $\mathcal{D}_{L}$ is $\partial + x_{8}(D_{A_{1}} + D_{A_{2}})$ while the differential for $\mathcal{D}_{L'}$ is $\partial' + y_{7}(D_{A_{2}} + D_{A_{3}})$. In addition, if we let $A$ be the region between $A_{1}$ and $A_{2}$ and $B$ be the region between $A_{2}$ and $A_{3}$ then there is a relationship between the areas in the resolved diagrams for $L$ and $L'$: $[A] = \Phi^{-1}([A']) + x_{8}$ and $[B']=\Phi([B]) + y_{7}$. 

\noindent To understand the relationship between the areas of discs in the resolutions comprising $\mathcal{D}_{L}$ and those for $\mathcal{D}_{L'}$, we start by dissecting along the arcs which form the boundary of the diagram \ref{fig:comp}. Using the additional regions gives a stably isomorphic complex as in lemma \ref{lem:decomp}. The area of regions outside the rectangle will be mapped correctly by $\Phi$ as deformation can be thought of as adding ``area'' for the arcs $A_{1}$ and $A_{2}$ for $L$ and $A_{2}$, $A_{3}$ for $L'$. We analyze what happens to each of the four regions in the rectangle. 
\begin{enumerate}
\item Region I: for a resolution of $L$ this region contributes $x_{1} + x_{8}$ after deforming with image $\Phi(x_{1} + x_{8}) = (y_{1} + y_{8}) + y_{8} = y_{1}$; for the corresponding resolution $L'$ region I abuts neither of the arcs $A_{2}$ or $A_{3}$, so it corresponds to multiplication by $y_{1}$ only. 
\item Region II: for a resolution of $L$, region II contributes $[A] = \Phi^{-1}([A'])+x_{8}$ before deforming, and $\Phi{-1}([A']) + x_{8} + x_{8} + x_{8}$ after deforming (one for each arc $A_{1}$ and $A_{2}$). $\Phi$ maps this contribution to $[A'] + y_{7}$. For the corresponding resolution of $L'$, before deformation we have multiplication by $[A']$, but afterwards, since region II abuts $A_{2}$, we have multiplication by $[A'] + y_{7}$ as well. 
\item Region III: for a resolution of $L$, before deformation, region III contributes $[B]$; after deformation it contributes $[B] + x_{8}$ as $A_{2}$ is in its boundary. $\Phi$ maps this contribution to $\Phi([B]) + y_{7} = [B']$. For $L'$ before deformation we have multiplication by $[B']$, and after deformation we obtain multiplication by $[B'] + y_{7} + y_{7} = [B']$, where one $y_{7}$ comes from $A_{2}$ and the other from $A_{3}$. 
\item Region IV: for $L$, both before and after deformation, region IV contributes $x_{4}$. Applying $\Phi$ gives $y_{4} + y_{7}$. For $L'$, before deformation, region IV contributes $y_{4}$, but after it contributes $y_{4} + y_{7}$, since the region abuts $A_{3}$. 
\end{enumerate}
\ \\
\noindent Consequently, $\Phi$ adjusts weights appropriately on of the four regions. A disc, $\disc{C}$, in a resolution can be thought of as a union of multiplicity 1 regions in the resolution diagram. In particular, it intersects the rectangular neighborhood, $N$, of the Reidemeister move in a subset of the four regions, each counted with multiplicity 1. We have already noted that the areas of regions outside the rectangle are mapped by $\Phi$ correctly into the area in the resolution of $L'$. If the disc only intersects the rectangle in one region, then one of the four calculations above applies. For the remaining cases suppose $\disc{C} \cap N$ is a union of the regions above counted with multiplicity one. Then the deformed area of $\disc{C}$ is the area of $\disc{C}$ plus the contribution of the arcs in $C$. This is the same as the sum of the deformed areas of the multiplicity one regions, since when two regions share an arc, the contribution of the arc will cancel in the sum, corresponding to the arc not being in $C$. On the other hand, if an arc is not shared by two regions, it will be in $C$ and thus the contribution of the arc should be counted, just as it is in the deformed area of the corresponding regions. In short, if $\disc{C}$ is the union of regions $R_{i}, \ldots, R_{k}$, with some outside of $N$, then its deformed area is the sum of the deformed areas of these regions.\\
\ \\
\noindent Thus, $\Phi$ maps the vertical differentials in $\mathcal{D}_{L}$ to the corresponding differential in $\mathcal{D}_{L'}$ since it maps the area of any disc correctly. After simplifying the top layers, $\Phi$ is a chain isomorphism of the whole complex, built out of an automorphism of the coefficient rings. Since the simplifications themselves are chain homotopy equivalences, and the automorphism is a chain isomorphism, the two complexes are (stably) chain homotopy equivalent (through $\F[x_{1}, \ldots, x_{n}, y_{7}]$).   

\section{The spanning tree deformation of $KT_{\ast}(L) \otimes \field{L}$}

\noindent In the remaining sections we will work over the {\em field} $\field{L}$, the field of fractions of $\ring{L}$. The complex $KT_{\ast}(L)$ is defined over $\field{L}$, and the arguments for invariance transfer directly. Thus the homotopy equivalence class of $KT_{\ast}(L)$ is invariant up to changes of marked point. Over $\field{L}$, however, we will be able to show that the invariance under change of marked point as well. In fact, the complex as a whole can be deformed using lemma \ref{lem:deform} into the one in the introduction. We start by deriving the form in the introduction, and then prove the irrelevance of changes in the marked point. \\ 
\ \\ 
\begin{theorem}
Over $\field{L}$, $\CT{L}$ is isomorphic to a deformation retract of $KT_{\ast}(L)$. In fact, the map $\mathcal{V}(L_{S}) \rightarrow |S|$ is a filtration on $KT_{\ast}(L)$. Then the associated spectral sequence has $E^{2}$-page isomorphic to $\CT{L}$, and collapses at the $E^{3}$-page.  
\end{theorem}    

\begin{center}
\begin{figure}
\treefig{0.7}{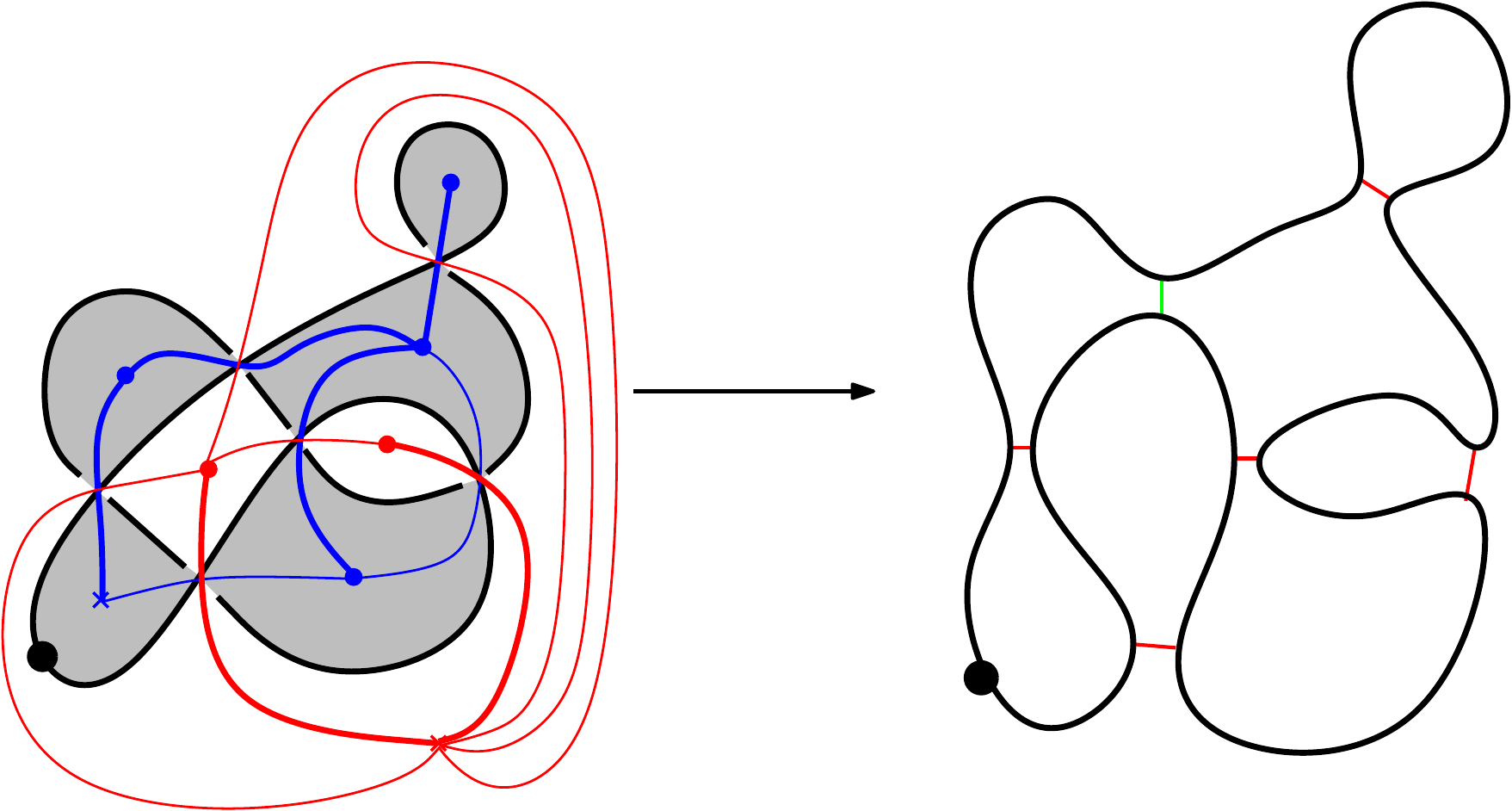} 
\caption{The dual spanning trees for the white and black Tait graphs on the left correspond to the single circle resolution diagram on the right.}
\label{fig:spanning}
\end{figure}
\end{center}

\noindent{\bf Proof:} For $S \subset \cross{L}$ with  $|S|=i$  of $i$, the complex, $\mathcal{V}(L_{S})$ is of the form $\F_{L} \otimes \mathcal{V}_{C_{1}} \otimes \cdots \otimes \mathcal{V}_{C_{k}}$ where $k \geq 1$ and $C_{1}$ is the marked circle. Since these complexes are over a field, the homology of the tensor product is the tensor product of the homologies. Since $[A_{C_{i}}]$ is invertible for all $i$ we have $\partial_{C_{i}}$ is an isomorphism from $\field{L}v_{+}$ to $\field{L}v_{-}$. Recall that $O(L)$ is the set of $S \subset \cross{L}$ such that the associated resolution $L_{S}$ consists of a single circle. We can decompose $KT_{\ast}(L)$ into a direct sum of three pieces 1) $K_{0} = \bigoplus_{S \in O(L)} \mathcal{V}(L_{S})$, 2) $K_{+} \cong \bigoplus_{S \not\in O(L)} \widetilde{v}_{0} \otimes v^{2}_{+} \otimes \mathcal{V}_{C_{2}} \otimes \cdots \mathcal{V}_{C_{k(S)}}$ and 3) $K_{-} \cong \bigoplus_{S \not\in O(L)} \widetilde{v}_{0} \otimes v^{2}_{+} \otimes \mathcal{V}_{C_{2}} \otimes \cdots \mathcal{V}_{C_{k(S)}}$. Then the differential induces an isomorphism $K_{+} \rightarrow K_{-}$ since the it contains a term $\widetilde{v}_{0} \otimes v^{2}_{+} \otimes W \stackrel{\cdot[\disc{C_{2}}]}{\longrightarrow} \widetilde{v}_{0} \otimes v^{2}_{-} \otimes W$ which is an isomorphism of vector spaces, and which is the only map preserving the value of $|S|$ with this image. Consequently, we can cancel $K_{+}$ and $K_{-}$ through this map, leaving $K_{0}$. Furthermore, since this map preserves the value of $|S|$, the resulting filtered chain complex has an isomorphic spectral sequence. In particular, all of the pages will be isomorphic.\\
\ \\
\noindent We now compute the perturbed boundary map. Since $K_{0}$ consists of single circle resolutions, the only generator for the chain group $\mathcal{V}(L_{S})[|S|]\{|S|\}$ occurs in Khovanov bigrading $(|S|, |S|)$ (for the unshifted complex). Note that any non-trivial contribution to the perturbed boundary map must increase $\delta = 2i - j$ only by $2$ since we canceled terms in the differential. Consequently, if $(\Delta i, \Delta j)$ is the change in the Khovanov bigrading, then $2\Delta\,i - \Delta\,j = 2$. If we start with a single circle resolution $S$, then for there to be a non-trivial map in the deformed complex to another single circle resolution $S'$ we will have $\Delta\,i = \Delta\,j = |S'| - |S|$. Consequently, $|S'| - |S| = 2$ and $|S'| = |S| + 2$. Thus the spectral sequence will have $E^{1}$-page with trivial differentials, potentially non-trivial $E^{2}$-page, but collapse at the $E^{3}$-page.\\
\ \\ 
\noindent To compute the differential at $E^{2}$ note that to get a single circle resolution, $S'$, with $|S'| - |S| = 2$ we must have $S'\backslash S = \{c_{1}, c_{2}\}$. Both resolutions $S \cup \{c_{1}\}$ and $S \cup \{c_{2}\}$ must have two circles in their diagrams. Now consider the image $\partial_{KH} \widetilde{v}_{S}$ in $\mathcal{V}(L_{S \cup \{c_{1}\}}) \cong \widetilde{v} \otimes \big(\field{F}v_{+} \oplus \field{F}v_{-}\big)$. In each resolution $\partial_{KH}(\widetilde{v}_{S}) = \widetilde{v} \otimes v_{-}$ which is canceled by $\frac{1}{\area{C_{2}}}\,\big(\widetilde{v} \otimes v_{+}\big)$.  Furthermore, under $\partial_{KH}$, $L_{S \cup \{c_{i}\}} = C_{\ast} \cup C_{2}$ is resolved to $L_{S'}$ by by merging $C_{\ast}$ and $C_{2}$ at the other crossing. Doing this for both $c_{1}$ and $c_{2}$ results in
$$
\partial^{\mathrm{pert}} \widetilde{v}_{S} = \left(\frac{1}{\area{C_{S\cup \{c_{1}\}, 2}}} + \frac{1}{\area{C_{S\cup\{c_{2}\},2}}}\right)\widetilde{v}_{S'}
$$
which is easily seen to agree with the differential for $\CT{L}$ described in the introduction.$\Diamond$\\
\ \\
\noindent We can interpret this complex in terms of a checkerboard coloring of the regions in the diagram $L$. The Tait graphs for $L$ are two planar graphs, one for each color. The black Tait graph has the black colored regions as vertices.  Each crossing $c$ of $L$ provides an edge joining the vertices containing  the two diagonally opposite black quadrants at $c$. We will usually draw the black Tait graph embedded in the union of the black regions and the projection of $L$ (see Figure \ref{fig:spanning}). The white Tait graph is defined similarly, using the white regions as vertices and all the crossings to provide edges. The colors will be used only to identify the graph, and will not otherwise be prescribed. Let $S \subset \cross{L}$ such that $L_{S}$ is a single circle. At each crossing $c$ the resolution bridges either the two black quadrants or the two white quadrants. If we let $T$ be the subset of $\cross{L}$ which bridge opposite black quadrants, we can consider the corresponding edges in the black Tait graph. The subgraph formed by these edges is a deformation retract of the union of the black regions in $L_{S}$, which is a disc. Consequently, $T$ determines a sub-tree of the black Tait graph, which is necessarily spanning since there is only one disc. Identically, the subset of $\cross{L}$ where the resolution $L_{S}$ bridges the white quadrants can be identified with a (dual) spanning tree for the white Tait graph. Furthermore, we will take the white and black regions on either side of the marked point, $p \in L$, as roots for the Tait graphs, and thus for all the spanning trees.\\
\ \\
\noindent On the other hand, given a partition of $\cross{L}$ into two sets $T$ and $T'$ with $T$ determining a spanning tree for the black Tait graph and $T'$ a dual spanning tree for the white Tait graph, we can resolve $L$ so that the edges in $T$ correspond to those crossings where the resolution bridges black quadrants, and the edges of $T'$ correspond to those crossings where the resolution bridges the white quadrants. The resulting digram $L_{S}$ will consist of a single circle (note that $S$ must then be determined by which crossings are resolved according to the rule in the introduction). Consequently, there is a one to one correspondence between spanning trees for the black Tait graph of a link projection and the generators of $\CT{L}$. For an illustration, see Figure \ref{fig:spanning}. We will often use these trees as generators spanning the chain groups. \\

\noindent Note also that when we change the resolution at a crossing of a single circle resolution, we cut off either a black disklike region from the black disk or a white disklike region from the white disk. To get back to a single circle resolution the second crossing change must rejoin these. Done in reverse order, the color of the cleaved off disk also changes. \\
\ \\
\noindent With this deformation equivalent complex in hand, we can resolve the difficulties surrounding invariance with respect to the choice of basepoint\footnote{The author does not know if this is possible in the totally twisted theory over $\ring{L}$.}
 
\begin{lemma}\label{lem:marked}
Let $(L,p)$ and $(L,p')$ be two marked projections of $\lk{L}$ which differ only in the marked point. Then there is a field isomorphism $I: \F_{L} \rightarrow \F_{L}$ such that $\CT[p]{L} \otimes_{I} \F_{L}$ is isomorphic to $\CT[p']{L}$. In short, the stable isomorphism class of $\CT{L}$ does not depend upon the choice of marked point on $L$. 
\end{lemma}

\noindent{\bf Proof of lemma \ref{lem:marked}:} Let $B \subset \regions{L}$ be the set of all black regions, and $W$ be the set of all white regions. Suppose the black region abutting the old marked point, $p$, is $s_{B}$ and the black region abutting the new marked point, $p'$, is $s_{B}'$ (these could be the same). Suppose similarly that $s_{W}$ is the white region abutting the old marked point and $s_{W}'$ is the white region abutting the new one. There are two cases to consider for the black regions and two generators $ T \in O_{i}(L)$ and $T'$ with $T' \in O_{i+2}(T,L)$. Let $T < b < T'$ where $b$ alters the resolution on the arc crossing the black regions.
\begin{enumerate}
\item Both $s_{B}$ and $s'_{B}$ are in the same region of $b(L)$: In this case $B^{p}_{T,T'}$ is the same as $B^{p'}_{T,T'}$, and neither contains $s_{B}$ or $s'_{B}$. 
\item Each of the two circles in $b(L)$ contains one of $s_{B}$ and $s_{B}'$. In this case $B^{p'}_{T,T'} $ $=B\backslash B^{p}_{T,T'}$. Furthermore, $[B^{p'}_{T,T'}]$ contains a single $x_{s_{B}}$ summand and $[B^{p}_{T,T'}]$ contains a single $x_{s_{B}'}$ summand. 
\end{enumerate}
A similar pattern holds for the white regions. We will define a field isomorphism of $\F_{L}$ by taking 
$$
\begin{array}{l}
x_{s_{B}'} \longrightarrow x_{s_{B}'} + [B] \\
x_{s_{W}'} \longrightarrow x_{s_{W}'} + [W] \\
x_{j} \longrightarrow x_{j} \\
\end{array}
$$
In the first case above, $[B^{p}_{T,T'}]$ is fixed by this automorphism, and equals $[B^{p'}_{T,T'}]$. In the second case, $[B^{p}_{T,T'}] = x_{s_{B}'} + [B^{p}_{T,T'}\backslash \{s_{B}'\}]$ is mapped to $x_{s_{B}'} + [B] + [B^{p}_{T,T'}\backslash \{s_{B}'\}]$ $= [B] + [B^{p}_{T,T'}]$ $= [B\backslash B^{p}_{T,T'}]$ $=[B^{p'}_{T,T'}]$. Thus the coefficient of $T'$ in $\partial^{p}_{i,L}$ coming from the black regions is mapped, under the automorphism, to the coefficient of $T'$ in $\partial^{p'}_{i,L}$ defined from the black regions. Mutatis mutandis, the result also holds for the white regions. $\Diamond$ \\
\ \\
\noindent Thus, at this point we have established that $\CT{L}$, up to stable homotopy equivalence, is a link invariant, and thus its homology is also a link invariant. 

\section{Verifying $\partial^{2} = 0$ for the spanning tree differential without reference to Khovanov homology}\label{sec:boundary}

\noindent One can prove that $\CT{L}$ is a chain complex directly from the explicit representation of its differential, and without the circuitous route through twisted Khovanov homology. This is done in the next proposition, which serves as a good introduction to the combinatorial complexities we avoided (somewhat) by using twisted Khovanov homology. It is easy to also prove invariance under the Reidemeister I moves directly from the definition of $\CT{L}$. It is also possible to directly prove invariance under the RII move, although this is a substantially more involved combinatorial proof. However, the author has not been able to prove RIII invariance without using twisted Khovanov homology. For now we content ourselves with proving, directly from the definition, that the boundary map for $\CT{L}$ really is a differential. 

\begin{prop}\label{lem:square}
For the map $\partial_L$ we have $\partial_L^{2} \equiv 0$. 
\end{prop}

\noindent{\bf Proof of lemma \ref{lem:square}:} Let $T \in O_{i}(L)$. Take $p$, the marked point, and move it to infinity. Then we may think of $T$ as
the $y$-axis in the plane, and the arcs from the resolutions as semi-circles whose ends lie on this axis, and which are wholly contained either in $x \geq 0$ or $x \leq 0$. We may choose the black region inside $T$ to correspond to the set $x \geq 0$. The endpoints of each arc $c$ cut the $y$-axis into three segments, $U_{\pm}(c)$, the segment unbounded towards $\pm\infty$, and $\gamma(c)$, the bounded segment. 

\begin{defn}
Two disjoint arcs $c_{1}$ and $c_{2}$ in $x \geq 0$ will be called {\it parallel} if $\gamma(c_{2})\subset\gamma(c_{1})$, or vice-versa. If $\gamma(c_{1}) \cap \gamma(c_{2}) = \emptyset$ then we will call the arcs {\it peers}. An arc $c$ in $x \geq 0$ and an arc $a$ in $x \leq 0$ will be said to {\it interleave} if $\gamma(a) \cap \gamma(c) \neq \emptyset$ but $\gamma(a) \not\subset \gamma(c)$ and $\gamma(c) \not\subset \gamma(a)$. 
\end{defn}

\noindent To show that $\partial^{2} = 0$ we compute
$$
\partial^{2} T = \sum_{T' \in O_{i+4}(T,L)} \left(\sum_{T' > T_{r} > T} \pairing{T'}{T_{r}}\pairing{T_{r}}{T}\right) T'
$$
We will show, for each $T' \in O_{i+4}(L)$, that $\sum_{T' > T_{r} > T} \pairing{T'}{T_{r}}\pairing{T_{r}}{T} = 0$ in $\F_{L}$. Each term in this sum corresponds to four arcs, each coded with a $0$, on the diagram for $T$ considered in $\R^{2}$, two arcs in $x \leq 0$ and two in $x \geq 0$. Let these arcs be $\{a, a'\}$ and $\{c, c'\}$ respectively. Then each $T_{r}$ in the summation corresponds to two interleaved arcs, $r=\{a, c\}$. We now forget the remainder of the arcs and concentrate only on these configurations.\\ \ \\

\noindent We will analyze configurations in the following cases: there is a labeling of the arcs in the white region as $a$ and $a'$ and the arcs
in the black region as $c$ and $c'$ such that
\begin{enumerate}
\item[]I. $a$ does not interleave with either $c$ or $c'$.
\item[]II. Each of $\{a,c\}$,$\{a',c'\}$, $\{a',c\}$ and $\{a,c'\}$ interleave.
\item[]III. $\{a,c\}$ and $\{a',c'\}$ interleave, but $\{a,c'\}$ and $\{a',c\}$ do not. 
\item[]IV. $\{a,c\}$ and $\{a',c'\}$ and $\{a',c\}$ interleave, but $\{a, c'\}$ do not.
\end{enumerate}
The remaining cases can be obtained by either by switching the roles of $a'$ and $a$, $c'$ and $c$, in the interleaving of the last case, or by arguing by symmetry between the white and black regions. 
We now analyze each of the cases.\\ \ \\
 
\noindent{\bf Case I:} Resolving $a$ results in a new circle component which cannot be rejoined to the other components by resolving along either $c$ or $c'$. Consequently, the result of resolving all four arcs is not a single circle, and there is no contribution to $\partial^{2}_L$.\\ \ \\

\noindent{\bf Case II:} $a$ and $a'$ each interleave with both $c$ and $c'$. If either $a$ and $a'$ are parallel or $c$ and $c'$ are parallel, then there resolving along all four arcs does not result in a single circle resolution, $T'$, and thus this case does not contribute to $\partial^{2}_{L} T$. To see this, suppose $c$ and $c'$ are parallel and $\gamma(c') \subset \gamma(c)$. Their mutual resolution results in a new circle component between the two arcs, which intersects the $y$-axis in segments $s =\gamma(c) \backslash \gamma(c')$. The endpoints of $a$ cannot be on the segments in $s$: if there was an endpoint in one of the segments, then for $a$ and $c'$ to interleave, the other endpoint would need to be in $\gamma(c')$, but then $a$ and $c$ would not interleave. The same argument applies to $a'$. Consequently, resolving along $a$ and $a'$ does not affect the new circle component and $T'$ is not a single circle. By symmetry, the if $a$ and $a'$ are parallel then this case does not contribute to $\partial^{2}$. So assume that $c$ and $c'$ are peers. Since $a$ and $a'$ interleave both, $a$ must have an endpoint in $\gamma(c)$ and in $\gamma(c')$ as these segments are disjoint. So must $a'$, and since $a$ and $a'$ are disjoint they will have to be parallel. Thus, one or both pairs $\{a,a'\}$ or $\{c,c'\}$ are parallel and this case does not contribute to $\partial_{L}^{2}$.\\ \ \\

\noindent{\bf Case III:} In this case, the arc pairs $\{a,c\}$ and $\{a',c'\}$ are independent. Let $T_{r}$ be the result of resolving along $\{a,c\}$ and $T_{r}'$ be the result of resolving along $\{a',c'\}$. If $c$ and $c'$ are pairs, then $B_{T,T_{r}}$ and $B_{T,T'_{r}}$ are disjoint. If we have resolved $c$ and then resolve $c'$, the region cut off is the same as if we resolve $T$ along $c'$, i.e. $B_{T_{r},T'} = B_{T,T_{r}'}$. Likewise $B_{T_{r}',T'} = B_{T,T_{r}}$. By symmetry, the same argument holds for $a$ and $a'$ when they are peers. Now suppose $c$ and $c'$ are parallel with $\gamma(c') \subset \gamma(c)$. Then 
$B_{T,T_{r}'} \subset B_{T,T_{r}}$. If we resolve first along $\{a,c\}$, then we rejoin $B_{T,T_{r}}$ to the unbounded black region, without affecting the
region cut out by the arc $c'$ since the endpoints of $a$ do not intersect $\gamma(c')$. Consequently, $[B_{T_{r},T'}] = [B_{T,T_{r}'}]$ since the formal variables are unchanged. Furthermore, if we first resolve $\{a',c'\}$, since $a'$ does not interleave $c$, one endpoint of $a'$ is in $\gamma(c')$ and the other is in $\gamma(c) \backslash \gamma(c')$. Thus the region $B_{T,T_{r}'}$ is rejoined to $B_{T,T_{r}}\backslash B_{T,T_{r}'}$ by $a'$, so that $[B_{T_{r}',T'}] = [B_{T, T_{r}}]$. A similar argument applies to the white regions.  Considerations of this type for each of the possible peer/parallel configurations yields that $[B_{T,T_{r_{2}}}] = [B_{T_{r_{1}}, T'}]$, $[B_{T_{r_{2}}, T'}] = [B_{T,T_{r_{1}}}]$ and likewise, $[W_{T,T_{r_{2}}}] = [W_{T_{r_{1}}, T'}]$ and $[W_{T_{r_{2}}, T'}] = [W_{T,T_{r_{1}}}]$, and it is straightforward to see that the contribution to $\pairing{\partial_{L}^{2}T}{T'}$ coming from these arc pairs cancels in the summation.\\ \ \\

\noindent{\bf Case IV:}   Let $T_{r}$ result when resolving $\{a,c\}$, $T_{r}'$ result when resolving $\{a',c'\}$ and $T_{s}$ result when resolving $\{a',c\}$. For each of these circles, resolving the remaining two arcs results in the same circle $T'$. 
\begin{enumerate}
\item For $T \rightarrow T_{r} \rightarrow T'$: First we cut off region $B_{c}$ and $W_{a}$ when changing from $T$ to $T_{r}$. Since $a$ interleaves with $c$, but not with $c'$ $B_{c}$ is rejoined to the same component of $\{x \geq 0\} \backslash c'$ as it was cut from. Consequently, resolving $c'$ on $T_{r}$ cuts off a region with formal representative $[B_{c'}]$. Meanwhile, resolving $c$ rejoins the region $W_{a}$ to 
\begin{enumerate}
\item $W_{a'}$ in the case that they are peers. When resolving $a'$ we will cut off a region with formal area $[W_{a}] + [W_{a'}]$,
\item the unbounded white region if $a$ and $a'$ are parallel with $\gamma(a) \subset \gamma(a')$. Thus resolving $a'$ will cut off the region between $a'$ and $a$ which has formal area $[W_{a'}] - [W_{a}]$. Since we are working in characteristics two, this equals $[W_{a}] + [W_{a'}]$.
\item $W_{a'}$ in the case that $a$ and $a'$ are parallel with $\gamma(a') \subset \gamma(a)$. Resolving $a'$ will cut off the region between $a$ and $a'$ with area $[W_{a}] - [W_{a'}]$ which is the same as $[W_{a}] + [W_{a'}]$. 
\end{enumerate}
For all three cases the contributions to $\partial_L^{2}$ are the same: 
$$
\left(\frac{1}{[B_{c}]} + \frac{1}{[W_{a}]}\right)\left(\frac{1}{[B_{c'}]} + \frac{1}{[W_{a}] + [W_{a'}]}\right)
$$ 
\item For $T \rightarrow T_{r}' \rightarrow T'$: First we cut off $B_{c'}$ and $W_{a'}$. By a similar argument as above, the exact configuration of parallel and peer arcs for $a$ and $a'$ and $c$ and $c'$ will not matter. So we can do the calculation when each pair of arcs is peer. This results in the following contribution:
$$
\left(\frac{1}{[B_{c'}]} + \frac{1}{[W_{a'}]}\right)\left(\frac{1}{[B_{c}] + [B_{c'}]} + \frac{1}{[W_{a}]}\right)
$$
\item For $T \rightarrow T_{s} \rightarrow T'$: First we cut off $B_{c}$ and $W_{a'}$. Again we can do the calculation only for the peer case:
$$
\left(\frac{1}{[B_{c}]} + \frac{1}{[W_{a'}]}\right)\left(\frac{1}{[B_{c}]+ [B_{c'}]} + \frac{1}{[W_{a}]+ [W_{a}']}\right)
$$
\end{enumerate}
We will now combine fractions, multiply, and simplify. To aid us, note that
$$
\frac{1}{XY} + \frac{1}{X(X+Y)} + \frac{1}{Y(X+Y)} = \frac{(X+Y) + Y + X}{XY(X+Y)} = 0
$$
in $\F_{L}$ for any non-zero elements $X$ and $Y$ we choose. Observe that in the products above, if we take only those terms using for the black regions we obtain such a sum with $X = B_{c}$ and $Y = B_{c'}$. Likewise if we take only those terms involving the white regions, we get such a sum with $X = W_{a}$ and $Y=W_{a'}$. So we only need to consider the sum of the cross-terms:
\begin{eqnarray*}
\frac{1}{[W_{a}][B_{c'}]} + \frac{1}{[B_{c}]([W_{a}] + [W_{a'}])} + \frac{1}{[W_{a'}]([B_{c}] + [B_{c'}])}  
\\ + \frac{1}{[B_{c'}][W_{a}]} + \frac{1}{[W_{a'}]([B_{c}]+ [B_{c'}])} + \frac{1}{[B_{c}]([W_{a}]+ [W_{a}'])}
\end{eqnarray*}
However, these come in canceling pairs, so the sum of the three products is zero.\\ \ \\

\noindent Consequently, for each of the four cases we have $\sum_{T' > T_{r} > T} \pairing{T'}{T_{r}}\pairing{T_{r}}{T} = 0$, and thus $\partial_{L}^{2} = 0$.
$\Diamond$

\section{Properties}\label{sec:prop}

\subsection{A long exact sequence}\label{sec:LES}

\begin{center}
\begin{figure}
\treefig{0.5}{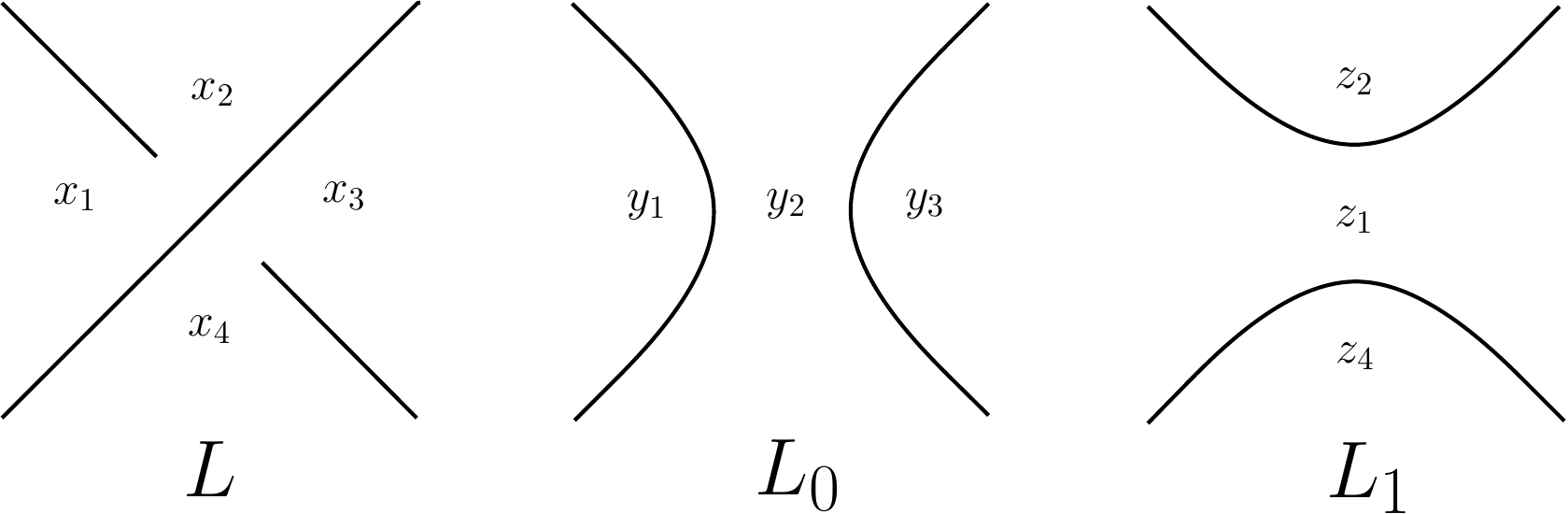} 
\caption{}
\label{fig:long}
\end{figure}
\end{center}

\noindent Let $c \in \cross{L}$ as depicted in Figure \ref{fig:long}. We have labeled the regions abutting this crossing $x_{1}, x_{2}, x_{3}$ and $x_{4}$. A priori some of these could be equal: $x_{1}$ could equal $x_{3}$ or $x_{2}$ could equal $x_{4}$. If so, one of the resolutions of $L$ at this crossing will result in a disconnected diagram. For the moment, we require that all four regions be distinct, i.e. that both the $0$ resolution and the $1$ resolution at $c$ result in connected link diagrams. 

\begin{theorem}
Let $L_{0}$ and $L_{1}$ be the diagrams found by resolving $L$ using the $c \not\in S$ and $c \in S$ rules, respectively. Then 
$$
\CT{L} \cong \mathrm{MC}\big(\CT{L_{0}}\otimes_{\F_{L_{0}}} \F_{L} \stackrel{\tau_{c}}{\longrightarrow} \CT{L_{1}} \otimes_{\F_{L_{1}}} \F_{L}\big)
$$
where
\begin{itemize}
\item  $\F_{L_{0}}$ acts on $\F_{L}$ by $y_{1}, y_{3} \rightarrow x_{1},x_{3}$ and $y_{2} \rightarrow x_{2} + x_{4}$ and $\F_{L_{1}}$ acts by  $z_{1} \rightarrow x_{1} + x_{3}$ and
$z_{2},z_{4} \rightarrow x_{2},x_{4}$.
\item If $T \in O_{i}(L_{0})$ then $$\tau_{c}(T) = \sum_{T' \in O_{i+1}(L_{1})} \pairing{T}{T'}_{\F_{L}} T'$$ and $T'$ is a single circle resolution of $L_{1}$ differing from $T$ at $c$ and some other crossing.
\end{itemize}
\end{theorem}

\noindent{\bf Proof:} $\CT{L}$ can be decomposed along those resolutions $T$ with $T(c) = 0$ and $T(c) = 1$, $\CT{L} \cong C_{0} \oplus C_{1}$. Those $T$ with $T(c) = 0$ restrict to resolutions on $L_{0}$ by actually resolving the crossing $c$. Likewise, those with $T(c) = 1$ provide single circle resolutions on $L_{1}$. Furthermore, all single circle resolutions for $L_{0}$ and $L_{1}$ arise in this manner. If we consider $\partial_{L}$ we can decompose into three parts $\partial_{0} \oplus \partial_{01} \oplus \partial_{1}$ where $\partial_{i}$ counts those pairs $T \Rightarrow T'$ which have $T(c) = T(c) = i$ and $\partial_{01}$ counts pairs where $T(c) = 0$ and $T'(c) = 1$. Then $(C_{0}, \partial_{0}) \cong \CT{L_{0}} \otimes_{\F_{L_{0}}} \F_{L}$, where the tensor product arises because two of the formal variables for $L$, $x_{2}, x_{4}$ occur in the same region in $L_{0}$. Using the results of the previous section, we see that this only changes the coefficient field. Likewise, $(C_{1}, \partial_{1}) \cong \CT{L_{1}} \otimes_{\F_{L_{1}}} \F_{L}$. The map $\tau_{c}$ comes from $\partial_{01}$, an its form is readily descried from that of the differential, $\partial_L$. Lastly we verify the shift in gradings: note that $C_{0} \cong \CT{L_{0}}$ with no shift, since $T(c) = 0$ has no effect on $\widetilde{\delta}$. However, $C_{1} \cong \CT{L_{1}}[1]$ since $T(c) = 1$ implies that the grading on $\CT{L_{1}}$, where no resolution occurs, will be shifted up when considered in $\CT{L}$. Namely, $C_{1;i} \cong \CT[i-1]{L_{1}}$. Since $\partial_L$ increases $\widetilde{\delta}$ by $2$, this is the correct shift for a mapping cone. $\Diamond$. \\ \ \\

\begin{prop}
Given a crossing $c$ there is a long exact sequence \\

\begin{equation}
\cdots \rightarrow \unHT[i-1]{L_{1}}\otimes \F_{L} \rightarrow \unHT[i]{L} \rightarrow \unHT[i]{L_{0}}\otimes \F_{L} \stackrel{\tau_{c,\ast}}{\longrightarrow} \unHT[i+1]{L_{1}}\otimes \F_{L} \rightarrow \cdots
\label{eqn:unshiftedseq}
\end{equation}
\ \\
\noindent When $L$ is oriented and $c$ is a positive crossing, then if $e=n_{+}(L) - n_{+}(L_{1})$ (for any orientation on $L_{1}$), then\\
\begin{equation}
\cdots \rightarrow \HTsub[i+e-1]{L}{1}\otimes \F_{L} \rightarrow \HT[i]{L} \rightarrow \HTsub[i+1]{L}{0}\otimes \F_{L} \stackrel{\tau_{c,\ast}}{\longrightarrow} \HTsub[i+e+1]{L}{1}\otimes \F_{L} \rightarrow \cdots
\label{eqn:posshiftedseq}
\end{equation}
\ \\
\noindent On the other hand, if $c$ is negative, let $f = n_{+}(L) - n_{+}(L_{0})$. Then \\
\begin{equation}
\cdots \rightarrow \HTsub[i-1]{L}{1}\otimes \F_{L} \rightarrow \HT[i]{L} \rightarrow \HTsub[i+f]{L}{0}\otimes \F_{L} \stackrel{\tau_{c,\ast}}{\longrightarrow} \HTsub[i+1]{L}{1}\otimes \F_{L} \rightarrow \cdots
\label{eqn:negshiftedseq}
\end{equation}
\end{prop}

\noindent{\bf Proof:} The exact sequence \ref{eqn:unshiftedseq} is an immediate consequence of the description of $\CT{L}$ as a mapping cone, using standard homological algebra. To verify \ref{eqn:posshiftedseq} assume that $c$ is positive. We will adjust subscripts to account for the shifting of $\unHT{L}$ by $[-n_{+}(L_{+})]$:
$$
\cdots \rightarrow \unHT[i+n_{+}-1]{L_{1}}\otimes \F_{L} \rightarrow \unHT[i+n_{+}]{L} \rightarrow \unHT[i+n_{+}]{L_{0}}\otimes \F_{L} \stackrel{\tau_{c,\ast}}{\longrightarrow} \unHT[i+n_{+}+1]{L_{1}}\otimes \F_{L} \rightarrow \cdots
$$
\ \\
$$
\cdots \rightarrow \unHT[i+n_{+}(L_{1}) + e-1]{L_{1}}\otimes \F_{L} \rightarrow \HT[i]{L} \rightarrow \unHT[i+n_{+}(L_{0}) + 1]{L_{0}}\otimes \F_{L} \stackrel{\tau_{c,\ast}}{\longrightarrow} \unHT[i+ n_{+}(L_{1}) + e + 1]{L_{1}}\otimes \F_{L} \rightarrow \cdots
$$
\ \\
$$
\cdots \rightarrow \HTsub[i + e-1]{L}{1}\otimes \F_{L} \rightarrow \HT[i]{L} \rightarrow \HTsub[i + 1]{L}{0}\otimes \F_{L} \stackrel{\tau_{c,\ast}}{\longrightarrow} \HTsub[i+e+1]{L}{1}\otimes \F_{L} \rightarrow \cdots
$$
\ \\
\noindent When $c$ is negative, we proceed as before\\
$$
\cdots \rightarrow \unHT[i+n_{+}-1]{L_{1}}\otimes \F_{L} \rightarrow \unHT[i+n_{+}]{L} \rightarrow \unHT[i+n_{+}]{L_{0}}\otimes \F_{L} \stackrel{\tau_{c,\ast}}{\longrightarrow} \unHT[i+n_{+}+1]{L_{1}}\otimes \F_{L} \rightarrow \cdots
$$
\ \\
\noindent But now $n_{+}(L) = n_{+}(L_{1})$ since we resolve a negative crossing. Furthermore, if we orient $L_{0}$ we may compute $f = n_{+}(L) - n_{+}(L_{0})$ and\\
$$
\cdots \rightarrow \unHT[i+n_{+}(L_{1}) -1]{L_{1}}\otimes \F_{L} \rightarrow \HT[i]{L} \rightarrow \unHT[i+n_{+}(L_{0}) + f]{L_{0}}\otimes \F_{L} \stackrel{\tau_{c,\ast}}{\longrightarrow} \unHT[i+ n_{+}(L_{1}) + 1]{L_{1}} \otimes \F_{L}\rightarrow \cdots
$$
\ \\
$$
\cdots \rightarrow \HTsub[i - 1]{L}{1}\otimes \F_{L} \rightarrow \HT[i]{L} \rightarrow \HTsub[i + f]{L}{0}\otimes \F_{L} \stackrel{\tau_{c,\ast}}{\longrightarrow} \HTsub[i+1]{L}{1}\otimes \F_{L} \rightarrow \cdots
$$ $\Diamond$

\section{The G\"{o}ritz matrix of $L$ (following \cite{MaOz})}\label{sec:goritz}

\noindent This section recalls some results concerning the signature and determinant of a link, which will be useful in the following sections. \\
\ \\
\noindent Let $L$ be an oriented link diagram and color the elements of $\regions{L}$ with the colors black and white, in checkerboard fashion. Let $\mathfrak{W} = \big\{W_{0}, \ldots, W_{n}\big\} \subset \regions{L}$ be the $n+1$ faces which are colored white. To each crossing $c \in \cross{L}$ assign we assign a value $\mu(c)$ according to:

\begin{center}
\treefig{0.35}{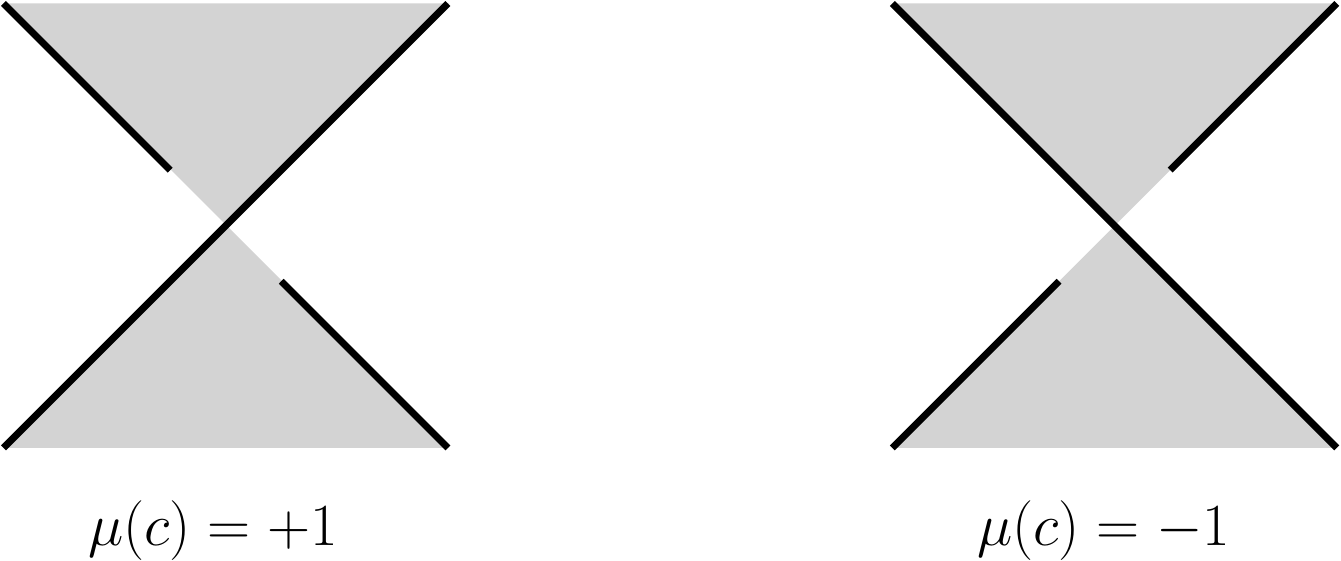}
\end{center}

\noindent Crossings with $\mu(c) = -1$ will be called compatible crossings, or $co$-crossings, since the black quadrants will be merged in any resolution where $c \in S$. The edges with $\mu = +1$ will be called incompatible, or $ic$-crossings. \\
\ \\
\noindent In addition, using the orientation we can assign a chirality: $n(c) = +1$ for positive crossings and $n(c) = -1$ for negative crossings. Finally, 
crossings will be called Type I if $\mu(c)n(c) = -1$ and Type II if $\mu(c)n(c) = +1$. \\
\ \\
\noindent Let
$$
\mu(L) = \sum_{\mathrm{c\ of\ Type\ II}} \mu(c)
$$
$$
g(W_{i}, W_{j}) =  -\sum_{c \in \overline{W_{i}} \cap \overline{W_{j}}} \mu(c)
$$
and
$$
g(W_{i}) = - \sum_{j \neq i} g(W_{i},W_{j})
$$
The G\"{o}ritz matrix of $L$ is the matrix $G(L)$ with $G_{ij}(L) = g(W_{i}, W_{j})$ for $ 0 \geq i \neq j \leq n$ and $G_{ii}(L) = g(W_{i})$ for $1 \leq i \leq n$. The matrix $G$ can be used to compute both the signature, $\sigma(\lk{L})$ and determinant $\mathrm{det}(\lk{L})$ (where the right handed trefoil is taken to have $\sigma = -2$) using
\begin{enumerate}
\item{\bf Gordon-Litherland Formula} $\sigma(\lk{L}) = \mathrm{sign}(G(L)) - \mu(L)$\\
\ \\
\item $\mathrm{det}(\lk{L}) = |\mathrm{det}(G(L))|$.
\end{enumerate}

\noindent If $L$ is an connected, reduced alternating diagram, then there is a simpler formulation. Checkerboard color the diagram so that every crossing is incompatible with the black regions. Then $\sigma(\lk{L}) = n - n_{+}$ where $n$ is as above. \\
\ \\ 
\section{Euler characteristic}\label{sec:euler}

\begin{theorem}
For $\lk{L}$, an oriented link in $S^{3}$, let
$$
P(\delta) = \sum_{i \in \Z} rk_{\F_{L}}\left(\HT[i]{L}\right)\delta^{i}
$$
then $det(\mathcal{L}) = |P(i)|$ where $i =\sqrt{-1}$.
\end{theorem}

\noindent{\bf Proof:} Let $L$ be a diagram for $\lk{L}$, and checkerboard color the faces of $L$, and let
$$
R(\delta) = \sum_{i \in \Z}rk_{\F_{L}}CT_{i+n_{+}(L)}(L) \delta^{i}
$$
be the Poincare\'e polynomial for the (shifted) chain groups. Then $|R(i)| = |P(i)|$, encoding that the Euler
characteristic can be determined from the chain groups. \\
\ \\
\noindent First, we gather some statistics for $W_{L}$, the white Tait graph for $L$. We let
\begin{itemize}
\item $V_{L}$ and $E_{L}$ are the number of vertices and edges in $W_{L}$, respectively;
\item $n_{0}$ the number of $co$-crossings, 
\item $\widetilde{e}(L) = n_{0} - V_{L} + 1$;
\item $\nu(L) = \widetilde{e}(L) - n_{+}(L)$
\end{itemize}
We recall the manner for relating $\det(L)$ to spanning trees through the Kirchoff Matrix-Tree theorem, \cite{Kauf}. We will start by working in the ring $\Z[\delta^{\pm 1}]$. Label each edge in $W_{L}$ with a $\kappa(c) = +1$ if $\mu(c) = +1$ for the corresponding crossing in $L$. Label the edge with $\kappa(c) = \delta$ if $\mu(c) = -1$ for the corresponding crossing. To each spanning tree, $T$, of $W_{L}$ let $w_{T}$ be the product of the labels attached to those edges in $T$. Then the Matrix-Tree theorem asserts
$$
\sum_{T \in trees(W_L)} w_{T} = \mathrm{det}\,[A]_{11}(S)
$$
where $[A]_{11}$ is the $(1,1)$-minor of the matrix $V_{L} \times V_{L}$ matrix $A$ formed using the elements
\begin{itemize}
\item $A_{ij} = -\sum_{c \in W_{i} \cap W_{j}} \kappa(c)$, i.e. the sum of labels of edges between vertex $i$ and vertex $j$
\item $A_{ii} = -\sum_{j \neq i} A_{ij}$;
\end{itemize}
Call the resulting polynomial $Q(\delta)$.\\
\ \\
\noindent If we specify $\delta \rightarrow 1$ then the determinant above computes the number of maximal spanning trees, \cite{Kauf}, for $W_{L}$. However, if we set $\delta = -1$, then the matrix $A$ is the G\"{o}ritz-matrix which has determinant, up to sign, equal to $det(\lk{L})$. Consequently, $|Q(-1)|=det(L)$. We now relate the polynomial $Q(\delta)$ to the polynomial $R(\delta)$. \\
\ \\
\noindent If $T$ is a tree in $W_L$, let $k$ be the number of edges in $T$ with label $ic$, and thus must be resolved using a $1$-resolution in the corresponding single circle resolution. In $Q(\delta)$ this contributes $\delta^{k}$. To obtain $\degrad{T}$ we must also count the number of $1$-resolutions on edges not in $T$. There are $V_{L} - 1$ edges in $T$, and thus $(V_{L} - 1) - k$ edges in $T$ adorned with $co$. Outside of $T$ there are $n_{0} - (V_{L} - 1) + k$ edges labeled with $co$. These, when outside a tree in $W_{L}$, receive a $1$ resolution in the corresponding single circle resolution. Thus, $T$ contributes
$\delta^{k} \cdot \delta^{n_{0} - V_{L} + 1 + k} = \delta^{2k} \cdot \delta^{\widetilde{e}(L)}$ in the polynomial for the unshifted complex, $\widetilde{R}(\delta) = \sum \mathrm{rk} CT_{i}{L} \delta^{i}$. Thus, $\widetilde{R}(\delta) = \delta^{\widetilde{e}(L)}Q(\delta^{2})$. Shifting alters the powers of $\delta$ by multiplying by $\delta^{-n_{+}(L)}$, and thus
$$
R(\delta) = \delta^{\nu(L)}Q(\delta^{2})
$$
We now plug in $\delta = i$ and take the complex modulus. This produces $|R(i)| = |(i)^\mu(L)Q(-1)| = det(\lk{L})$. The conclusion then follows from $|R(i)| = |P(i)|$ $\Diamond$\\
\ \\
\noindent Because $\partial_{L}$ shifts grading by $2$, plugging in $i$ above results in an Euler characteristic calculation. However, we cannot avoid the square root of $-1$, since after shifting by $-n_{+}$, we cannot ensure that all the exponents in $P(S)$ have even parity, a priori.

\section{Results for mirrors, connect sums, and quasi-alternating links}\label{sec:quasi}

\noindent We record two results which facilitate the calculation of the homology for knots and links built out of mirrors and connect sums. We then do a calculation which identifies the homology exactly for the class of quasi-alternating links. The results precisely mirror those for other knot homology theories. 

\begin{theorem}
Let $\lk{L}$ be an oriented link, then $\HT[i]{L} \cong \HT[-i]{\mir{L}}$.
\end{theorem}

\noindent{\bf Proof:} Let $L$ be a diagram for $\lk{L}$ and $\mir{L}$ be the mirror digram. We can use the same set of regions in forming $\F_{L}$ and $\F_{\mir{L}}$, and can thus identify the coefficient fields. Furthermore, $n_{+}(\mir{L}) = n_{-}(L)$ and $n_{-}(\mir{L}) = n_{+}(L)$. Every resolution $T \in O_{i}(L)$ corresponds to a resolution $\mir{T} \in O_{|\cross{L}| - i}(\mir{L})$ where $c \in T$ for the resolution of $L$ corresponds to $c \not\in \mir{T}$ for the resolution of $\mir{L}$.  A tree $T' \in O_{i+2}(T,L)$ corresponds to a tree $\mir{T}' \in O_{|\cross{L}| - i - 2}(\mir{L})$ and $\mir{T} \in O_{|\cross{L}| - i}(\mir{T'}, \mir{L})$. Thus $\mir{T}$ can appear in $\partial_{\mir{L}}\mir{T}'$. $T' = T \cup \{c_{1}, c_{2}\}$ with $c_{1}, c_{2} \not\in T$ corresponds to  $\mir{T'} = \mir{T} \backslash \{c_{1}, c_{2}\}$, or $\mir{T} = \mir{T'} \cup \{c_{1},c_{2}\}$ with $c_{1},c_{2} \not\in \mir{T'}$. When these alterations are performed the same two regions will be cut off, $B_{\mir{T}',\mir{T}} = B_{T,T'}$ and $W_{\mir{T}', \mir{T}} = W_{T,T'}$, and thus
$$
\pairing{\partial_{L} T}{T'} = \pairing{\mir{T}}{\partial_{\mir{L}}T'}
$$
Consequently, the differential $\partial_{\mir{L}}$ corresponds to the cohomology differential $\partial_{L}^{\ast}$. Thus the unshifted cohomology for $L$ in degree $i$ is isomorphic to $\unHT[|\cross{L}| - i]{\mir{L}}$. Since we are working with coefficients in a field, we have
$$
\begin{array}{c}
\HT[i]{L} \cong \unHT[i+n_{+}(L)]{L} \cong \unHT[|\cross{L}| - (i + n_{+}(L))]{\mir{L}} \cong \unHT[n_{-}(L) - i]{\mir{L}} \\
 \ \\
\cong \unHT[-i + n_{+}(\mir{L})]{\mir{L}} \cong \HT[-i]{\mir{L}}  
\end{array}
$$ 
$\Diamond$ \\ \ \\

\begin{theorem}
Let $\lk{L}_{1}$, $\lk{L}_{2}$ be two non-split oriented links, and let $\lk{L} = \lk{L}_{1} \# \lk{L}_{2}$, in some manner. Then
$$
\HT[k]{L} \cong \oplus_{i+j=k} \HTsub[i]{L}{1} \otimes \HTsub[j]{L}{2}
$$ 
where $\cong$ denotes stable equivalence.
\end{theorem}

\noindent{\bf Proof:} Let $L$ be a standard connect sum diagram in the plane for $\lk{L}$, where the portion of $L$ with $x < 0$ is a diagram for $\lk{L}_{1}$ with an arc removed and the portion of $L$ with $x> 0$ is a diagram for $\lk{L}_{2}$. We will choose our marked point for $L$ to lie on one of the two arcs intersecting $x=0$. For $L_{1}$ and $L_{2}$, we choose the marked point to lie on the removed arc. Let the regions for $L_{1}$ correspond to $x_{1}, \ldots, x_{k}$ and the regions for $L_{2}$ correspond to $y_{1}, \ldots, y_{l}$. $x_{1}$ and $y_{1}$ should correspond to the bounded regions abutting the marked points. We will ignore the unbounded region in all diagrams, as it will always abut the marked point, and thus not play a role in the calculations. In the diagram for $L$, there are $k + l -1$ regions which we will label $z_{1}, z_{2}, \ldots, z_{k}, \ldots, z_{l+k-1}$. $z_{1}$ corresponds to $x_{1} + y_{1}$, while $z_{i} \sim x_{i}$ for $2 \leq i \leq k$ and $z_{i} \sim y_{i-k+1}$ for $k + 1 \leq i \leq l + k - 1$ (where $\sim$ means corresponds to the ``same'' region). These identifications will be implicitly used in the stable equivalence.\\
\ \\
\noindent If we color black the bounded region abutting the basepoint in each diagram, then $B_{L}$ is $B_{L_{1}}$ and $B_{L_{2}}$ fused at their basepoint. Every spanning tree in $B_{L}$ is thus the fusion of a spanning tree for $B_{L_{1}}$ and $B_{L_{2}}$. As $\widetilde{\delta}$ and $n_{+}$ will both add under connect sums, if $T \in O_{i}(L_{1})$ and $T' \in O_{j}(L_{2})$ then $T \# T' \in O_{i+j}(L)$. To compute $\partial(T \# T')$ we need to consider trees in $B_{L}$ where we have removed one edge of $T\#T'$, and reconnected the resulting pieces with an edge of $B_{L} \backslash (T \# T')$. If the removal occurs in $B_{L_{1}}$, we have a disconnected component in $\{x < 0\}$. There are no edges which cross $x=0$, so the replacement must also occur with an edge from $B_{L_{1}} \backslash T$. The same argument applies if the edge removed occurs in $\{x>0\}$. These are precisely the trees which occur in either $\partial_{L_{1}}(T)$ or $\partial_{L_{2}}(T')$. Furthermore, due to the placement of the basepoints, if we measure all coefficients using $z_{i}$'s, then all the coefficients will also be the same as the connect sum only altered the region abutting the basepoint. Thus $\partial_{L}(T\# T') = \partial_{L_{1}}(T)\# T' + T \# \partial_{L_{2}} T'$ (where we extend $\#$ linearly). Consequently, the chain complex for $L$ is the tensor product of chain complexes for $L_{1}$ and $L_{2}$. Since we are working over a field, the result follows.   
$\Diamond$\\
\ \\
\noindent We now turn to computing the homology exactly for certain links, which can then be building blocks as above. Recall that a link $\lk{L}$ is called quasi-alternating if it is in the set $\mathcal{Q}$, the smallest set of links such that
\begin{itemize}
\item The unknot is in $\mathcal{Q}$;
\item If $\lk{L}$ has a diagram $L$ containing a crossing $c$ such that the two resolutions at $c$, $L_{0}$ and $L_{1}$ represent links $\lk{L}_{0}, \lk{L}_{1} \in \mathcal{Q}$ with
$det(\lk{L}) = det(\lk{L}_{0}) + det(\lk{L}_{1})$, then $L \in \mathcal{Q}$
\end{itemize}
Alternating links are quasi-alternating, and $det(\lk{L}) > 0$  when $\lk{L}$ is quasi-alternating. As  for knot homology theories, for quasi-alternating links, the homology is especially simple.

\begin{theorem}[]\ref{thm:quasi}
If $\lk{L}$ represents a quasi-alternating link with a connected diagram, then $\HT[i]{L} \cong 0$ when $i \neq \sigma(L)$ and has rank $\mathrm{det(L)}$ when $i = \sigma(L)$. 
\end{theorem}

\noindent{\bf Proof:} We first prove this result for alternating links. Pick a non-split, reduced diagram $L$. Checkerboard color the faces of $L$ so that $W_{L}$ a connected graph whose edges all correspond to $ic$-crossings. Suppose there are $n+1$ white faces. Thus, for each tree $T$ representing a generator for $\CT{L}$ we have $\degrad{T} = n$ since there are $(n+1) - 1$ edges in any spanning tree of $W_{L}$ and each edge in $T$ receives a $1$ resolution, while each edge not in $T$ receives a $0$ resolution. Consequently, every tree $T$ occurs in a single grading, and $\partial_L = 0$. Therefore, the homology is supported in only one grading, and has rank equal to the number of spanning trees for $W_{L}$. The number of spanning trees of $W_{L}$ is equal to $det(L)$. Furthermore, by the result of Gordon and Litherland, we may compute the signature of $\lk{L}$ by $\sigma(\lk{L}) = n - n_{+}(L) = \degrad{T} - n_{+}(L)$. This is the grading for $T$ in $\HT{L}$. Consequently, all the homology is in the grading given by $\sigma(\lk{L})$. \\ \ \\

\noindent For the quasi-alternating links, we will prove the result by induction on the number of crossings. This reproves the result for alternating links. First, note that the unknot has the stated property. The usual unknot diagram has a single resolution which is connected, and it has $\degrad{T} = 0$. Shifting by $-n_{+} = 0$ shows that the homology for the unknot is rank $1$, the determinant of the unknot, supported in grading $0$, the signature of the unknot. We  now proceed as in the analogous argument for reduced Khovanov homology found in \cite{MaOz}. Suppose that the result holds for all quasi-alternating links possessing a diagram with fewer than $n$ crossings, and that $\lk{L}$ is a quasi-alternating link such that there is a diagram $L$ as in the specification of $\mathcal{Q}$ with crossing $c$ such that resolving at $c$ results in quasi-alternating links. We use the following lemma, when $\lk{L}$ is an oriented link:

\begin{lemma}[Manolescu-Ozsv\'ath]
Suppose that $det(\lk{L}_{v})$, $det(\lk{L}_{h}) > 0$ and $det(\lk{L}_{+}) = det(\lk{L}_{v}) + det(\lk{L}_{h})$. Then 
$$
\sigma(\lk{L}_{v}) - \sigma(\lk{L}_{+}) = 1
$$
and
$$
\sigma(\lk{L}_{h}) - \sigma(\lk{L}_{+}) = e'
$$
where $e' = n_{-}(L_{h}) - n_{-}(L_{+})$.
\end{lemma}

\noindent Here $L_{+}$ is the diagram where the crossing $c$ is made to be positive, i.e. either $L$ or $L$ with the crossing $c$ switched, and $L_{h}$, $L_{v}$ are the $1$ and $0$ resolutions of $L_{+}$. Furthermore, $L_{h}$ is oriented arbitrarily. Therefore, $n_{+}(L_{v}) = n_{+}(L_{+}) - 1$ and $n_{+}(L_{h}) = n_{+}(L_{+}) + n_{-}(L_{+}) - 1 - n_{-}(L_{h}) = n_{+}(L_{+}) - e' - 1$. \\ \ \\

\noindent Now consider the long exact sequence in (\ref{eqn:posshiftedseq}). Notice that the $e'$ in the lemma equals the $e - 1$ in (\ref{eqn:posshiftedseq}).
If we let $i = \sigma(\lk{L}_{+}) + j$ then \\
\begin{equation}
\cdots \rightarrow \HTsub[\sigma(\lk{L}_{h}) + j]{L}{h}\otimes \F_{L} \rightarrow \HT[\sigma(\lk{L}) + j]{L} \rightarrow \HTsub[\sigma(\lk{L}_{v}) + j]{L}{v}\otimes \F_{L} \stackrel{\tau_{c,\ast}}{\longrightarrow} 
\HTsub[\sigma(\lk{L}_{h})+j+2]{L}{h}\otimes \F_{L} \rightarrow \cdots
\end{equation}
\ \\
\noindent Consequently, we may proceed by induction. If $\HTsub[\ast]{L}{h}$ is supported only in grading $\sigma(\lk{L}_{h})$ and $\HTsub[\ast]{L}{v}$ is supported only in $\sigma(\lk{L}_{v})$ then $\HT[\ast]{L}$ will only be supported in $\sigma(\lk{L})$. \\ \ \\

\noindent If the crossing $c$ of $L$ is negative, then we may apply the proceeding argument to $\mir{L}$ and use the symmetry under mirrors on $L_{0}$ and $L_{1}$ to conclude that $\HT{\mir{L}}$ is supported only in grading $\sigma(\mir{\lk{L}}) = -\sigma(\lk{L})$. Using the symmetry under mirrors again gives the desired result. $\Diamond$

\end{document}